\documentclass[a4,11pt]{article}
\usepackage{graphicx} 
\usepackage{float}
\usepackage{amsmath,amsthm}
\usepackage{amsfonts}
\usepackage{latexsym}
\usepackage{amssymb}
\usepackage{setspace}
\usepackage{comment}

\makeatletter
 
  \@addtoreset{equation}{section}
 \makeatother

\begin{document}
\title{Probabilistic Approach to Mean Field Games and 
Mean Field Type Control Problems with Multiple Populations
~\footnote{
Forthcoming in {\it Minimax Theory and its Applications}.
All the contents expressed in this research are solely those of the author and do not represent any views or 
opinions of any institutions. The author is not responsible or liable in any manner for any losses and/or damages caused by the use of any contents in this research.
}
}

\author{Masaaki Fujii\footnote{Quantitative Finance Course, Graduate School of Economics, The University of Tokyo. }
}
\date{ 
This version:  2 November, 2020
}
\maketitle



\newtheorem{definition}{Definition}[section]
\newtheorem{assumption}{Assumption}[section]
\newtheorem{condition}{$[$ C}
\newtheorem{lemma}{Lemma}[section]
\newtheorem{proposition}{Proposition}[section]
\newtheorem{theorem}{Theorem}[section]
\newtheorem{remark}{Remark}[section]
\newtheorem{example}{Example}[section]
\newtheorem{corollary}{Corollary}[section]
\def\n{{\bf n}}
\def\A{{\bf A}}
\def\B{{\bf B}}
\def\C{{\bf C}}
\def\D{{\bf D}}
\def\E{{\bf E}}
\def\F{{\bf F}}
\def\G{{\bf G}}
\def\H{{\bf H}}
\def\I{{\bf I}}
\def\J{{\bf J}}
\def\K{{\bf K}}
\def\L{{\bf L}}
\def\M{{\bf M}}
\def\N{{\bf N}}
\def\O{{\bf O}}
\def\P{{\bf P}}
\def\Q{{\bf Q}}
\def\R{{\bf R}}
\def\S{{\bf S}}
\def\T{{\bf T}}
\def\U{{\bf U}}
\def\V{{\bf V}}
\def\W{{\bf W}}
\def\X{{\bf X}}
\def\Y{{\bf Y}}
\def\Z{{\bf Z}}
\def\cala{{\cal A}}
\def\calb{{\cal B}}
\def\calc{{\cal C}}
\def\cald{{\cal D}}
\def\cale{{\cal E}}
\def\calf{{\cal F}}
\def\calg{{\cal G}}
\def\calh{{\cal H}}
\def\cali{{\cal I}}
\def\calj{{\cal J}}
\def\calk{{\cal K}}
\def\call{{\cal L}}
\def\calm{{\cal M}}
\def\caln{{\cal N}}
\def\calo{{\cal O}}
\def\calp{{\cal P}}
\def\calq{{\cal Q}}
\def\calr{{\cal R}}
\def\cals{{\cal S}}
\def\calt{{\cal T}}
\def\calu{{\cal U}}
\def\calv{{\cal V}}
\def\calw{{\cal W}}
\def\calx{{\cal X}}
\def\caly{{\cal Y}}
\def\calz{{\cal Z}}
%
\def\sskip{\hspace{0.5cm}}
\def\simleq{ \raisebox{-.7ex}{\em $\stackrel{{\textstyle <}}{\sim}$} }
\def\leqsim{ \raisebox{-.7ex}{\em $\stackrel{{\textstyle <}}{\sim}$} }
\def\ep{\epsilon}
\def\half{\frac{1}{2}}
\def\iku{\rightarrow}
\def\Iku{\Rightarrow}
\def\ikup{\rightarrow^{p}}
\def\inclusion{\hookrightarrow}
\def\cadlag{c\`adl\`ag\ }
\def\up{\uparrow}
\def\down{\downarrow}
\def\doti{\Leftrightarrow}
\def\douti{\Leftrightarrow}
\def\dochi{\Leftrightarrow}
\def\douchi{\Leftrightarrow}%
\def\yy{\\ && \nonumber \\}
\def\y{\vspace*{3mm}\\}
\def\nn{\nonumber}
\def\be{\begin{equation}}
\def\ee{\end{equation}}
\def\bea{\begin{eqnarray}}
\def\eea{\end{eqnarray}}
\def\beas{\begin{eqnarray*}}
\def\eeas{\end{eqnarray*}}
%
\def\hd{\hat{D}}
\def\hv{\hat{V}}
\def\hsd{{\hat{d}}}
\def\hx{\hat{X}}
\def\hsx{\hat{x}}
\def\bsx{\bar{x}}
\def\bsd{{\bar{d}}}
\def\bx{\bar{X}}
\def\ba{\bar{A}}
\def\bb{\bar{B}}
\def\bc{\bar{C}}
\def\bv{\bar{V}}
\def\balpha{\bar{\alpha}}
\def\bbalpha{\bar{\bar{\alpha}}}
\def\combi{\l(\begin{array}{c}\alpha\\ \beta \end{array}\r)}
\def\f{^{(1)}}
\def\s{^{(2)}}
\def\ss{^{(2)*}}
\def\l{\left}
\def\r{\right}
\def\a{\alpha}
\def\b{\beta}
\def\L{\Lambda}


\def\E{{\bf E}}
\def\P{{\bf P}}
\def\Q{{\bf Q}}
\def\R{{\bf R}}

\def\cadlag{{c\`adl\`ag~}}

\def\calf{{\cal F}}
\def\calp{{\cal P}}
\def\calq{{\cal Q}}
\def\wtW{\widetilde{W}}
\def\wtB{\widetilde{B}}
\def\wtPsi{\widetilde{\Psi}}
\def\wt{\widetilde}
\def\mbb{\mathbb}
\def\ol{\overline}
\def\ul{\underline}

\def\hTheta{\hat{\Theta}}
\def\hPhi{\hat{\Phi}}
\def\L2nu{{\mbb{L}^2(\nu)}}
\def\esup{{\rm ess}\sup}

\def\LDis{\frac{\bigl.}{\bigr.}}
\def\ep{\epsilon}
\def\del{\delta}
\def\Del{\Delta}
\def\part{\partial}
\def\wh{\widehat}
\def\vep{\varepsilon}
\def\ul{\underline}
\def\ol{\overline}

\def\nn{\nonumber}
\def\be{\begin{equation}}
\def\ee{\end{equation}}
\def\bea{\begin{eqnarray}}
\def\eea{\end{eqnarray}}
\def\beas{\begin{eqnarray*}}
\def\eeas{\end{eqnarray*}}
\def\tbf{\textbf}
\def\bg{\boldsymbol}

\def\bull{$\bullet~$}

\newcommand{\Slash}[1]{{\ooalign{\hfil/\hfil\crcr$#1$}}}

\begin{abstract}
In this work, we systematically investigate mean field games and mean field type control problems with multiple populations. We study the mean field limits of the three different situations; (i) every agent is non-cooperative; (ii) the agents within each population are cooperative; and (iii) the agents in some populations are cooperative. We provide several sets of sufficient conditions for the existence of a mean field equilibrium for each case. We also show that, under appropriate conditions,  each mean field solution actually provides an approximate Nash equilibrium for the corresponding game with a large but finite number of agents.
\end{abstract}
\vspace{5mm}

\section{Introduction}
In pioneering works of Lasry \& Lions~\cite{Lions-1, Lions-2, Lions-3} and  Huang, Malhame \& Caines~\cite{Caines-Huang}, 
the two groups of researchers independently proposed a powerful technique to produce an approximate Nash equilibrium for
stochastic differential games among a large number of agents with symmetric interactions.
Importantly, each agent is assumed to be affected by the states of the other agents only through 
their empirical distribution.
In the large population limit, the problem is shown to result in two highly coupled nonlinear partial 
differential equations (PDEs), the one is of the Hamilton-Jacobi-Bellman type, which takes 
care of the optimization problem, while the other is of the Kolmogorov type guaranteeing  the consistent time evolution 
of the distribution of the individual states of the agents.
The greatest benefit of the mean-field game approach is to render notoriously intractable problems of stochastic differential games among 
many agents into simpler stochastic optimal control problems.
For details of the analytical approach and its various applications, one may consult  the monographs by 
Bensoussan, Frehse \& Yam~\cite{Bensoussan-mono}, Gomes, Pimentel \& Voskanyan~\cite{Gomes-reg} and also Kolokoltsov \& Malafeyev~\cite{Kolokoltsov-mono}.

In a series of works \cite{Carmona-Delarue-MKV, Carmona-Delarue-MFG, Carmona-Delarue-MFTC}, 
Carmona \& Delarue developed a probabilistic approach to these problems, 
where  forward-backward stochastic differential equations (FBSDEs) of McKean-Vlasov type instead of PDEs were shown to be 
the relevant objects for investigation.
In particular, they provided the sufficient conditions for the existence of an equilibrium for mean field games with 
the cost functions of quadratic growth in \cite{Carmona-Delarue-MFG}.
In the case of cooperative agents who adopt the common feedback control function, they showed in \cite{Carmona-Delarue-MFTC} that
the large population limit results in the optimization problem with respect to a controlled McKean-Vlasov SDE.
Using the notion of so-called L-derivative, which is a type of derivatives for functions defined on the space of probability measures,
they solved the problem by a new class of FBSDEs of McKean-Vlasov type.
A probabilistic but weak formulation of  mean-field games based on the concept of relaxed control 
has been studied in \cite{Carmona-Lacker, Carmona-Delarue-Lacker, Djete,Lacker-1, Lacker-2, Lacker-3},
which does not produce any equation characterizing the equilibrium solution but can significantly weaken
the regularity assumptions we need. 
In particular, we refer to Carmona, Delarue \& Lacker~\cite{Carmona-Delarue-Lacker} for problems in the presence of common noise.
For interested readers for the probabilistic approach, the two volumes of excellent monograph~\cite{Carmona-Delarue-1, Carmona-Delarue-2}
are now available. 

The above explained developments of the mean field game theory have opened a new horizon for 
variety of interesting multi-agent problems.  In particular, there appeared  large amount of literature
on, optimal trading, optimal liquidation, optimal exploitation of exhaustible resources, price formation 
in an electricity market, zero-sum game, systemic risk analysis of interbank network, etc.
See, for example, \cite{Matoussi, Carmona-Fouque, Carmona-Tan,Casgrain, Djehiche-E, Evangelista, Feron,  Fu-Horst-1, Fu-Horst-2, Fu,  Gomes-Mohr, Gueant-Oil, Lehalle, Noce, Parise, Sun}. A quite popular approach adopted in many of these references 
is to assume that the relevant price process is decomposed into two parts, one is 
an exogenous process which is independent of the agents' action, and the other representing the price impact which is often proportional to the average trading speed of the agents. 
Various examples of interesting applications to macroeconomic problems can be found in \cite{Lions-Eco, Gomes-eco},
where the economic interactions are studied in the continuum limit of agents. The mean field game theory with a major player~\cite{Bensoussan-Chau, Bensoussan-Chau-2, Bensoussan-Cass, Carmona-Zhu, Huang-Jaimungal, Nguyen} is another important branch of research. Although the impact from an individual minor player becomes insignificant as the size of population grows,  the major player does not lose its impact to the system.
The literature provides an important framework of mean field games for the applications to the financial as well as economic problems 
in the presence of a giant bank, a central bank,  or a monopolistic producer of certain goods or services, etc.

In this work, we are interested in  mean field games and mean field type control problems 
in the presence of multiple populations.
Here, the same cost functions
as well as the coefficient functions of the state dynamics are shared among the agents {\it within} each population, 
but they can be different population by population.
Mean field games with multiple populations arise naturally in most of the practical applications, 
and have been studied by many researchers following the analytic approach.  
In fact, a related problem was already considered in the first original work of \cite{Caines-Huang}.
Lachapelle \& Wolfram~\cite{Wolfram} modeled a congestion problem of pedestrian crowds, 
and Achdou, Bardi \& Cirant~\cite{Achdou-segregation} studied the issue of urban settlements
and residential choice using the mean-field game representation.
Feleqi~\cite{Feleqi} and Cirant~\cite{Cirant} dealt with ergodic mean field games of multiple populations under different boundary conditions.
Among these works, the recent publication by Bensoussan, Huang \& Lauriere \cite{Bensoussan-mfg-paper}
is most closely related to the current work.
They systematically studied the problems of mean-field games and mean-field type control problems
with multiple populations. In particular,  their analysis includes 
the case where the agents within each population are cooperative but compete with those in the other populations.

In the current paper, differently from the existing works, we have adopted the probabilistic approach and
closely followed the procedures developed in \cite{Carmona-Delarue-MFG, Carmona-Delarue-MFTC}. 
In addition to the mean field games of multiple populations, we have studied the situation where the 
agents in each population are cooperative as in \cite{Bensoussan-mfg-paper}, and yet another situation 
which is a mixture of the first two cases: the agents in some populations are cooperative within their own
but those in the other populations are not.
The presence of multiple populations induces a system of FBSDEs of McKean-Vlasov type.
Although it is a {\it coupled} system of FBSDEs due to the interactions among different populations, the couplings appear only through 
the mean field interactions i.e., the distribution of the state of the representative agent of each population.
This feature allows us to solve a matching problem corresponding to the state of equilibrium
by Schauder's fixed point theorem in a quite similar manner to \cite{Carmona-Delarue-MFG}. 
In \cite{Bensoussan-mfg-paper}, the mean field equilibrium is assumed to exist.
Combined with suitable convexity and differentiability,  the authors derived the coupled system of PDEs
as the {\it necessary condition} for optimality.  On the other hand, although we need more technical assumptions than in \cite{Bensoussan-mfg-paper}, 
in each of the three cases mentioned above, we have found several sets of {\it sufficient (instead of necessary) conditions} 
for the existence of the mean field equilibrium, in particular the one which allows the cost functions of quadratic 
growth both in the state variable as well as in its distribution so that
it is applicable to some of the popular linear quadratic problems.
Another important advantage in adopting the probabilistic approach is that we can discuss quantitatively the 
relation between the mean-field limits and the corresponding games with finite population size.
In \cite{Bensoussan-mfg-paper}, as a common feature of the analytic approach,  the problems are discussed only in the continuum limit of 
agents. On the other hand, thanks to the powerful technique called the {\it propagation-of-chaos}, we have proved that
each mean field solution provides an approximate Nash equilibrium for the corresponding game
with finite number of agents. By Glivenko-Cantelli convergence theorem in the Wasserstein distance,
we have actually obtained  a non-asymptotic estimate in terms of the population size $N$.
Our analysis also highlights an interesting difference between the game in which  all the 
agents are non-cooperative and the one in which  the agents are cooperative in some populations.
In the non-cooperative case, the impact from each agent dissipates in the large population limit. 
On the other hand, when the agents are cooperative, the aggregate impact from their
common strategy does not  dissipates in the same limit.  
We shall observe that this feature in the latter requires us to make more stringent assumptions on the coefficients functions
so that the mean-field limit gives an approximate Nash equilibrium for the game of finite population.
This additional difficulty looks somewhat similar to the situation 
for the mean field games with a major player. In fact, in general,  the optimization problem for the major player is known to become a McKean-Vlasov type~\cite{Carmona-Zhu}. Carrying out more detailed comparison as well as  finding a general characterization of equilibrium
to allow more flexible interactions in the coefficients remain interesting research topics for the future.

Finally, let us make some comment on the issues related to the common noise.
In the current paper, we study the problems involving only the idiosyncratic noises for simplicity.
However, if the common noise has only finite number of states, the equilibrium can be constructed 
as a simple superposition of equilibria, each of which corresponds to the equilibrium for 
a given sate of the common noise (See Chapter 3 in \cite{Carmona-Delarue-2}). 
Hence, after approximating the common noise by a finite discretization,
the problem can be handled by the same method provided in this work. 
In general, taking the limit of finer discretization requires  delicate arguments of weak convergence \cite{Carmona-Delarue-2}.  
As a related but different approach based on the relaxed control, we refer to the work~\cite{Carmona-Delarue-Lacker}.
On the other hand, we can directly apply Peng-Wu's continuation method~\cite{Peng-Wu} 
to obtain a unique strong solution for mean field games even in the presence of  common noise,
if the coefficients functions are Lipschitz continuous also in the measure arguments
and additionally satisfy the appropriate monotone conditions. 
See \cite{Fujii-Takahashi-1, Fujii-Takahashi-2} for applications to a market clearing equilibrium in the presence of common noise and multiple populations. 

The organization of the paper is as follows: after explaining notation in Section~\ref{sec-notation},
we study the mean field problems in the first half of the paper; in Section~\ref{sec-mfg} (i) the case of non-cooperative agents,
in Section~\ref{sec-mftc} (ii) the case where the agents are cooperative within each population,
and in Section~\ref{sec-mftc-mfg} (iii) the agents in some populations are cooperative but those in the other populations are not.
In the second half of the paper, we investigate the corresponding problem with a finite number of agents;
we treats in Section~\ref{sec-mfg-fa}  the case (i), 
in Section~\ref{sec-mftc-fa}  the case (ii),  and finally in Section~\ref{sec-mftc-mfg-fa} we treats the case (iii).
Although we set the number of populations to two in the main analysis, this is just for notational convenience.
We shall see that the analysis can be easily generalized to any finite number of populations.
Finally, we conclude in Section~\ref{sec-conclusion}.
\section{Notations}
\label{sec-notation}
Throughout the paper, we work on some complete probability space  $(\Omega,\calf,\mbb{P})$
equipped with a right-continuous and complete filtration $\mbb{F}=(\calf_t)_{t\in[0,T]}$
supporting two independent $d$-dimensional
standard Brownian motions $\bigl(\bg{W}^1=(W^1_t)_{t\in[0,T]}, \bg{W}^2=(W^2_t)_{t\in[0,T]}\bigr)$
as well as two independent random variables $\xi^1, \xi^2 \in \mbb{L}^2(\Omega,\calf_0,\mbb{P};\mbb{R}^d)$.
For each $i\in\{1,2\}$, $\mbb{F}^i:=(\calf_t^i)_{t\in[0,T]}\subset \mbb{F}$ is a complete and right-continuous 
filtration generated by $(\xi^i, \bg{W}^i)$. Here, $T>0$ is a given terminal time.
To lighten the notation, unless otherwise stated, we use indices $i$ and $j$ specifically to represent an element in $\{1,2\}$,  and 
we always suppose that $j\neq i$ when they appear in the same expression. Moreover, 
we use the symbol $C$ to represent a general nonnegative constant which may change line by line.
When we want to emphasize that $C$ depends only on some specific variables, say $a$ and $b$, we 
use the symbol $C(a,b)$. We let $||\cdot||_2$ denote the $\mbb{L}^2(\Omega,\calf,\mbb{P};\mbb{R}^d)$-norm.
We use the following notations for frequently encountered spaces:\\
$\bullet~\mbb{S}^2$ is the set of $\mbb{R}^d$-valued continuous processes $\bg{X}$ satisfying
\be
||X||_{\mbb{S}^2}:=\mbb{E}\bigl[\sup_{t\in[0,T]}|X_t|^2\bigr]^\frac{1}{2}<\infty~. \nn 
\ee
$\bullet~\mbb{S}^\infty$ is the set of $\mbb{R}^d$-valued essentially bounded continuous processes $\bg{X}$
satisfying
\be
||X||_{\mbb{S}^\infty}:=\bigl|\bigl| \sup_{t\in[0,T]} |X_t| \bigr|\bigr|_{\infty}<\infty~. \nn
\ee
$\bullet~\mbb{H}^2$ is the set of $\mbb{R}^{d\times d}$-valued progressively measurable processes $Z$ satisfying
\be
||Z||_{\mbb{H}^2}:=\mbb{E}\Bigl[\Bigl(\int_0^T |Z_t|^2dt\Bigr)\Bigr]^\frac{1}{2}<\infty~.\nn
\ee
$\bullet~\call(X)$ denotes the law of a random variable $X$. \\
$\bullet~\calm_f^1(\mbb{R}^d)$ is the set of finite signed measures $\mu$ on $(\mbb{R}^d,\calb(\mbb{R}^d))$ such that
$\int_{\mbb{R}^d} |x| d|\mu|(x)<\infty$.  \\
$\bullet~\calp(\mbb{R}^d)$ is the set of probability measures on $(\mbb{R}^d,\calb(\mbb{R}^d))$. \\
$\bullet~\calp_p(\mbb{R}^d)$ with $p\geq 1$ is the subset of $\calp(\mbb{R}^d)$ with finite $p$-th moment; i.e.,
the set of $\mu\in \calp(\mbb{R}^d)$ satisfying
\bea
M_p(\mu):=\Bigl(\int_{\mbb{R}^d}|x|^p \mu(dx)\Bigr)^\frac{1}{p}<\infty~.\nn
\eea
We always assign $\calp_p(\mbb{R}^d)$ with $(p\geq 1)$ the $p$-Wasserstein distance $W_p$,
which makes the space $\calp_p(\mbb{R}^d)$ a complete separable metric space.
As an important property, for any $\mu, \nu\in \calp_p(\mbb{R}^d)$,  we have
\bea
W_p(\mu,\nu)={\inf}\Bigl\{ \mbb{E}[|X-Y|^p]^\frac{1}{p}; \call(X)=\mu, \call(Y)=\nu\Bigr\}~. \nn
\eea
For more details,  see Chapter 5 in \cite{Carmona-Delarue-1}
or Chapter 3 in \cite{Bogachev}.

\section{Mean Field Games with Multiple Populations}
\label{sec-mfg}
In this section, we consider a mean-field limit of  a game among a large number of non-cooperative agents
in the presence of two populations. Here, each agent competes with 
all the other agents but shares the common cost functions as well as coefficient functions of the state dynamics
within each population. As we shall see,  extending to the general situation with a finite number of populations is straightforward.
Corresponding problem with a finite number of agents and its relation to the mean-field problem discussed in this section
will be investigated in Section~\ref{sec-mfg-fa}.
\subsection{Definition of the Mean Field Problem}
Before specifying detailed assumptions, let us formulate the problem of finding 
an equilibrium in the limiting framework. It proceeds in the following three steps.\\
(i) Fix any two deterministic flows of probability measures $(\bg{\mu}^i=(\mu^i_t)_{t\in[0,T]})_{i\in\{1,2\}}$ 
given on $\mbb{R}^d$.\\
(ii) Solve the two optimal control problems
\be
\inf_{\bg{\alpha}^1\in \mbb{A}_1}J_1^{\bg{\mu}^1,\bg{\mu}^2}(\bg{\alpha}^1), \qquad 
\inf_{\bg{\alpha}^2\in \mbb{A}_2}J_2^{\bg{\mu}^2,\bg{\mu}^1}(\bg{\alpha}^2) 
\label{mfg-optimal}
\ee
over some  admissible strategies $\mbb{A}_i$ $(i\in\{1,2\})$, where 
\bea
&&J_1^{\bg{\mu}^1,\bg{\mu}^2}(\bg{\alpha}^1):=\mbb{E}\Bigl[\int_0^T f_1(t,X_t^1,\mu_t^1,\mu_t^2,\alpha_t^1)dt+
g_1(X_T^1,\mu_T^1,\mu_T^2)\Bigr]~,\nn \\
&&J_2^{\bg{\mu}^2,\bg{\mu}^1}(\bg{\alpha}^2):=\mbb{E}\Bigl[\int_0^T f_2(t,X_t^2,\mu_t^2,\mu_t^1,\alpha_t^2)dt+
g_2(X_T^2,\mu_T^2,\mu_T^1)\Bigr]~,\nn
\eea
subject to the $d$-dimensional diffusion dynamics:
\bea
&&dX_t^1=b_1(t,X_t^1,\mu_t^1,\mu_t^2,\alpha_t^1)dt+\sigma_1(t,X_t^1,\mu_t^1,\mu_t^2)dW_t^1~,\nn \\
&&dX_t^2=b_2(t,X_t^2,\mu_t^2,\mu_t^1,\alpha_t^2)dt+\sigma_2(t,X_t^2,\mu_t^2,\mu_t^1)dW_t^2~,\nn
\eea
for $t\in[0,T]$ with $\bigl(X_0^i=\xi^i\in \mbb{L}^2(\Omega,\calf_0,\mbb{P};\mbb{R}^d)\bigr)_{1\leq i\leq 2}$. 
For each population $i\in\{1,2\}$, we suppose that $\mbb{A}_i$ is the set of  $A_i$-valued $\mbb{F}^i$-progressively
measurable processes $\bg{\alpha}^i$ satisfying
$\mbb{E}\int_0^T |\alpha^i_t|^2 dt<\infty$
where $A_i\subset \mbb{R}^k $ is closed and convex. \\
(iii) Find a pair of probability flows $(\bg{\mu}^1, \bg{\mu}^2)$ as a solution to the matching problem:
\bea
\forall t\in[0,T],\qquad \mu_t^1=\call(\hat{X}_t^{1,\bg{\mu}^1,\bg{\mu}^2}), \quad
\mu_t^2=\call(\hat{X}_t^{2,\bg{\mu}^2,\bg{\mu}^1})~,
\label{mfg-fixed-point}
\eea
where $(\hat{X}^{i,\bg{\mu}^i,\bg{\mu}^j})_{i\in\{1,2\}, j\neq i}$ are the solutions to the optimal control problems in  (ii).

\begin{remark}
It is just for convenience to use the common dimension $d$ (as well as $k$ for $A_i$)
for both populations.  Note also that since $\mu^i$ is deterministic and $\alpha^i$ is $\mbb{F}^i$-adapted,
$X^1$ and $X^2$ are independent. Hence, there is no gain of information by considering the joint law $\call(X^1,X^2)$.

\end{remark}

\subsection{Optimization for  given flows of probability measures}
The main assumptions in this section are as follows:
\begin{assumption}{\rm{\tbf{(MFG-a)}}} 
$L, K\geq 0$ and $\lambda>0$ are some constants. For $1\leq i\leq 2$, the measurable functions $b_i:[0,T]\times \mbb{R}^d\times \calp_2(\mbb{R}^d)^2\times A_i\rightarrow \mbb{R}^d$,
$\sigma_i:[0,T]\times \mbb{R}^d\times \calp_2(\mbb{R}^d)^2\rightarrow \mbb{R}^{d\times d}$,
$f_i:[0,T]\times \mbb{R}^d\times \calp_2(\mbb{R}^d)^2\times A_i\rightarrow \mbb{R}$,
and $g_i:\mbb{R}^d\times \calp_2(\mbb{R}^d)^2\rightarrow \mbb{R}$ satisfy the following conditions: \\
\tbf{(A1)} The functions $b_i$ and $\sigma_i$ are affine in $(x,\alpha)$ in the sense that, for any
$(t,x,\mu,\nu,\alpha)\in[0,T]\times \mbb{R}^d\times \calp_2(\mbb{R}^d)^2\times A_i$,
\bea
&&b_i(t,x,\mu,\nu,\alpha):=b_{i,0}(t,\mu,\nu)+b_{i,1}(t,\mu,\nu)x+b_{i,2}(t,\mu,\nu)\alpha~, \nn \\
&&\sigma_i(t,x,\mu,\nu):=\sigma_{i,0}(t,\mu,\nu)+\sigma_{i,1}(t,\mu,\nu)x~,\nn
\eea
where $b_{i,0}, b_{i,1}, b_{i,2}, \sigma_{i,0}$ and $\sigma_{i,1}$ defined on $[0,T]\times \calp_2(\mbb{R}^d)^2$
are $\mbb{R}^d$, $\mbb{R}^{d\times d}, \mbb{R}^{d\times k}, \mbb{R}^{d\times d}$ and $\mbb{R}^{d\times d\times d}$-valued
measurable functions, respectively. \\
\tbf{(A2)} For any $t\in[0,T]$, the functions $\calp_2(\mbb{R}^d)^2\ni (\mu,\nu)
\mapsto (b_{i,0},b_{i,1},b_{i,2}, \sigma_{i,0}, \sigma_{i,1})(t,\mu,\nu)$ are continuous in $W_2$-distance.
Moreover for any $(t,\mu,\nu)\in [0,T]\times \calp_2(\mbb{R}^d)^2$, 
\bea
&&|b_{i,0}(t,\mu,\nu)|, |\sigma_{i,0}(t,\mu,\nu)|\leq K+L\bigl(M_2(\mu)+M_2(\nu)\bigr)~, \nn \\
&&|b_{i,1}(t,\mu,\nu)|, |b_{i,2}(t,\mu,\nu)|, |\sigma_{i,1}(t,\mu,\nu)|\leq L~. \nn
\eea
\tbf{(A3)} The function $\mbb{R}^d\times A_i\ni (x,\alpha)\mapsto f_i(t,x,\mu,\nu,\alpha)\in \mbb{R}$ is once 
continuously differentiable with $L$-Lipschitz derivatives, i.e. for any $t\in[0,T], \mu,\nu\in \calp_2(\mbb{R}^d)$,
 $x,x^\prime\in \mbb{R}^d$, $\alpha,\alpha^\prime \in A_i$, 
\be
|\part_{(x,\alpha)}f_i(t,x^\prime ,\mu,\nu,\alpha^\prime )-\part_{(x,\alpha)}f_i(t,x,\mu,\nu,\alpha)|
\leq L\bigl(|x^\prime-x|+|\alpha^\prime-\alpha|\bigr)~,\nn
\ee 
where $\part_{(x,\alpha)}f_i$ denotes the gradient in the joint variables $(x,\alpha)$. $f_i$ also 
satisfies the $\lambda$-convexity:
\bea
f_i(t,x^\prime,\mu, \nu, \alpha^\prime)-f_i(t,x,\mu,\nu,\alpha)-\langle (x^\prime-x,\alpha^\prime-\alpha), \part_{(x,\alpha)}
f_i(t,x,\mu,\nu,\alpha)\rangle \geq \lambda|\alpha^\prime-\alpha|^2~.\nn 
\eea
\tbf{(A4)}
For any $(t,x,\mu,\nu,\alpha)\in [0,T]\times \mbb{R}^d\times\calp_2(\mbb{R}^d)^2
\times A_i$, 
\bea
&&|\part_{(x,\alpha)}f_i(t,x,\mu,\nu,\alpha)|\leq K+L\bigl(|x|+|\alpha|+M_2(\mu)+M_2(\nu)\bigr)~.\nn\\
&&|f_i(t,x,\mu,\nu,\alpha)|\leq K+L\bigl(|x|^2+|\alpha|^2+M_2(\mu)^2+M_2(\nu)^2\bigr)~.\nn
\eea
Moreover, for any $(t,x,\alpha)\in[0,T]\times \mbb{R}^d\times A_i$, the functions $\calp_2(\mbb{R}^d)^2\ni (\mu,\nu) \mapsto 
f_i (t,x,\mu,\nu,\alpha)$ and  $\calp_2(\mbb{R}^d)^2\ni (\mu,\nu) \mapsto 
\part_{(x,\alpha)}f_i (t,x,\mu,\nu,\alpha)$ are continuous in $W_2$-distance. \\
\tbf{(A5)} For any $\mu,\nu\in \calp_2(\mbb{R}^d)$, the function $\mbb{R}^d\ni x\mapsto g_i(x,\mu,\nu)\in \mbb{R}$ 
is convex.  It is also once continuously differentiable with $L$-Lipschitz derivatives, i.e. $\forall x,x^\prime \in \mbb{R}^d, 
\mu,\nu\in\calp_2(\mbb{R}^d)$, it holds
\be
|\part_x g_i(x^\prime,\mu,\nu)-\part_x g_i(x,\mu,\nu)|\leq L|x^\prime-x|~.\nn
\ee
For any $x\in \mbb{R}^d$, the functions $\calp_2(\mbb{R}^d)^2\ni (\mu,\nu)\mapsto g_i(x,\mu,\nu)$
and $\calp_2(\mbb{R}^d)^2\ni (\mu,\nu)\mapsto \part_x g_i(x,\mu,\nu)$ are
continuous in $W_2$-distance. Moreover, the growth conditions
\bea
&&|\part_x g_i(x,\mu,\nu)|\leq K+L\bigl(|x|+M_2(\mu)+M_2(\nu)\bigr)~, \nn \\
&&|g_i(x,\mu,\nu)|\leq K+L\bigl(|x|^2+M_2(\mu)^2+M_2(\nu)^2\bigr)~, \nn
\eea
are satisfied.
\end{assumption}

We first consider the optimal control problem $(\ref{mfg-optimal})$ for given 
deterministic flows of probability measures. The corresponding Hamiltonian for each population 
$H_i:[0,T]\times \mbb{R}^d\times \calp_2(\mbb{R}^d)^2\times \mbb{R}^{d}\times \mbb{R}^{d\times d}\times  
A_i\rightarrow \mbb{R}$ is defined by:
\bea
H_i(t,x,\mu,\nu,y,z,\alpha):=\langle b_i(t,x,\mu,\nu,\alpha),y\rangle+{\rm tr}[\sigma_i(t,x,\mu,\nu)^\top z]
+f_i(t,x,\mu,\nu,\alpha)~.
\label{def-H}
\eea
Since $\sigma_i$ is independent of the control parameter, the minimizer $\hat{\alpha}_i(t,x,\mu,\nu,y)$ of the Hamiltonian 
$H_i$ can also be defined by a reduced Hamiltonian $H_i^{(r)}$:
\bea
\hat{\alpha}_i(t,x,\mu,\nu,y):={\rm{argmin}}_{\alpha\in A_i} H_i^{(r)}(t,x,\mu,\nu,y, \alpha) 
\label{def-alpha}
\eea
where
\be
H_i^{(r)}(t,x,\mu,\nu,y,\alpha):=\langle b_i(t,x,\mu,\nu,\alpha),y\rangle
+f_i(t,x,\mu,\nu,\alpha).~\nn
\ee
The following result regarding the regularity of $\hat{\alpha}_i$ is a straightforward extension of 
Lemma 2.1~\cite{Carmona-Delarue-MFG}.
\begin{lemma}
\label{lemma-alpha}
Under Assumption {\rm{\tbf{(MFG-a)}}}, for all $(t,x,\mu,\nu,y)\in[0,T]\times \mbb{R}^d\times \calp_2(\mbb{R}^d)^2\times \mbb{R}^d$,
there exists a unique minimizer $\hat{\alpha}_i(t,x,\mu,\nu,y)$ of $H_i^{(r)}$,
where the map $[0,T]\times \mbb{R}^d\times \calp_2(\mbb{R}^d)^2\times \mbb{R}^d\ni (t,x,\mu,\nu,y)\mapsto \hat{\alpha}_i(t,x,\mu,\nu,y)\in A_i$
is measurable. There exist constants $C$ depending only on $(L,\lambda)$ and $C^\prime$ depending additionally on $K$ such that,
for any $t\in[0,T], x,x^\prime, y, y^\prime \in \mbb{R}^d, \mu,\nu \in\calp_2(\mbb{R}^d)$, 
\bea
&&|\hat{\alpha}_i(t,x,\mu,\nu,y)|\leq C^\prime+C\bigl(|x|+|y|+M_2(\mu)+M_2(\nu)\bigr)\nn \\
&&|\hat{\alpha}_i(t,x,\mu,\nu,y)-\hat{\alpha}_i(t,x^\prime,\mu,\nu,y^\prime)|\leq C\bigl(|x-x^\prime|+|y-y^\prime|\bigr)~. \nn 
\eea
Moreover, for any $(t,x,y)\in[0,T]\times \mbb{R}^d\times \mbb{R}^d$, the map
$\calp_2(\mbb{R}^d)^2\ni (\mu, \nu)\mapsto \hat{\alpha}_i(t,x,\mu,\nu,y)$ is continuous with respect to $W_2$-distance:
\bea 
&&|\hat{\alpha}_i(t,x,\mu,\nu,y)-\hat{\alpha}_i(t,x,\mu^\prime, \nu^\prime, y)| \nn \\
&&\leq (2\lambda)^{-1}\Bigl(|b_{i,2}(t,\mu,\nu)-b_{i,2}(t,\mu^\prime, \nu^\prime)||y|+|\part_\alpha f_i(t,x,\mu,\nu,\hat{\alpha}_i)
-\part_\alpha f_i(t,x,\mu^\prime, \nu^\prime,\hat{\alpha}_i)|\Bigr)\nn
\eea
where $\hat{\alpha}_i:=\hat{\alpha}_i(t,x,\mu,\nu,y)$.
\begin{proof}
To lighten the notation, let us write $\rho=(\mu, \nu)\in \calp_2(\mbb{R}^d)^2$.
Since the function $A_i\ni \alpha \mapsto H_i^{(r)}(t,x,\rho,y,\alpha)$ is 
strictly convex and once continuously differentiable, $\hat{\alpha}_i(t,x,\rho,y)$ is given by 
the unique solution to the variational inequality:
\bea
\forall \beta\in A_i, \quad \langle \beta-\hat{\alpha}_i(t,x,\rho,y), \part_\alpha H_i^{(r)}(t,x,\rho,y,\hat{\alpha}_i(t,x,\rho,y))\rangle \geq 0~.
\label{alpha-vi}
\eea 
By strict convexity, the measurability is a consequence of the gradient descent algorithm (Lemma 3.3~\cite{Carmona-Delarue-1}).

With an arbitrary point $\beta_i\in A_i$, the $\lambda$-convexity implies that 
\bea
H_i^{(r)}(t,x,\rho,y,\beta_i)&\geq& H_i^{(r)}(t,x,\rho,y,\hat{\alpha}_i) \nn \\
&\geq &H_i^{(r)}(t,x,\rho,y,\beta_i)+\langle \hat{\alpha}_i-\beta_i, \part_\alpha H_i^{(r)}(t,x,\rho,y,\beta_i)\rangle+\lambda|\hat{\alpha}_i-\beta_i|^2~,\nn
\eea
where $\hat{\alpha}_i:=\hat{\alpha}_i(t,x,\rho,y)$. Hence we have
\be 
|\hat{\alpha}_i-\beta_i|\leq \lambda^{-1}\bigl(|b_{i,2}(t,\rho)||y|+|\part_\alpha f_i(t,x,\rho,\beta_i)|\bigr)~.
\label{alpha-growth}
\ee
This gives the first growth condition.
Next, with $\hat{\alpha}_i:=\hat{\alpha}(t,x,\rho,y)$ and $\hat{\alpha}_i^\prime:=\hat{\alpha}_t(t,x^\prime,\rho,y^\prime)$,
the optimality condition implies
\bea
\langle \hat{\alpha}_i^\prime-\hat{\alpha}_i, \part_\alpha H_i^{(r)}(t,x,\rho,y,\hat{\alpha}_i)-\part_\alpha H_i^{(r)}(t,x^\prime,\rho,y^\prime,
\hat{\alpha}_i^\prime)\rangle\geq 0~. \nn
\eea
This inequality, together with the $\lambda$-convexity, gives
\bea
\langle \hat{\alpha}_i^\prime-\hat{\alpha}_i,~ b_{i,2}(t,\rho)\cdot (y-y^\prime)+\part_\alpha f_i(t,x,\rho,\hat{\alpha}_i)-
\part_\alpha f_i(t,x^\prime,\rho,\hat{\alpha}_i)\rangle \geq 2\lambda |\hat{\alpha}_i-\hat{\alpha}_i^\prime|^2, \nn
\eea
and thus
$|\hat{\alpha}_i-\hat{\alpha}_i^\prime|\leq (2\lambda)^{-1}\bigl(|b_{i,2}(t,\rho)||y-y^\prime|+
|\part_\alpha f_i(t,x,\rho,\hat{\alpha}_i)-\part_\alpha f_i(t,x^\prime,\rho,\hat{\alpha}_i)|\bigr)$.
This proves the Lipschitz continuity in $(x,y)$. The continuity with respect to the measure arguments follows exactly in the same way.
\end{proof}
\end{lemma}

For given flows $\bg{\mu}^1, \bg{\mu}^2\in \calc([0,T];\calp_2(\mbb{R}^d))$, the adjoint equation of the optimal control problem 
$(\ref{mfg-optimal})$  for each population $1\leq i\leq 2$ is given by
\bea
&&dX_t^i=b_i(t,X_t^i,\mu_t^i,\mu_t^j,\hat{\alpha}_i(t,X_t^i,\mu_t^i,\mu_t^j, Y_t^i))dt+
\sigma_i(t,X_t^i,\mu_t^i,\mu_t^j)dW_t^i~,\nn \\
&&dY_t^i=-\part_x H_i(t,X_t^i,\mu_t^i,\mu_t^j,Y_t^i,Z_t^i,\hat{\alpha}_i(t,X_t^i,\mu_t^i,\mu_t^j,Y_t^i))dt+Z_t^i dW_t^i~,
\label{adjoint-non-random}
\eea
with $j\neq i$, $X_0^i=\xi^i\in \mbb{L}^2(\Omega,\calf_0, \mbb{P};\mbb{R}^d)$ and $Y_T^i=\part_x g_i(X_T^i, \mu_T^i,\mu_T^j)$,
which is a $C(L,\lambda)$-Lipschitz FBSDE.
Notice that $H_i$ must be the full Hamiltonian instead of reduced one due to the state dependence in $\sigma_i$.
Here, $\part_x H_i$ has the form:
\be
\part_x H_i(t,x,\mu, \nu, y,z,\alpha):=b_{i,1}(t,\mu,\nu)^\top y+{\rm{tr}}[\sigma_{i,1}(t,\mu,\nu)^\top z]+\part_x f_i(t,x,\mu, \nu,\alpha)~.\nn
\ee
\begin{theorem}
\label{th-fbsde-existence}
Under Assumption {\rm{\tbf{(MFG-a)}}}, for any flows $\bg{\mu}^1, \bg{\mu}^2\in \calc([0,T];\calp_2(\mbb{R}^d))$, the adjoint FBSDE $(\ref{adjoint-non-random})$ 
for each $i\in\{1,2\}$ 
has a unique solution $(\hat{X}_t^{i}, \hat{Y}_t^{i}, \hat{Z}_t^{i})_{t\in[0,T]}\in \mbb{S}^2\times \mbb{S}^2\times \mbb{H}^2$.
Moreover, there exits a  measurable function $u_i^{\bg{\mu}^i,\bg{\mu}^j}:[0,T]\times \mbb{R}^d\rightarrow \mbb{R}^d$
such that,  with some constant $C$ depending only on $(L,\lambda)$, 
\bea
\forall t\in[0,T], \forall x,x^\prime\in \mbb{R}^d, \quad |u_i^{\bg{\mu}^i,\bg{\mu}^j}(t,x)-u_i^{\bg{\mu}^i,\bg{\mu}^j}(t,x^\prime)|
\leq C|x-x^\prime|
\label{mfg-decoupling-Lip}
\eea
and also that  $\forall t\in[0,T]$, $\hat{Y}_t^i=u_i^{\bg{\mu}^i,\bg{\mu}^j}(t,\hat{X}_t^i)$,  $\mbb{P}$-a.s.

If we set $\hat{\bg{\alpha}}^i=(\hat{\alpha}^i_t=\hat{\alpha}_i(t,\hat{X}_t^{i}, \mu^i_t, \mu_t^j,\hat{Y}_t^{i}))_{t\in[0,T]}$,
then for any $\bg{\beta^i}=(\beta_t^i)_{t\in[0,T]}\in \mbb{A}_i$, it holds:
\bea
J_i^{\bg{\mu}^i,\bg{\mu}^j}(\hat{\bg{\alpha}}^i)+\lambda \mbb{E}\int_0^T |\beta_t^i-\hat{\alpha}_t^i|^2 dt
\leq J_i^{\bg{\mu}^i, \bg{\mu}^j}(\bg{\beta}^i)~.
\label{J-convexity}
\eea
\begin{proof}
The last claim regarding the sufficiency of the stochastic maximal principle is well known.
Indeed, if a solution $(\hat{X}^i_t,\hat{Y}^i_t,\hat{Z}^i_t)_{t\in[0,T]}$ exists, then 
using the convexity $g_i(X_T^{i}, \mu_T^i, \mu_T^j)-g_i(\hat{X}_T^{i}, \mu_T^i, \mu_T^j)\geq \langle
(X_T^i-\hat{X}_T^i), \part_x g_i(\hat{X}_T^i,\mu_T^i,\mu_T^j)\rangle$, 
evaluating the expectation $\mbb{E}[\langle X_T^i-\hat{X}_T^i,\hat{Y}_T^i\rangle ]$ by the Ito formula,
making use of the $\lambda$-convexity of the Hamiltonian, we get the desired result. Here, $(X^i_t)_{t\in[0,T]}$ with $X^i_0=\xi^i$ denotes the solution of the SDE $(\ref{adjoint-non-random})$ 
with $\bg{\beta}_i$ instead of $\hat{\bg{\alpha}}_i$ as its control.
See, for example, Theorem 6.4.6 in \cite{Pham}.

The existence of a unique solution to the adjoint FBSDE as wall as the Lipschitz continuous decoupling field follows from a
straightforward extension of Lemma 3.5 in \cite{Carmona-Delarue-MFG}.   
First, since the adjoint FBSDE $(\ref{adjoint-non-random})$ is $C(L,\lambda)$-Lipschitz continuous 
in $(x,y,z)$ and $\sigma_i$ is independent of $Z^i$, Theorem 1.1 in \cite{Delarue-FBSDE} guarantees
the existence of a unique solution for small time $T\leq c$, where $c=c(L,\lambda)$ is a constant depending only on $(L,\lambda)$.
Thus, for a general $T$, we still have the unique solvability on $[T-\del,T]$ 
with $0<\del\leq c$ and any initial condition $\xi^i\in \mbb{L}^2(\Omega,\calf_{t_0},\mbb{P};\mbb{R}^d)$ at $t_0\in[T-\del,T]$.
We let $(X_t^{i,t_0,\xi^i}, Y_t^{i,t_0,\xi^i}, Z_t^{i,t_0,\xi^i})_{t\in[t_0,T]}$ denote this solution.
Following the proof of Theorem 2.6 in \cite{Delarue-FBSDE} (see also Proposition 4.8 in \cite{Carmona-Delarue-1}),
we can establish the existence and uniqueness on the whole $[0,T]$ by connecting the short-term solutions
provided we have
\bea
\forall x, y\in \mbb{R}^d, \quad |Y^{i,t_0,x}_{t_0}-Y^{i,t_0,y}_{t_0}|\leq C|x-y|~,
\label{Y-Lip}
\eea 
for some $C$ independent of $t_0$ and $\del$.  Here, by Blumenthal's zero-one law, $Y_{t_0}^{i,t_0,x}$ and $Y_{t_0}^{i,t_0,y}$
are deterministic. We are now going to prove $(\ref{Y-Lip})$. Let us put
\be
\hat{J}_i^{t_0,x}:=\mbb{E}\Bigl[\int_{t_0}^T f_i(t,X_t^{i,t_0,x},\mu_t^i,\mu_t^j, \hat{\alpha}^{i, t_0,x}_t)dt+
g_i(X_T^{i,t_0,x},\mu^i_T,\mu^j_T)\Bigr]~,\nn
\ee
with $\hat{\alpha}^{i,t_0,x}_t:=\hat{\alpha}_i(t,X_t^{i,t_0,x},\mu_t^i,\mu_t^j,Y_t^{i,t_0,x})$. Then similar arguments
deriving the relation $(\ref{J-convexity})$ gives
\bea
\langle y-x, Y_{t_0}^{t_0, x}\rangle+\hat{J}_i^{t_0,x}+\lambda \mbb{E}\int_{t_0}^T |\hat{\alpha}^{i,t_0,y}_t-\hat{\alpha}^{i,t_0,x}_t|^2 dt
\leq J_i^{t_0,y}~.
\label{th31-eq}
\eea
Exchanging the role of $x$ and $y$ in $(\ref{th31-eq})$ and adding the two inequalities, we obtain that
\bea
2\lambda \mbb{E}\int_{t_0}^T |\hat{\alpha}_t^{i,t_0,x}-\hat{\alpha}_t^{i,t_0,y}|^2 dt\leq \langle x-y, Y_{t_0}^{i,t_0,x}-Y_{t_0}^{i,t_0,y}\rangle~.
\label{th31-y}
\eea

Now, treating the controls $\hat{\alpha}^{i,t_0,x}, \hat{\alpha}^{i,t_0,y}$ as well as the forward  variables $X^{i,t_0,x}, X^{i,t_0, y}$ as
external inputs, we apply the standard stability result of Lipschitz BSDEs (e.g., see Theorem 4.2.3 in \cite{Zhang-BSDE})
to obtain the estimate for $\mbb{E}[\sup_{t \in [t_0,T]}|Y_t^{i,t_0,x}-Y_t^{i,t_0,y}|^2]$.
Then, applying the standard stability result of  Lipschitz SDEs (e.g., see Theorem 3.24 in \cite{Zhang-BSDE}) to this estimate gives
\bea
&&\mbb{E}\bigl[\sup_{t\in[t_0,T]}|X_t^{i,t_0,x}-X_t^{i,t_0,y}|^2+\sup_{t\in[t_0,T]}|Y_t^{i,t_0,x}-Y_t^{i,t_0,y}|^2\bigr]\nn \\
&&\quad \leq C(L)\Bigl(|x-y|^2+\mbb{E}\int_{t_0}^T |\hat{\alpha}_t^{i,t_0,x}-
\hat{\alpha}_t^{i,t_0,y}|^2 dt\Bigr)~.\nn
\eea
Now the inequality $(\ref{th31-y})$ proves the relation $(\ref{Y-Lip})$ with $C$ depending only on $(L,\lambda)$, 
and hence also the existence of a unique solution for general $T$.
The decoupling field is defined by $u_i^{\bg{\mu}^i,\bg{\mu}^j}:[0,T]\times \mbb{R}^d\ni(t,x) \mapsto Y_t^{i,t,x}$,
and the  representation $\mbb{P}\bigl(\forall t\in[0,T]$,  $\hat{Y}_t^i=u_i^{\bg{\mu^i},\bg{\mu}^j}(t,\hat{X}_t^i)\bigr)=1$ follows from 
the uniqueness of the solution as well as its continuity (Corollary 1.5 in \cite{Delarue-FBSDE}). Its Lipschitz continuity is 
a direct result of $(\ref{Y-Lip})$.
\end{proof}
\end{theorem}

\begin{remark}
In the remainder, we often use the simpler notation $u_i$ for the decoupling field without the superscripts $(\bg{\mu}^i,\bg{\mu}^j)$.
\end{remark}

\begin{lemma}
\label{lemma-stability}
Suppose that two set of functions $(b_i, \sigma_i, f_i, g_i)$ and $(b_i^\prime, \sigma_i^\prime, f_i^\prime,g_i^\prime)$
satisfy Assumption {\rm{\tbf{(MFG-a)}}}. For given inputs $\xi^i, \xi^{\prime, i}\in \mbb{L}^2(\Omega, \calf_0, \mbb{P};\mbb{R}^d)$
and $(\bg{\mu}^i, \bg{\mu}^j), (\bg{\mu}^{\prime, i},\bg{\mu}^{\prime, j})\in \calc([0,T];\calp_2(\mbb{R}^d))^2$, 
let us denote the corresponding solution to $(\ref{adjoint-non-random})$ by $(X^i_t,Y^i_t,Z^i_t)_{t\in[0,T]}$
and $(X^{\prime,i}_t, Y^{\prime, i}_t, Z^{\prime, i}_t)_{t\in[0,T]}$, respectively.
Then, there exists a constant $C$ depending only on $(L,\lambda)$ such that
\bea
\label{stability-mfg}
&&\mbb{E}\Bigl[\sup_{t\in[0,T]}|X_t^{i}-X_t^{\prime, i}|^2+\sup_{t\in[0,T]}|Y_t^i-Y_t^{\prime, i}|^2
+\int_0^T|Z_t^i-Z_t^{\prime,i}|^2 dt\Bigr]\nn \\
&&\leq C\mbb{E}\Bigl[|\xi^i-\xi^{\prime,i}|^2
+|\part_x g_i(X_T^i, \mu_T^i, \mu_T^j)-\part_x g_i^\prime(X_T^i, \mu_T^{\prime,i},\mu_T^{\prime,j})|^2\nn \\
&&\quad+\int_0^T|b_i(t,X_t^i, \mu_t^i, \mu_t^j, \hat{\alpha}_i(t,X_t^i, \mu_t^i,\mu_t^j,Y_t^i))
-b_i^\prime (t,X_t^i, \mu_t^{\prime, i}, \mu_t^{\prime, j}, \hat{\alpha}_i^\prime (t,X_t^i, \mu_t^{\prime,i},\mu_t^{\prime,j},Y_t^i))|^2dt \nn\\
&&
\quad+\int_0^T |\sigma_i(t,X_t^i,\mu_t^i,\mu_t^j)-\sigma_{i}^\prime(t,X_t^i,\mu_t^{\prime,i}, \mu_t^{\prime,j})|^2 dt \nn
+\int_0^T \Bigl(|\part_x H_i(t,X_t^i, \mu_t^i, \mu_t^j, Y_t^i, Z_t^i, \hat{\alpha}_i(t,X_t^i, \mu_t^i, \mu_t^j, Y_t^i)) \nn \\
&&\hspace{15mm}-\part_x H_i^\prime (t,X_t^i, \mu_t^{\prime, i},\mu_t^{\prime, j}, Y_t^i, Z_t^i, \hat{\alpha}_i^\prime (t,X_t^i,\mu_t^{\prime, i},\mu_t^{\prime,j}, Y_t^i))|^2\Bigr) dt \Bigr],
\eea
where the functions $H_i, H_i^\prime:[0,T]\times \mbb{R}^d\times\calp_2(\mbb{R}^d)^2\times \mbb{R}^d\times \mbb{R}^{d\times d}\times A_i
\rightarrow \mbb{R}$
are the Hamiltonians $(\ref{def-H})$ associated with the coefficients $(b_i, \sigma_i, f_i)$ and $(b_i^\prime, \sigma_i^\prime, f_i^\prime)$, respectively, 
and $\hat{\alpha}_i, \hat{\alpha}_i^\prime$ are their minimizers.
In particular,  there is another constant $C^\prime$ depending additionally on  $K$ such that
\bea
\mbb{E}\Bigl[\sup_{t\in[0,T]}|X_t^i|^2+\sup_{t\in[0,T]}|Y_t^i|^2+\int_0^T |Z_t^i|^2dt\Bigr]\leq C\Bigl(||\xi^i||_2^2+\sup_{t\in[0,T]}\sum_{j=1}^2M_2(\mu_t^j)^2\Bigr)+C^\prime, 
\label{mfg-l2-growth}
\eea
and, for any $t\in[0,T]$, 
\be
|Y_t^i|\leq C\Bigl(|X_t^i|+\sup_{t\in[0,T]}\sum_{j=1}^2 M_2(\mu_t^j)\Bigr)+C^\prime, ~\mbb{P}\text{-a.s.}
\label{decoupling-growth}
\ee
\begin{proof}
Since the FBSDE $(\ref{adjoint-non-random})$ has Lipschitz continuous coefficients, 
 it is standard to show that there exists some constant $c$ depending only 
on $(L,\lambda)$ such that the estimate $(\ref{stability-mfg})$ holds for small $T\leq c$.
In particular, by applying Ito formula to $|Y_t^i-Y^{\prime, i}_t|^2$, we see
that $\mbb{E}[\sup_{t\in[0,T]}|Y_t^i-Y^{\prime,i}_t|^2+\int_0^T |Z_t^i-Z_t^{\prime,i}|^2]$
is bounded by the terms related to the backward equation in $(\ref{stability-mfg})$ plus  the term $C\mbb{E}[\sup_{t\in[0,T]}|X_t^i-X_t^{\prime, i}|^2]$.
On the other hand, the similar calculation shows that $\mbb{E}[\sup_{t\in[0,T]}|X_t^i-X_t^{\prime,i}|^2]$
is bounded by the remaining terms in $(\ref{stability-mfg})$ plus $C T \mbb{E}[\sup_{t\in[0,T]}|Y_t^{i}-Y_t^{i,\prime}|^2]$
where $C$ depends only on the Lipschitz constants of the system. Hence, for small $T\leq c$, we obtain the 
desired estimate.

For general $T$, the estimate is a result of connecting the short-term estimates (Lemma 4.9~\cite{Carmona-Delarue-1}).
Since the same technique will be used  also in Lemma~\ref{lemma-stability-mkv}, let us explain it here in details.
We first divide the interval $[0,T]$ into a finite number of  subintervals $([T_{k-1},T_{k}])_{1\leq k\leq N}$ with $T_0=0, T_N=T$ and 
$T_k-T_{k-1}\leq c$ for each $k$. The estimate for $\Theta(T_{k-1},T_k):=\mbb{E}\bigl[\sup_{t\in[T_{k-1},T_k]}|X_t^i-X_t^{\prime,i}|^2+
\sup_{t\in[T_{k-1},T_k]}|Y_t^i-Y_t^{\prime,i}|^2+\int_{T_{k-1}}^{T_k}|Z_t^i-Z_t^{\prime,i}|^2 dt\bigr]$ on each interval 
can be written in the from
\bea
\Theta(T_{k-1},T_k)&\leq& C\mbb{E}\Bigl[|X_{T_{k-1}}^i-X_{T_{k-1}}^{\prime,i}|^2+|(u_i-u_i^\prime)(T_k,X_{T_k}^i)|^2\nn \\
&&\qquad+\int_{T_{k-1}}^{T_k}|(b_i-b_i^\prime,\sigma_i-\sigma_i^\prime,\part_x H_i-\part_x H_i^\prime)(s,X_s^i,Y_s^i,Z_s^i)|^2 ds\Bigr]\nn
\eea
where we have used the notation $u_i:=u_i^{\bg{\mu}^i,\bg{\mu}^j}, u_i^\prime:=u_i^{\prime, \bg{\mu}^{\prime, i}, \bg{\mu}^{\prime,j}}$
and omitted the arguments regarding $(\bg{\mu}^j,\bg{\mu}^{\prime,j})_{1\leq j\leq 2}$ to 
lighten the expression.  The Lipschitz continuity in $(\ref{mfg-decoupling-Lip})$
is crucial to derive this expression.  For $k=N$, it gives
$
\Theta(T_{N-1},T_N)\leq C\mbb{E}\Bigl[|X_{T_{N-1}}^i-X_{T_{N-1}}^{\prime,i}|^2+|\del_T|^2+\int_{T_{N-1}}^{T}|\del h_s|^2 ds\Bigr]
$
where $\del_T:=(\part_x g_i-\part_x g_i^\prime)(T,X_T^i)$ and $\del h_s:=(b_i-b_i^\prime,\sigma_i-\sigma_i^\prime,
\part_x H_i-\part_x H_i^\prime)(s,X_s^i, Y_s^i,Z_s^i)$.
This means, in particular, 
$\mbb{E}\bigl[|(u_i-u_i^\prime)(T_{N-1}, X_{T_{N-1}}^i)|^2\bigr]\leq 
 C\mbb{E}\Bigl[|X_{T_{N-1}}^i-X_{T_{N-1}}^{\prime,i}|^2+|\del_T|^2+\int_{T_{N-1}}^{T}|\del h_s|^2 ds\Bigr]$.
Since it holds for any {\it initial value} $X_{T_{N-1}}^{\prime, i}$, we obtain
$\mbb{E}\bigl[|(u_i-u_i^\prime)(T_{N-1}, X_{T_{N-1}}^i)|^2\bigr]\leq 
 C\mbb{E}\Bigl[|\del_T|^2+\int_{T_{N-1}}^{T_N}|\del h_s|^2 ds\Bigr]$ by choosing $X_{T_{N-1}}^{\prime,i}=X_{T_{N-1}}^i$.
This estimate then implies $\Theta(T_{N-2},T_{N-1})\leq C\mbb{E}\Bigl[|X_{T_{N-2}}^i-X_{T_{N-2}}^{\prime,i}|^2+
|\del_T|^2+\int_{T_{N-2}}^T|\del h_s|^2 ds\Bigr]$.  By iteration, we get for any $k$,
\bea
\Theta(T_k,T_{k+1})\leq C\mbb{E}\Bigl[|X_{T_k}^i-X_{T_k}^{\prime,i}|^2+|\del_T|^2+\int_{T_k}^T |\del h_s|^2 ds\Bigr]~.
\label{stable-eq-1}
\eea
Moreover, by iterating the relation $\mbb{E}\bigl[|X^i_{T_k}-X^{\prime, i}_{T_k}|^2\bigr]\leq 
C\mbb{E}\Bigl[|X_{T_{k-1}}^i-X_{T_{k-1}}^{\prime,i}|^2+|\del_T|^2+\int_{T_{k-1}}^T |\del h_s|^2 ds\Bigr]$, we get
\bea
\mbb{E}\bigl[|X_{T_k}^i-X_{T_k}^{\prime,i}|^2\bigr]\leq C\mbb{E}\Bigl[|\xi^i-\xi^{\prime, i}|^2+
|\del_T|^2+\int_0^T|\del h_s|^2 ds\Bigr]~.
\label{stable-eq-2}
\eea
Inserting the estimate $(\ref{stable-eq-2})$ into $(\ref{stable-eq-1})$ and summing over $k$, 
we obtain the desired estimate.

In order to obtain the growth estimate, we put, for any $(t,x,\mu,\nu,\alpha)\in [0,T]\times \mbb{R}^d\times \calp_2(\mbb{R}^d)^2\times A_i$,
\bea
&&b_i(t,x,\mu,\nu,\alpha)=\sigma_i(t,x,\mu,\nu)=g_i(x,\mu,\nu)=0~,\quad f_i(t,x,\mu, \nu,\alpha)=\lambda |\alpha|^2~, \nn
\eea
and $\xi^i=0$, which then satisfies Assumption \tbf{(MFG-a)} and makes $(X^i,Y^i,Z^i)$ identically zero. Plugging them into $(\ref{stability-mfg})$,
we obtain the estimate:
\bea
&&\mbb{E}\Bigl[\sup_{t\in[0,T]}|X_t^{\prime,i}|^2+\sup_{t\in[0,T]}|Y_t^{\prime,i}|^2
+\int_0^T|Z_t^{\prime,i}|^2 dt\Bigr]\nn \\
&&\leq C\Bigl(||\xi^{\prime, i}||_2^2+|\part_x g_i^\prime(0,\mu_T^{\prime,i},\mu_T^{\prime,j})|^2 
+\int_0^T |\part_x f_i^\prime (t,0,\mu_t^{\prime,i}, \mu_t^{\prime,j}, \hat{\alpha}_i^\prime (t,0,\mu_t^{\prime, i}, \mu_t^{\prime,j},0))|^2 dt \nn \\
&&+\int_0^T\bigl[|b_i^\prime (t,0,\mu^{\prime,i}_t,\mu^{\prime,j}_t,\hat{\alpha}_i^\prime(t,0,\mu^{\prime,i}_t,\mu^{\prime,j}_t,0))|^2
+|\sigma_i^\prime(t,0,\mu^{\prime,i}_t,\mu^{\prime,j}_t)|^2\bigr]dt\Bigr)~.
\label{mfg-gw-estimate}
\eea
Now, by symmetry, the desired estimate $(\ref{mfg-l2-growth})$ holds for $(X^i, Y^i, Z^i)$.
Finally, using the initial condition $X_t^i=0$ at time $t$ yields
\be
|Y_t^{i,t,0}|\leq \mbb{E}\Bigl[\sup_{s\in[t,T]}|Y_s^{i,t,0}|^2\Bigr]^\frac{1}{2}\leq C\sup_{s\in[t,T]}\sum_{j=1}^2 M_2(\mu_s^j)+C^\prime.\nn
\ee
Now, by  the Lipschitz continuity $(\ref{mfg-decoupling-Lip})$(or equivalently $(\ref{Y-Lip})$), we have $|Y_t^i-Y_t^{i,t,0}|=|Y_t^{i,t,X_t^{i}}-Y_t^{i,t,0}|\leq C|X_t^i|$. This proves the growth estimate $(\ref{decoupling-growth})$.
\end{proof}
\end{lemma}

\subsection{MFG equilibrium under boundedness assumptions}
In the  preceding subsections, we have seen that, for given deterministic flows of probability measures 
$\bg{\mu}^1, \bg{\mu}^2\in \calc([0,T];\calp_2(\mbb{R}^d))$, the solution to each optimal control problem of $(\ref{mfg-optimal})$
is characterized by the uniquely solvable FBSDE $(\ref{adjoint-non-random})$.
Hence finding an equilibrium condition  $(\ref{mfg-fixed-point})$ results in finding a solution to 
the following system of FBSDEs of McKean-Vlasov (MKV) type: for $i,j\in\{1,2\}, j\neq i$;
\bea
&&dX_t^i=b_i\bigl(t,X_t^i, \call(X_t^i),\call(X_t^j),\hat{\alpha}_i(t,X_t^i, \call(X_t^i),\call(X_t^j), Y_t^i)\bigr)dt+\sigma_i\bigl(t,X_t^i,
\call(X_t^i),\call(X_t^j)\bigr)dW_t^i, \nn \\
&&dY_t^i=-\part_x H_i\bigl(t,X_t^i, \call(X_t^i), \call(X_t^j), Y_t^i, Z_t^i, \hat{\alpha}_i(t,X_t^i,\call(X_t^i), \call(X_t^j), Y_t^i)\bigr)dt
+Z_t^i dW_t^i, 
\label{mfg-fbsde}
\eea
with $X_0^i=\xi^i\in \mbb{L}^2(\Omega,\calf_0,\mbb{P};\mbb{R}^d)$ and $Y_T^i=\part_x g_i(X_T^i, \call(X_T^i),\call(X_T^j))$.
Although two MKV-type FBSDEs are now coupled, their interactions appear only through the laws of the two populations.
Thanks to this property, we can still apply a similar strategy developed by Carmona \& Delarue \cite{
Carmona-Delarue-MFG,Carmona-Delarue-MKV}.
A crucial tool to prove the existence of an equilibrium is the Schauder's fixed point theorem \cite{Schauder} generalized by 
Tychonoff~\cite{Tychonoff}\footnote{See, for example Shapiro~\cite{Shapiro},  for pedagogical introduction of  
fixed-point theorems and relevant references.}. 
The following form is taken from Theorem 4.32 in \cite{Carmona-Delarue-1}.
\begin{theorem} (Schauder FPT)
\label{th-Schauder}
Let $(V, ||\cdot||)$ be a normed linear vector space and $E$ be a nonempty closed convex subset of $V$.
Then, any continuous mapping from $E$ into itself which has a relatively compact range has a 
fixed point.
\end{theorem}

In this subsection, we prove the existence of a solution to the system of FBSDEs $(\ref{mfg-fbsde})$
under additional assumptions.
\begin{assumption} {\rm{\tbf{(MFG-b)}}} 
For $1\leq i\leq 2$, 
there exist some  element $0_{A_i}\in A_i$ and a constant $\Lambda$ such that,
for any $(t,\mu,\nu)\in[0,T]\times \calp_2(\mbb{R}^d)^2$,
\bea
&&|b_{i,0}(t,\mu,\nu)|, |\sigma_{i,0}(t,\mu,\nu)|\leq \Lambda~, \nn \\
&&|\part_x g_i(0,\mu,\nu)|, |\part_{(x,\alpha)} f_i (t,0,\mu, \nu,0_{A_i})|\leq \Lambda~.\nn 
\eea
\end{assumption}
Here is the main result of this subsection.
\begin{theorem}
\label{th-mfg-bounded}
Under Assumptions {\rm{\tbf{(MFG-a,b)}}}, the system of FBSDEs $(\ref{mfg-fbsde})$ (and hence the matching problem~$(\ref{mfg-fixed-point})$) 
 is solvable for any $\xi^1,\xi^2\in\mbb{L}^2(\Omega,\calf_0, \mbb{P};\mbb{R}^d)$.

\begin{proof}
With slight abuse of notation, we let $(X^{i,\bg{\rho}}_t, Y^{i,\bg{\rho}}_t, Z^{i,\bg{\rho}}_t)_{t\in[0,T]}$ denote the solution to the FBSDE $(\ref{adjoint-non-random})$ for a given flows $\bg{\rho}:=(\bg{\mu}^1,\bg{\mu}^2)\in \calc([0,T];\calp_2(\mbb{R}^d))^2$ and the
initial condition $X_0^{i,\bg{\rho}}=\xi^i$. By Theorem~\ref{th-fbsde-existence}, we can define a map:
\be
\Phi:\calc([0,T];\calp_2(\mbb{R}^d))^2\ni (\bg{\mu}^1, \bg{\mu}^2)\mapsto (\call(X_t^{1,\bg{\rho}})_{t\in[0,T]},
\call(X_t^{2,\bg{\rho}})_{t\in[0,T]})\in \calc([0,T];\calp_2(\mbb{R}^d))^2~.\nn
\ee
In the following, we are going to check the conditions necessary for the application of Schauder FPT to this map.
As a linear vector space $V$ in the FPT, we use the product space $\calc([0,T];\calm_f^1(\mbb{R}^d))^2$
equipped with the supremum of the Kantorovich-Rubinstein norm:
\bea 
\label{KR-norm}
&&||(\bg{\mu}^1,\bg{\mu}^2)||:=\sum_{i=1}^2 \sup_{t\in[0,T]}||\mu^i_t||_{\rm{KR}_*},  \\
&&{\text{with}}: \quad ||\mu||_{{\rm{KR}}_*}:=|\mu(\mbb{R}^d)|+\sup\Bigl\{ \int_{\mbb{R}^d}l(x)\mu(dx);~l\in {\rm Lip}_1(\mbb{R}^d), l(0)=0\Bigr\}~,\nn
\eea
where ${\rm{Lip}}_1(\mbb{R}^d)$ is the set of $1$-Lipschitz continuous functions on $\mbb{R}^d$. Importantly, the norm $||\cdot ||_{\rm{KR}_*}$
is known to coincide with the 1-Wasserstein distance $W_1$ on $\calp_1(\mbb{R}^d)$ (Corollary 5.4 in \cite{Carmona-Delarue-1}).
Of course, the reason to use a space of signed measures is to make it linear.

From $(\ref{alpha-growth})$ with $(\beta_i=0_{A_i})$, it is immediate to see that $|\hat{\alpha}_i(t,0,\rho,0)|\leq C(\lambda, \Lambda)$.
Hence,  by using the estimate $(\ref{mfg-gw-estimate})$ in Lemma~\ref{lemma-stability},  we get
$\mbb{E}\bigl[\sup_{t\in[0,T]}|Y_t^{i, \bg{\rho}}|^2\bigr]\leq C(L,\lambda,\Lambda)(1+||\xi^i||_2^2)$.
Then the Lipschitz continuity of the decoupling field implies that 
$|Y_t^{i,\bg{\rho}}|\leq C\bigl(1+|X_t^{i,\bg{\rho}}|\bigr)$
with $C$ independent of $\bg{\rho}$. 
Therefore, again by $(\ref{alpha-growth})$, we have
$
|\hat{\alpha}_i(t,X_t^{i,\bg{\rho}},\rho_t,Y_t^{i,\bg{\rho}})|\leq C( 1+|X_t^{i,\bg{\rho}}|)$.
Now, it is standard to check that
$\mbb{E}\bigl[|X_t^{i,\bg{\rho}}-X_s^{i,\bg{\rho}}|^2\bigr]\leq C|t-s|$
and hence 
\be
W_2(\call(X_t^{i,\bg{\rho}}),\call(X_s^{i,\bg{\rho}}))\leq C|t-s|^\frac{1}{2}~
\label{x-continuity}
\ee
holds uniformly in $\bg{\rho}$. 
Since $\hat{\alpha}_i$ is of linear growth in $X_t^{i,\bg{\rho}}$,
it is also straightforward to obtain 
$\mbb{E}\Bigl[\sup_{t\in[0,T]}|X_t^{i,\bg{\rho}}|^4|\calf_0\Bigr]^\frac{1}{2}\leq  C\bigl(1+|\xi^i|^2)$
uniformly in $\bg{\rho}$. This inequality guarantees the uniform square integrability. In fact, the following estimate holds uniformly in 
$\bg{\rho}$ with any $a\geq 1$;\footnote{See p. 259 in \cite{Carmona-Delarue-1}.}
\bea
\mbb{E}\Bigl[\sup_{t\in[0,T]}|X_t^{i,\bg{\rho}}|^2\bold{1}_{\{\sup_{t\in[0,T]}|X_t^{i,\bg{\rho}}|\geq a\}}\Bigr]
\leq C\Bigl(a^{-1}+\mbb{E}\bigl[|\xi^i|^2\bold{1}_{\{|\xi^i|\geq \sqrt{a}\}}\bigr]\Bigr)^\frac{1}{2}.
\label{uniform-square-integrability}
\eea
Since the relation will be used repeatedly in the following, let us explain it here.
For any $D\in \calf$ and $\ep>0$, we have, by Cauchy-Schwarz inequality,
\bea
&&\mbb{E}\bigl[\sup_{t\in[0,T]}|X_t^{i,\bg{\rho}}|^2\bold{1}_{D}\bigr]\leq C\mbb{E}\bigl[(1+|\xi^i|^2)\mbb{P}(D|\calf_0)^\frac{1}{2}\bigr]\nn \\
&&\quad \leq C\Bigl(\ep+\ep^{-1}\mbb{E}\bigl[(1+|\xi^i|^2)\bold{1}_D\bigr]\Bigr)
\leq C\mbb{E}\bigl[(1+|\xi^i|^2)\bold{1}_D\bigr]^\frac{1}{2}. \nn
\eea
In the last inequality, we have maximized in $\ep$.  Here, $C$ depends on 
$||\xi^i||_2$ but not on $\bg{\rho}$.  We also have
\bea
\sup_{D\in \calf; \mbb{P}(D)\leq Ca^{-2}}\mbb{E}\bigl[(1+|\xi^i|^2)\bold{1}_D\bigr]&\leq&Ca^{-2}+\mbb{E}\bigl[|\xi^i|^2(\bold{1}_{\{|\xi^i|\leq \sqrt{a}\}}+\bold{1}_{\{|\xi^i|\geq \sqrt{a}\}})\bold{1}_D\bigr]\nn \\
&\leq &2Ca^{-1}+\mbb{E}\bigl[|\xi^i|^2\bold{1}_{\{|\xi^i|\geq \sqrt{a}\}}\bigr]~. \nn
\eea
Since $\mbb{P}\Bigl(\sup_{t\in[0,T]}|X_t^{i,\bg{\rho}}|\geq a\Bigr)\leq Ca^{-2}$ by Chevyshev's inequality,
the estimate $(\ref{uniform-square-integrability})$ is now established.

The above estimate suggests us to restrict the map $\Phi$ to the following domain:
\bea
&&\cale:=\Bigl\{(\bg{\mu}^1, \bg{\mu}^2)\in \calc([0,T];\calp_2(\mbb{R}^d))^2;\nn \\
&&\hspace{10mm} \forall a\geq 1, 1\leq i\leq 2, \sup_{t\in[0,T]}\int_{|x|\geq a} |x|^2 \mu^i_t(dx)
\leq C\Bigl(a^{-1}+\mbb{E}\bigl[|\xi^i|^2\bold{1}_{\{|\xi^i|\geq \sqrt{a}\}}\bigr]\Bigr)^\frac{1}{2}\Bigr\},\nn
\eea
which is a closed and convex subset of $C([0,T];\calm_f^1(\mbb{R}^d))^2$.
Choosing $C$ sufficiently large, we can make $\Phi$ a self-map on $\cale$.
By the estimate $(\ref{uniform-square-integrability})$ and Corollary 5.6 in \cite{Carmona-Delarue-1}, 
there exists a compact subset $\calk\subset \calp_2(\mbb{R}^d)^2$ such that $\forall t\in[0,T], ~[\Phi(\bg{\rho})]_t\in \calk$
for any $\bg{\rho}\in \cale$. Combined with the equicontinuity $(\ref{x-continuity})$, Arzela-Ascoli theorem 
implies that the image $\Phi(\cale)$ is a relatively compact subset of $\calc([0,T];\calp_2(\mbb{R}^d))^2$, 
and in particular of $\calc([0,T];\calp_1(\mbb{R}^d))^2$.

Finally, by Lemma~\ref{lemma-stability} and also by the continuity of coefficients in the measure arguments in $W_2$-distance,
the dominated convergence theorem implies that
$\mbb{E}\bigl[\sup_{t\in[0,T]}|X_t^{i,\bg{\rho}}-X_t^{i,\bg{\rho}^\prime}|^2\bigr]\rightarrow 0$
when $1\leq i\leq 2, \forall t\in[0,T], W_2(\mu^i_t,\mu_t^{\prime,i})\rightarrow 0$.
Note that, by Theorem 5.5 in \cite{Carmona-Delarue-1}, when $\bg{\rho}$ converge with respect to the norm $||\cdot||$ in $(\ref{KR-norm})$ under 
the restriction to the domain $\cale$, $\bg{\rho}$ actually converges in $W_2$-distance.
This proves the continuity of the map $\Phi$. 
Now the existence of a fixed point (not necessarily unique) of the map $\Phi$ is guaranteed by Schauder FPT, 
which provides a solution to the system of FBSDEs $(\ref{mfg-fbsde})$.
\end{proof}
\end{theorem}

\subsection{MFG equilibrium for small $T$ or small coupling}
It is well known that the Lipschitz continuity of coefficient functions is not enough 
to make a coupled FBSDE well-defined for a general terminal time $T$. See, for example,
the discussions in \cite{Delarue-FBSDE,Ma-Yong} and,  in particular,  Section 3.2.3 of \cite{Carmona-Delarue-1}.
We have already used convexity conditions to overcome this problem in the last section.
In order to allow the quadratic cost functions relevant for popular Linear-Quadratic problems, we want to relax 
Assumption {\rm{\tbf{(MFG-b)}}}. This is exactly what Carmona \& Delarue have done in \cite{Carmona-Delarue-MFG}
for single population. Although we can follow the same route, it requires much stronger assumptions than {\rm{\tbf{(MFG-a)}}.
Unfortunately, the conditions required in \cite{Carmona-Delarue-MFG} preclude most of the interesting interactions among different populations through their state dynamics.
In this work, in order to allow flexible interactions among populations and also to be
complementary to the result in \cite{Carmona-Delarue-MFG}, we focus on the problems with small $T$.
Requiring small $T$ is a reasonable trade-off for quadratic interactions by considering the fact that, even for a deterministic LQ-problem, the 
relevant Riccati equation  may diverge within a finite time.
After the analysis for small $T$, we provide another solution which allows general $T$ but requires the couplings between FSDE and BSDE are small enough.

\begin{theorem}
\label{th-mfg-main}
Under Assumption {\rm{\tbf{(MFG-a)}}},  there exists some positive constant $c$ depending only on $(L,\lambda)$
such that, for any $T\leq c$, the system of FBSDEs $(\ref{mfg-fbsde})$ (and hence the matching problem~$(\ref{mfg-fixed-point})$) 
is solvable for any $\xi^1, \xi^2 \in\mbb{L}^2(\Omega,\calf_0, \mbb{P};\mbb{R}^d)$.
\begin{proof}

For any $n\in \mbb{N}$ and $\mu\in\calp_2(\mbb{R}^d)$, let us define $\phi_n\circ \mu$ as a push-forward of $\mu$
by the map $\displaystyle \mbb{R}^d\ni x\mapsto \frac{n x}{\max(M_2(\mu),n)}$. In other words, for any random variable $X$
with $\call(X)=\mu$, the law of $\displaystyle \frac{n X}{\max(M_2(\mu), n)}$ is given by $\phi_n\circ \mu$.
Obviously, $M_2(\phi_n\circ \mu)\leq n$ and the map $\mu\mapsto \phi_n\circ \mu$ is continuous with respect to $W_2$-distance.
Using this map, we introduce a sequence of approximated functions
\bea
&&(b_{i,0}^n,b_{i,2}^n,\sigma^n_{i,0})(t,\mu,\nu):=(b_{i,0},b_{i,2},\sigma_{i,0})(t,\phi_n\circ \mu, \phi_n\circ \nu)~,\nn \\
&&f_i^n(t,x,\mu, \nu, \alpha):=f_i(t,x,\phi_n\circ \mu, \phi_n\circ \nu, \alpha), \quad
g_i^n(x,\mu,\nu):=g_i(x,\phi_n\circ \mu,\phi_n\circ \nu), \nn
\eea
and accordingly define
\bea
&&b_i^n(t,x,\mu,\nu,\alpha):=b_{i,0}^n(t,\mu,\nu)+b_{i,1}(t,\mu,\nu)x+b_{i,2}^n(t,\mu,\nu)\alpha~,\nn \\
&&\sigma_i^n(t,x,\mu,\nu):=\sigma_{i,0}^n(t,\mu,\nu)+\sigma_{i,1}(t,\mu,\nu)x~,\nn \\
&&H_i^n(t,x,\mu,\nu,y,z,\alpha):=\langle b_i^n(t,x,\mu,\nu,\alpha),y\rangle+{\rm tr}[\sigma_i^n(t,x,\mu,\nu)^\top z]
+f_i^n(t,x,\mu,\nu,\alpha)~,\nn
\eea
for any $(t,x,\mu,\nu,y,z,\alpha)\in[0,T]\times \mbb{R}^d\times \calp_2(\mbb{R}^d)^2\times \mbb{R}^d\times \mbb{R}^{d\times d} 
\times A_i$. Since 
$\part_\alpha H_i^n(t,x,\mu, \nu, y,z,\alpha)=b_{i,2}(t,\phi_n\circ\mu,\phi_n\circ \nu)^\top y+\part_\alpha f_i(t,x,\phi_n\circ \mu, \phi_n\circ \nu,\alpha)$,
the minimizer given as a solution to the variational inequality $(\ref{alpha-vi})$ satisfies
\be
\hat{\alpha}_i^n(t,x,\mu,\nu,y)=\hat{\alpha}_i(t,x,\phi_n\circ \mu,\phi_n\circ \nu, y)~, 
\label{alpha-n-simple}
\ee 
where $\hat{\alpha}_i$ is the minimizer of the original Hamiltonian $H_i$. 
The regularization for $b_{i,2}$ is done solely to obtain the simple expression $(\ref{alpha-n-simple})$
for the minimizer.

The new coefficient functions $(b_i^n, \sigma_i^n, f_i^n,g_i^n)$ clearly
satisfy  {\rm{\tbf{(MFG-a,b)}}} for each $n$.
Thus Theorem~\ref{th-mfg-bounded} guarantees  that there exists a solution to the following system of 
FBSDEs of MKV-type with $i,j\in\{1,2\}$, $j\neq i$:
\bea
\label{mfg-fbsde-approximated}
&&\hspace{-5mm}dX_t^{i,n}=b_i^n(t,X_t^{i,n},\call(X_t^{i,n}),\call(X_t^{j,n}),\hat{\alpha}_i^n(t,X_t^{i,n},\call(X_t^{i,n}),\call(X_t^{j,n}), Y_t^{i,n}))dt \nn\\
&&\qquad+\sigma_i^n(t,X_t^{i,n},\call(X_t^{i,n}),\call(X_t^{j,n}))dW_t^i~, \\
&&\hspace{-5mm}dY_t^{i,n}=-\part_x H_i^n(t,X_t^{i,n},\call(X_t^{i,n}),\call(X_t^{j,n}), Y_t^{i,n},Z_t^{i,n},
\hat{\alpha}_i^n(t,X_t^{i,n},\call(X_t^{i,n}),\call(X_t^{j,n}),Y_t^{i,n}))dt+Z_t^{i,n}dW_t^i~,\nn 
\eea
with $X_0^{i,n}=\xi^i$ and $Y_T^{i,n}=\part_x g_i^n(X_T^{i,n},\call(X_T^{i,n}),\call(X_T^{j,n}))$.
Treating $(\call(X_t^{i,n}))_{1\leq i\leq 2}$ as inputs,  we can see that there exist 
some constants $C=C(L,\lambda)$ and $C^\prime=C^\prime(L,\lambda, K)$ such that
\be
|Y_t^{i,n}|\leq C\Bigl(|X_t^{i,n}|+\sup_{s\in[t,T]}\sum_{j=1}^2 M_2(\call(X_s^{j,n}))\Bigr)+C^\prime.\nn
\ee 
from the  growth estimate in Lemma~\ref{lemma-stability}.
It then follows from Lemma~\ref{lemma-alpha} that
\bea
|\hat{\alpha}_i^n(t,X_t^{i,n},\call(X_t^{i,n}),\call(X_t^{j,n}),Y_t^{i,n})|\leq  C\Bigl(|X_t^{i,n}|+\sup_{s\in[t,T]}\sum_{j=1}^2 M_2(\call(X_s^{j,n}))\Bigr)+C^\prime,\nn
\eea
uniformly in $n$. Then, for any $t\in[0,T]$, it is easy to check that
\bea
\mbb{E}\bigl[|X_t^{i,n}|^2\bigr]&\leq& C\mbb{E}\Bigl[|\xi^i|^2+\int_0^t\bigl[|b_i(s,X_s^{i,n},\call(X_s^{i,n}),\call(X_s^{j,n}),\hat{\alpha}_i^n(t))|^2+
|\sigma_i(t,X_t^{i,n},\call(X_t^{i,n}),\call(X_t^{j,n}))|^2 dt\Bigr]\nn \\
&\leq &C^\prime+C\Bigl(||\xi^i||_2^2+T\sup_{s\in[0,T]}\sum_{j=1}^2M_2(\call(X_s^{j,n}))^2\Bigr)+C\int_0^t\sum_{j=1}^2 \mbb{E}\bigl[|X_s^{j,n}|^2\bigr]ds~,
\label{xsq-ineq-mfg}
\eea 
with $C=C(L,\lambda)$. Applying Gronwall's inequality to the summation over $1\leq i\leq 2$,
we get
\bea
\sup_{t\in[0,T]}\sum_{i=1}^2 \mbb{E}\bigl[|X_t^{i,n}|^2\bigr]\leq C^\prime+C\Bigl(||\bg{\xi}||_2^2+T\sup_{t\in[0,T]}\sum_{i=1}^2 
M_2(\call(X_t^{i,n}))^2\Bigr)~,\nn
\eea
with $\bg{\xi}:=(\xi^1,\xi^2)$. Therefore, there exists a  constant $c$ depending only on $(L,\lambda)$ such that, 
for any $T\leq c$,
\be
\sup_{t\in[0,T]}\sum_{i=1}^2\mbb{E}\bigl[|X_t^{i,n}|^2\bigr]\leq C^\prime \bigl(1+||\bg{\xi}||^2_2\bigr) 
\label{mfg-uniform-bound}
\ee 
uniformly in $n$.

Let us assume $T\leq c$ in the remainder.
From $(\ref{mfg-uniform-bound})$, we can show straightforwardly that
$\mbb{E}\bigl[|X_t^{i,n}-X_s^{i,n}|^2\bigr]\leq C|t-s|$
and $\mbb{E}\bigl[\sup_{t\in[0,T]}|X_t^{i,n}|^4|\calf_0\bigr]^\frac{1}{2}
\leq C \bigl(1+|\xi^i|^2\bigr)$ hold uniformly in $n$.
Just as in $(\ref{uniform-square-integrability})$, we have for any $a\geq 1$, 
\bea
\sup_{n\geq 1} \mbb{E}\Bigl[\sup_{t\in[0,T]}|X_t^{i,n}|^2\bold{1}_{\{\sup_{t\in[0,T]}|X_t^{i,n}|\geq a\}}\Bigr]\leq
C\Bigl(a^{-1}+\mbb{E}\bigl[|\xi^i|^2\bold{1}_{\{|\xi^i|\geq \sqrt{a}\}}\bigr]\Bigr)^\frac{1}{2}. \nn
\eea
Hence, combined with the equicontinuity, we conclude that 
$(\call(X_t^{1,n})_{t\in[0,T]}, \call(X_t^{2,n})_{t\in[0,T]})_{n\geq 1}$ is a relatively compact subset of $\calc([0,T];\calp_2(\mbb{R}^d))^2$.
Therefore, there exists some $(\bg{\mu}^1, \bg{\mu}^2) \in \calc([0,T];\calp_2(\mbb{R}^d))^2$ such that,
upon extracting some subsequence (still denoted by $n$),  
\be
\lim_{n\rightarrow\infty} \sup_{t\in[0,T]} W_2(\call(X_t^{i,n}),\mu^i_t)=0, ~i\in\{1,2\}~.\nn 
\ee
Let us define $(X_t^i,Y_t^i,Z_t^i)_{t\in[0,T],1\leq i\leq 2}$ as the solutions to the FBSDEs $(\ref{adjoint-non-random})$
with those $(\bg{\mu}^1, \bg{\mu}^2)$ as inputs. The convergence $\mbb{E}\bigl[\sup_{t\in[0,T]}|X_t^{i,n}-X_t^i|^2\bigr]\rightarrow 0$
as $n\rightarrow \infty$ can be shown by the stability result in Lemma~\ref{lemma-stability}. 
Note that, thanks to the boundedness of $(\ref{mfg-uniform-bound})$, there exists $n_0\in \mbb{N}$ such that
we can replace all the approximated coefficients $(b_i^n, \sigma_i^n, f_i^n, g_i^n, \hat{\alpha}_i^n)$ by 
the original ones $(b_i,\sigma_i,f_i,g_i,\hat{\alpha}_i)$ for any $n\geq n_0$. 
The convergence  $\mbb{E}\bigl[\sup_{t\in[0,T]}|X_t^{i,n}-X_t^i|^2\bigr]\rightarrow 0$
then follows easily by the dominated convergence theorem.
By the inequality
$\sup_{t\in[0,T]}W_2(\call(X_t^i),\mu_t^i)\leq \sup_{t\in[0,T]}\bigl( W_2(\call(X_t^i), \call(X_t^{i,n}))+W_2(\call(X_t^{i,n}),\mu_t^i)\bigr)$, 
the above convergence implies  $\bg{\mu}^i=\call(X_t^i)_{t\in[0,T]}$ $1\leq i\leq 2$.  Therefore, $(X_t^i,Y_t^i,Z_t^i)_{t\in[0,T],1\leq i\leq 2}$
is actually a wanted solution to the system of FBSDEs $(\ref{mfg-fbsde})$.
\end{proof}
\end{theorem}

Another simple method to allow the quadratic cost functions is making the couplings between FSDE and BSDE small enough.
\begin{theorem}
\label{th-mfg-small-coupling}
Under Assumption {\rm{\tbf{(MFG-a)}}} and a given $T$, the system of FBSDEs $(\ref{mfg-fbsde})$ 
(and hence the matching problem $(\ref{mfg-fixed-point})$) is solvable for any $\xi^1,\xi^2\in \mbb{L}^2(\Omega,\calf_0,\mbb{P};\mbb{R}^d)$ 
if $\lambda^{-1}||b_{i,2}||_{\infty}$, $1\leq i\leq 2$ are
small enough.
\begin{proof}

By the growth estimate for $\hat{\alpha}_i$ in $(\ref{alpha-growth})$, it is straightforward to check that
the term involving  $\sup_{s\in[0,T]}M_2(\call(X_s^{j,n}))^2$
in $(\ref{xsq-ineq-mfg})$ is proportional to $\lambda^{-1}||b_{i,2}||_{\infty}$. Thus, by making $\lambda^{-1}||b_{i,2}||_{\infty}$ small enough
for a given $T$, we obtain the same estimate $(\ref{mfg-uniform-bound})$. 
The remaining procedures for the proof are exactly the same as in Theorem~\ref{th-mfg-main}.
\end{proof}
\end{theorem}

\begin{remark}
As one can see,  there is no difficulty to generalize all the analyses in Section~\ref{sec-mfg} for  
any finite number of populations $1\leq i\leq m$.
It results in a search for  a fixed point in the map $\calc([0,T];\calp_2(\mbb{R}^d))^m\ni (\bg{\mu}^i)_{i=1}^m\mapsto 
(\call(X_t^i)_{t\in[0,T]})_{i=1}^m\in \calc([0,T];\calp_2(\mbb{R}^d))^m$,
which can be done in the same way.
\end{remark}

\section{Games among Cooperative Populations}
\label{sec-mftc}
In this section, we try to establish the existence of a mean-field equilibrium between two competing populations
within each of which the agents share the same cost functions as well as the 
coefficient functions of the state dynamics. 
The difference from the situation studied in Section~\ref{sec-mfg} is that 
the agents within each population now cooperate by using the common feedback strategy, say,  under 
the command of  a central planner. This results in a control problem of  McKean-Vlasov type in the large population limit. 
See Section~\ref{sec-mftc-fa} (and also Chapter 6 in \cite{Carmona-Delarue-1}) to understand
the details how the large  population limit of cooperative agents induces a control problem of MKV type.
The current problem has been discussed in Section 3 in Bensoussan et.al.\cite{Bensoussan-mfg-paper}
under the name of {\it Nash Mean Field Type Control Problem},  
where the necessary conditions of the optimality are provided in the form of a master equation.
In this section, we adopt the probabilistic approach developed in Carmona \& Delarue (2015) \cite{Carmona-Delarue-MFTC}, 
 and then provide several sets of sufficient conditions for the 
existence of an equilibrium.

\subsection{Definition of Nash Mean Field Type Control Problem}
Let us first formulate the problem to be studied in this section. \\
(i) Fix any two deterministic flows of probability measures $(\bg{\mu}^i=(\mu^i_t)_{t\in[0,T]})_{i\in\{1,2\}}$ given on $\mbb{R}^d$. \\
(ii) Solve the two optimal control problems of McKean-Vlasov type
\bea
\inf_{\bg{\alpha}^1\in \mbb{A}_1}J_1^{\bg{\mu}^2}(\bg{\alpha}^1), \quad
\inf_{\bg{\alpha}^2\in \mbb{A}_2}J_2^{\bg{\mu}^1}(\bg{\alpha}^2) 
\label{mftc-problem}
\eea
over some admissible strategies $\mbb{A}_i~(i\in\{1,2\})$, where
\bea
&&J_1^{\bg{\mu}^2}(\bg{\alpha}^1):=\mbb{E}\Bigl[\int_0^T f_1(t,X_t^1,\call(X_t^1),\mu_t^2,\alpha_t^1)dt+g_1(X_T^1, \call(X_T^1),
\mu_T^2)\Bigr]~, \nn \\
&&J_2^{\bg{\mu}^1}(\bg{\alpha}^2):=\mbb{E}\Bigl[\int_0^T f_2(t,X_t^2,\call(X_t^2),\mu_t^1,\alpha_t^2)dt+g_2(X_T^2,\call(X_T^2),\mu_T^1)\Bigr]~,\nn
\eea
subject to the $d$-dimensional diffusion dynamics of McKean-Vlasov type:
\bea
&&dX_t^1=b_1(t,X_t^1,\call(X_t^1),\mu_t^2,\alpha_t^1)dt+\sigma_1(t,X_t^1,\call(X_t^1),\mu_t^2)dW_t^1~,\nn \\
&&dX_t^2=b_2(t,X_t^2,\call(X_t^2),\mu_t^1,\alpha_t^2)dt+\sigma_2(t,X_t^2,\call(X_t^2),\mu_t^1)dW_t^2~, \nn
\eea
for $t\in[0,T]$ with $\bigl(X_0^i=\xi^i\in \mbb{L}^2(\Omega,\calf_0,\mbb{P};\mbb{R}^d)\bigr)_{1\leq i\leq 2}$.
For each population $i\in\{1,2\}$, we suppose, as before, that $\mbb{A}_i$ is the set of $A_i$-valued $\mbb{F}^i$-progressively 
measurable processes $\bg{\alpha}^i$ satisfying $\mbb{E}\int_0^T |\alpha_t^i|^2 dt<\infty$ and
$A_i\subset \mbb{R}^k$ is closed and convex.  \\
(iii) Find a pair of probability flows $(\bg{\mu}^1, \bg{\mu}^2)$ as a solution to the matching problem:
\bea
\forall t\in[0,T], \quad \mu_t^1=\call(\hat{X}_t^{1,\bg{\mu}^2}), \quad
\mu_t^2=\call(\hat{X}_t^{2,\bg{\mu}^1})~, 
\label{mftc-matching}
\eea
where $(\hat{X}^{i,\bg{\mu}^j})_{i\in\{1,2\}, j\neq i}$ are the solutions to the optimal control problems in (ii).

\subsection{Optimization for given flows of probability measures}
In this subsection, we consider the step (ii) in the above formulation.
Before giving the set of main assumptions, let us mention the notion of differentiability for functions
defined on the space of probability measures.
We adopt the notion of L-differentiability used in \cite{Carmona-Delarue-MFTC}, which was first introduced by Lions in his lecture
at the {\it College de France} (see the lecture notes summarized in \cite{Cardaliaguet-note}),
where the differentiation is based on the {\it lifting} of functions $\calp_2(\mbb{R}^d)\ni \mu\rightarrow u(\mu)$
to functions $\wt{u}$ defined on a Hilbert space $\mbb{L}^2(\Omega,\calf,\mbb{P};\mbb{R}^d)$ by 
$\wt{u}(X):=u(\call(X))$ with $X\in \mbb{L}^2(\Omega, \calf, \mbb{P};\mbb{R}^d)$
over some probability space $(\Omega,\calf,\mbb{P})$ with $\Omega$ being a Polish space and
$\mbb{P}$ an atomless probability measure.

\begin{definition}(Definition 5.22 in \cite{Carmona-Delarue-1})
\label{def-LD}
A function u on $\calp_2(\mbb{R}^d)$ is said to be L-differentiable at $\mu_0\in \calp_2(\mbb{R}^d)$
if there exists a random variable $X_0$ with law $\mu_0$ such that the lifted function $\wt{u}$ is Frechet differentiable
at $X_0$.
\end{definition}
By Proposition 5.24~\cite{Carmona-Delarue-1},  if $u$ is L-differentiable at $\mu_0$ in the sense of 
Definition~\ref{def-LD}, then $\wt{u}$ is differentiable at any $X_0^\prime$ with $\call(X_0^\prime)=\mu_0$
and the law of the pair $(X_0^\prime, D\wt{u}(X_0^\prime))$ is independent of the choice of the random variable $X_0^\prime$. 
Thus the L-derivative may be denoted by $\part_\mu u(\mu_0)(\cdot):\mbb{R}^d\ni x\mapsto \part_\mu u(\mu_0)(x)\in \mbb{R}^d$,
which is uniquely defined $\mu_0$-almost everywhere on $\mbb{R}^d$.  It satisfies, according to the definition, that:
\bea
u(\mu)=u(\mu_0)+\mbb{E}\bigl[\langle X-X_0, \part_\mu u(\call(X_0))(X_0)\rangle \bigr]+o(||X-X_0||_2)~,\nn
\eea
whenever the random variables $X$ and $X_0$ have the distributions
$\call(X)=\mu, ~\call(X_0)=\mu_0$. 
For example, if the function $u$ is of the form
$u(\mu):=\int_{\mbb{R}^d}h(x)\mu(dx)$
for some function $h:\mbb{R}^d\rightarrow \mbb{R}$, we have 
$\wt{u}(X)=\mbb{E}[h(X)]$ by a random varibale $X$ with $\call(X)=\mu$. If the function $h$ is differentiable, 
the definition of L-derivative implies that $\part_\mu u(\mu)(\cdot)= \part_x h(\cdot)$.
For details of L-derivatives, their regularity properties and examples, see Section 6 in \cite{Cardaliaguet-note} and 
Chapter 5 in \cite{Carmona-Delarue-1}.
We now give the main assumptions in this section:\footnote{We  slightly abuse the notation to lighten the expression.}

\begin{assumption}{\rm{\tbf{(MFTC-a)}}} $L, K\geq 0$ and $\lambda>0$ are some constants.
For $1\leq i\leq 2$, 
the measurable functions $b_i:[0,T]\times \mbb{R}^d\times \calp_2(\mbb{R}^d)^2\times A_i\rightarrow \mbb{R}^d$,
$\sigma_i:[0,T]\times \mbb{R}^d\times \calp_2(\mbb{R}^d)^2\rightarrow \mbb{R}^{d\times d}$,
$f_i:[0,T]\times \mbb{R}^d\times \calp_2(\mbb{R}^d)^2\times A_i\rightarrow \mbb{R}$,
and $g_i:\mbb{R}^d\times \calp_2(\mbb{R}^d)^2\rightarrow \mbb{R}$ satisfy the following conditions: \\
{\rm{\tbf{(A1)}}} The functions $b_i$ and $\sigma_i$ are affine in $(x,\alpha,\bar{\mu})$ in the sense that, 
for any $(t,x,\mu,\nu,\alpha)\in [0,T]\times\mbb{R}^d\times \calp_2(\mbb{R}^d)^2\times A_i$,
\bea
&&b_i(t,x,\mu,\nu,\alpha):=b_{i,0}(t,\nu)+b_{i,1}(t,\nu)x+\bar{b}_{i,1}(t,\nu)\bar{\mu}+b_{i,2}(t,\nu)\alpha~, \nn \\
&&\sigma_i(t,x,\mu,\nu):=\sigma_{i,0}(t,\nu)+\sigma_{i,1}(t,\nu)x+\bar{\sigma}_{i,1}(t,\nu)\bar{\mu}~,\nn
\eea
where  $\bar{\mu}:=\int_{\mbb{R}^d}x\mu(dx)$,  and $b_{i,0}, b_{i,1}, \bar{b}_{i,1}, b_{i,2}, \sigma_{i,0}, \sigma_{i,1}$ and $\bar{\sigma}_{i,1}$ defined on 
$[0,T]\times\calp_2(\mbb{R}^d)$ are $\mbb{R}^d$, $\mbb{R}^{d\times d}$, $\mbb{R}^{d\times d}$, $\mbb{R}^{d\times k}$,
$\mbb{R}^{d\times d}$, $\mbb{R}^{d\times d\times d}$ and $\mbb{R}^{d\times d\times d}$-valued measurable functions, respectively.\\
{\rm{\tbf{(A2)}}}
For any $t\in[0,T]$, the functions $\calp_2(\mbb{R}^d)\ni \nu\mapsto (b_{i,0},b_{i,1},\bar{b}_{i,1},b_{i,2},\sigma_{i,0},\sigma_{i,1}, 
\bar{\sigma}_{i,1})(t,\nu)$ are continuous in $W_2$-distance. Moreover for any $(t,\nu)\in[0,T]\times \calp_2(\mbb{R}^d)$, 
\bea
&&|b_{i,0}(t,\nu)|, |\sigma_{i,0}(t,\nu)|\leq K+LM_2(\nu)~, \nn \\
&&|b_{i,1}(t,\nu)|, |\bar{b}_{i,1}(t,\nu)|, |b_{i,2}(t,\nu)|, |\sigma_{i,1}(t,\nu)|, |\bar{\sigma}_{i,1}(t,\nu)|\leq L~.\nn
\eea
{\rm{\tbf{(A3)}}} For any $t\in[0,T],x,x^\prime\in \mbb{R}^d, \mu,\mu^\prime, \nu\in \calp_2(\mbb{R}^d)$ and $\alpha, \alpha^\prime \in A_i$,
the functions $f_i$ and $g_i$ satisfy the quadratic growth conditions
\bea
&&|f_i(t,x,\mu, \nu,\alpha)|\leq K+L\bigl(|x|^2+|\alpha|^2+M_2(\mu)^2+M_2(\nu)^2\bigr)~, \nn \\
&&|g_i(x,\mu,\nu)|\leq K+L\bigl(|x|^2+M_2(\mu)^2+M_2(\nu)^2\bigr)~, \nn
\eea
and the local Lipschitz continuity
\bea
&&|f_i(t,x^\prime,\mu^\prime,\nu,\alpha^\prime)-f_i(t,x,\mu,\nu,\alpha)|+
|g_i(x^\prime,\mu^\prime,\nu)-g_i(x,\mu,\nu)| \nn \\
&&\leq\Bigl(K+L\bigl[|(x^\prime,\alpha^\prime)|+|(x,\alpha)|
+M_2(\mu^\prime)+M_2(\mu)+M_2(\nu)]\Bigr) \bigl[|(x^\prime,\alpha^\prime)-(x,\alpha)|+W_2(\mu^\prime,\mu)\bigr]~.\nn
\eea
{\rm{\tbf{(A4)}}} The functions $f_i$ and $g_i$ are once continuously differentiable in $(x,\alpha)$ and $x$ respectively,
and their derivatives are $L$-Lipschitz continuous with respect to $(x,\alpha,\mu)$ and $(x,\mu)$ i.e.
\bea
&&|\part_{(x,\alpha)}f_i(t,x^\prime,\mu^\prime,\nu,\alpha^\prime)-\part_{(x,\alpha)}f_i(t,x,\mu,\nu,\alpha)|+
|\part_x g_i(x^\prime,\mu^\prime,\nu)-\part_x g_i(x,\mu,\nu)|\nn \\
&&\leq L\bigl(|x^\prime-x|+|\alpha^\prime-\alpha|+W_2(\mu^\prime,\mu)\bigr)~, \nn
\eea
for any  $t\in[0,T],x,x^\prime\in \mbb{R}^d, \mu,\mu^\prime, \nu\in \calp_2(\mbb{R}^d), 
\alpha,\alpha^\prime \in A_i$.  The derivatives also satisfy the growth condition
\bea
|\part_{(x,\alpha)}f_i(t,x,\mu,\nu,\alpha)|+|\part_x g_i(x,\mu,\nu)|\leq K+L\bigl(|x|+|\alpha|+M_2(\mu)+M_2(\nu)\bigr)~.\nn
\eea
Moreover, the derivatives $\part_{(x,\alpha)}f_i$ and $\part_x g_i$ are continuous also in $\nu$ with respect to the $W_2$-distance.
\\
{\rm{\tbf{(A5)}}} The functions $f_i$ and $g_i$ are L-differentiable with respect to the 
first measure argument $\mu$ and they satisfy that, 
for any $t\in[0,T],x,x^\prime\in\mbb{R}^d,\mu,\mu^\prime,\nu\in\calp_2(\mbb{R}^d), \alpha,\alpha^\prime \in A_i$
and any random variables $X, X^\prime$  with $\call(X)=\mu$, $\call(X^\prime)=\mu^\prime$, 
$L$-Lipschitz continuity in $\mbb{L}^2$ i.e.,
\bea
&&||\part_\mu f_i(t,x^\prime,\mu^\prime,\nu,\alpha^\prime)(X^\prime)
-\part_\mu f_i(t,x,\mu,\nu,\alpha)(X)||_2+
||\part_\mu g_i(x^\prime,\mu^\prime,\nu)(X^\prime)-\part_\mu g_i(x,\mu,\nu)(X)||_2\nn \\
&&~\leq L\bigl(|x^\prime-x|+|\alpha^\prime-\alpha|+||X^\prime-X||_2\bigr)~,\nn
\eea
as well as the following growth condition:
\bea
||\part_\mu f_i(t,x,\mu,\nu,\alpha)(X)||_2+||\part_\mu g_i(x,\mu,\nu)(X)||_2\leq K+L\bigl(|x|+|\alpha|+M_2(\mu)+M_2(\nu)\bigr)~.\nn
\eea
Moreover, the maps $\calp_2(\mbb{R}^d)\ni \nu\mapsto \part_\mu f_i (t,x,\mu,\nu,\alpha)(v)$ and
$\calp_2(\mbb{R}^d)\ni \nu\mapsto \part_\mu g_i (x,\mu,\nu)(v)$ are  continuous with respect to  $W_2$-distance
for any $(t,x,\mu,\alpha)\in[0,T]\times \mbb{R}^d\times \calp_2(\mbb{R}^d)\times A_i$, and $\mu$-a.e. $v\in \mbb{R}^d$. \\
{\rm{\tbf{(A6)}}} For any $t\in[0,T], x, x^\prime\in \mbb{R}^d, \mu, \mu^\prime,\nu\in\calp_2(\mbb{R}^d), \alpha, \alpha^\prime\in A_i$,
and any random variables $X, X^\prime$ with $ \call(X)=\mu, \call(X^\prime)=\mu^\prime$, the functions $f_i$ and $g_i$ satisfy the convexity relations:
\bea
&&f_i(t,x^\prime,\mu^\prime,\nu,\alpha^\prime)-f_i(t,x,\mu,\nu,\alpha)-\langle (x^\prime-x, \alpha^\prime-\alpha),
\part_{(x,\alpha)}f_i(t,x,\mu,\nu,\alpha)\rangle\nn \\
&&\qquad -\mbb{E}\bigl[\langle X^\prime-X, \part_\mu f_i(t,x,\mu,\nu,\alpha)(X)\rangle]\geq \lambda|\alpha^\prime-\alpha|^2~, \nn \\
&&g_i(x^\prime,\mu^\prime,\nu)-g_i(x,\mu,\nu)-\langle x^\prime-x, \part_x g_i(x,\mu,\nu)\rangle
-\mbb{E}\bigl[\langle X^\prime-X, \part_\mu g_i(x,\mu,\nu)(X)\rangle\bigr]\geq 0~.\nn
\eea
\end{assumption}

\begin{remark}
By Lemma 3.3 in \cite{Carmona-Delarue-MFTC}, the Lipschitz continuity in {\rm{\tbf{(A5)}}} above implies that
we can modify $\part_\mu f_i(t,x,\mu,\nu,\alpha)(\cdot)$ 
and $\part_\mu g_i(x,\mu,\nu)(\cdot)$ on a $\mu$-negligible set in 
such a way that, $\forall v,v^\prime\in \mbb{R}^d$
\bea
&&|\part_\mu f_i(t,x,\mu,\nu,\alpha)(v^\prime)-\part_\mu f_i(t,x,\mu,\nu,\alpha)(v)| \leq L|v^\prime-v|~, \nn \\
&&|\part_\mu g_i(x,\mu,\nu)(v^\prime)-\part_\mu g_i(x,\mu,\nu)(v)|\leq L|v^\prime-v|~, \nn
\eea
for any $(t,x,\mu,\nu,\alpha)\in [0,T]\times \mbb{R}^d\times \calp_2(\mbb{R}^d)^2\times A_i$.
In the remainder of the work,  we always use these Lipschitz continuous versions.
\end{remark}

As before, we first consider the optimal control problem $(\ref{mftc-problem})$ for given 
deterministic flows of probability measures. The Hamiltonian for each population $H_i:[0,T]\times \mbb{R}^d\times 
\calp_2(\mbb{R}^d)^2\times \mbb{R}^d\times \mbb{R}^{d\times d}\times A_i
\ni (t,x,\mu,\nu,y,z,\alpha)\mapsto H_i(t,x,\mu,\nu,y,z,\alpha)\in \mbb{R}$ and  
its minimizer $\hat{\alpha}_i:[0,T]\times \mbb{R}^d\times \calp_2(\mbb{R}^d)^2\times \mbb{R}^d\ni (t,x,\mu,\nu,y)\mapsto 
\hat{\alpha}_i(t,x,\mu,\nu,y)\in A_i$ are defined in the same way as $(\ref{def-H})$ and $(\ref{def-alpha})$ 
with the coefficients replaced by those given in the current section.

\begin{lemma}
\label{lemma-alpha-mftc}
Under Assumption {\rm{\tbf{(MFTC-a)}}}, for all $(t,x,\mu,\nu,y)\in[0,T]\times \mbb{R}^d\times \calp_2(\mbb{R}^d)^2\times \mbb{R}^d$,
there exists a unique minimizer $\hat{\alpha}_i(t,x,\mu,\nu,y)$ of $H_i^{(r)}$,
where the map $[0,T]\times \mbb{R}^d\times \calp_2(\mbb{R}^d)^2\times \mbb{R}^d\ni (t,x,\mu,\nu,y)\mapsto \hat{\alpha}_i(t,x,\mu,\nu,y)\in A_i$
is measurable. There exist constants $C$ depending only on $(L,\lambda)$ and $C^\prime$ depending additionally on $K$ such that,
for any $t\in[0,T], x,x^\prime, y, y^\prime \in \mbb{R}^d, \mu,\nu \in\calp_2(\mbb{R}^d)$, 
\bea
&&|\hat{\alpha}_i(t,x,\mu,\nu,y)|\leq C^\prime+C\bigl(|x|+|y|+M_2(\mu)+M_2(\nu)\bigr)\nn \\
&&|\hat{\alpha}_i(t,x,\mu,\nu,y)-\hat{\alpha}_i(t,x^\prime,\mu,\nu,y^\prime)|\leq C\bigl(|x-x^\prime|+|y-y^\prime|\bigr)~. \nn 
\eea
Moreover, for any $(t,x,y)\in[0,T]\times \mbb{R}^d\times \mbb{R}^d$, the map
$\calp_2(\mbb{R}^d)^2\ni (\mu, \nu)\mapsto \hat{\alpha}_i(t,x,\mu,\nu,y)$ is continuous with respect to 
$W_2$-distance:
\bea 
&&|\hat{\alpha}_i(t,x,\mu,\nu,y)-\hat{\alpha}_i(t,x,\mu^\prime, \nu^\prime, y)| \nn \\
&&\leq (2\lambda)^{-1}\Bigl(LW_2(\mu,\mu^\prime)+|b_{i,2}(t,\nu)-b_{i,2}(t, \nu^\prime)||y|+|\part_\alpha f_i(t,x,\mu^\prime,\nu,\hat{\alpha}_i)
-\part_\alpha f_i(t,x,\mu^\prime, \nu^\prime,\hat{\alpha}_i)|\Bigr)\nn
\eea
where $\hat{\alpha}_i:=\hat{\alpha}_i(t,x,\mu,\nu,y)$.
\begin{proof}
It can be shown exactly in the same way as Lemma~\ref{lemma-alpha}.
\end{proof}
\end{lemma}

The control problem $(\ref{mftc-problem})$ 
for each population $1\leq i\leq 2$ with a given flow of probability measure $\bg{\mu}^j\in\calc([0,T];\calp_2(\mbb{R}^d)), j\neq i$ 
is actually the special case studied in \cite{Carmona-Delarue-MFTC} and Section 6.4 in \cite{Carmona-Delarue-1}.
In fact, we have removed the control $\bg{\alpha}_i$ dependency from the diffusion coefficient $\sigma_i$.\footnote{It then makes possible to derive
the stability relation irrespective of the size of Lipschitz constant for $Z$.}
The relevant adjoint equations for the optimal control problem of MKV-type $(\ref{mftc-problem})$ are given by
with $i,j\in \{1,2\}, j\neq i$:
\bea
\label{adjoint-mftc}
&&dX_t^i=b_i(t,X_t^i,\call(X_t^i),\mu_t^j,\hat{\alpha}_i(t,X_t^i,\call(X_t^i),\mu_t^j,Y_t^i))dt+\sigma_i(t,X_t^i, \call(X_t^i),\mu^j_t)dW_t^i~, \nn \\
&&dY_t^i=-\part_x H_i(t,X_t^i,\call(X_t^i),\mu_t^j, Y_t^i,Z_t^i,\hat{\alpha}_i(t,X_t^i,\call(X_t^i),\mu_t^j,Y_t^i))dt \\
&&\quad\qquad -\wt{\mbb{E}}\bigl[\part_\mu H_i(t,\wt{X}_t^i,\call(X_t^i),\mu_t^j,\wt{Y}_t^i, \wt{Z}_t^i, \hat{\alpha}_i
(t,\wt{X}_t^i,\call(X_t^i),\mu_t^j, \wt{Y}_t^i))(X_t^i)\bigr]dt+Z_t^i dW_t^i~, \nn
\eea
with $X_0^i=\xi^i\in \mbb{L}^2(\Omega,\calf_0,\mbb{P};\mbb{R}^d)$ and 
$Y_T^i=\part_x g_i(X_T^i,\call(X_T^i),\mu_T^j)+\wt{\mbb{E}}\bigl[\part_\mu g_i(\wt{X}_T^i,\call(X_T^i),\mu_T^j)(X_T^i)\bigr]$.
Here, $(\wt{\Omega},\wt{\calf}, \wt{\mbb{P}})$ denotes a copy of $(\Omega,\calf,\mbb{P})$ and every random variable 
with tilde, such as $\wt{X}$, denotes a clone of $X$ on $(\wt{\Omega},\wt{\calf},\wt{\mbb{P}})$.
The expectation under $\wt{\mbb{P}}$ is denoted by $\wt{\mbb{E}}$.
More explicitly, one can write $(\ref{adjoint-mftc})$ as 
\bea
&&dX_t^i=\bigl(b_{i,0}(t,\mu_t^j)+b_{i,1}(t,\mu_t^j)X_t^i+\bar{b}_{i,1}(t,\mu_t^j)\mbb{E}[X_t^i]+
b_{i,2}(t,\mu_t^j)\hat{\alpha}_i(t,X_t^i,\call(X_t^i),\mu_t^j,Y_t^i)\bigr)dt\nn \\
&&\quad+\bigl(\sigma_{i,0}(t,\mu_t^j)+\sigma_{i,1}(t,\mu_t^j)X_t^i+\bar{\sigma}_{i,1}(t,\mu_t^j)\mbb{E}[X_t^i]\bigr)dW_t^i~,\nn \\
&&dY_t^i=-\bigl(b_{i,1}(t,\mu_t^j)^\top Y_t^i+\sigma_{i,1}(t,\mu_t^j)^\top Z_t^i+\part_x f_i(t,X_t^i,\call(X_t^i),\mu_t^j,
\hat{\alpha}_i(t,X_t^i,\call(X_t^i),\mu_t^j,Y_t^i))\bigr)dt\nn\\
&&\quad-\bigl(\bar{b}_{i,1}(t,\mu_t^j)^\top\mbb{E}[Y_t^i]+\bar{\sigma}_{i,1}(t,\mu_t^j)^\top \mbb{E}[Z_t^i]+
\wt{\mbb{E}}\bigl[\part_\mu f_i(t,\wt{X}_t^i,\call(X_t^i),\mu_t^j,\hat{\alpha}_i(t,\wt{X}_t^i,\call(X_t^i),\mu_t^j,
\wt{Y}_t^i))(X_t^i)\bigr]dt\nn\\
&&\quad +Z_t^i dW_t^i~,\nn
\eea
which is a $C(L,\lambda)$-Lipschitz FBSDE of McKean-Vlasov type. Note that due to Lemma~\ref{lemma-alpha-mftc}, $\hat{\alpha}_i$
is Lipschitz continuous not only in $(X^i,Y^i)$ but also in $\call(X^i)$.
For each $i\in\{1,2\}$, it is important to notice that the Lipschitz constant is independent of the
given flow $\bg{\mu}^j, j\neq i$. 
Since Assumption {\rm{\tbf{(MFTC-a)}}} satisfies every solvability condition 
used in \cite{Carmona-Delarue-MFTC}, we have the following results:\footnote{In \cite{Carmona-Delarue-MFTC}, $A=\mbb{R}^k$
is assumed. However, there is no difficulty for extending a general closed and convex subset $A\subset \mbb{R}^k$, which is
actually the case studied  in Chapter 6 in \cite{Carmona-Delarue-1}.} 
\begin{theorem}
\label{th-mftc-existence}
Under Assumption {\rm{\tbf{(MFTC-a)}}},
the adjoint FBSDE $(\ref{adjoint-mftc})$ of each $i\in\{1,2\}$ 
has a unique solution $(\hat{X}_t^i,\hat{Y}_t^i,\hat{Z}_t^i)_{t\in[0,T]}\in \mbb{S}^2\times \mbb{S}^2\times \mbb{H}^2$
for any flow $\bg{\mu}^j \in\calc([0,T];\calp_2(\mbb{R}^d))$ and  any initial condition $\xi^i\in \mbb{L}^2(\Omega, \calf_0, \mbb{P};\mbb{R}^d)$.
If we set $\hat{\bg{\alpha}}_i=\bigl(\hat{\alpha}^i_t=\hat{\alpha}_i(t,\hat{X}_t^i,\call(X_t^i),\mu_t^j,\hat{Y}_t^i)\bigr)_{t\in[0,T]}$,
then it gives the optimal control. In particular, the inequality  $J_i^{\bg{\mu}^j}(\hat{\bg{\alpha}}^i)+\lambda \mbb{E}\int_0^T |\beta_t^i-\hat{\alpha}^i_t|^2 dt \leq J_i^{\bg{\mu}^j}(\bg{\beta}^i)$ holds
for any $\bg{\beta}^i\in \mbb{A}_i$.
\begin{proof}
This is the direct result of Theorem 4.7 (sufficiency) and 
Theorem 5.1 (unique solvability) by Carmona \& Delarue (2015)\cite{Carmona-Delarue-MFTC}, 
where the sufficiency is proved in a parallel way to Theorem~\ref{th-fbsde-existence}, and 
the unique solvability is based on  {\it the continuation method} developed by Peng \& Wu (1999)~\cite{Peng-Wu}.
\end{proof}
\end{theorem}

\begin{lemma}
\label{lemma-Lip-mftc}
Under the same conditions used in Theorem~\ref{th-mftc-existence}, 
for any $t\in[0,T]$ and any $\xi^i\in\mbb{L}^2(\Omega,\calf_t,\mbb{P};\mbb{R}^d)$, 
there exists a unique solution, denoted by $(X_s^{i,t,\xi^i},Y_s^{i,t,\xi^i}, Z_s^{i,t,\xi^i})_{t\leq s\leq T}$, of $(\ref{adjoint-mftc})$
on $[t,T]$ with $X_t^{i,t,\xi^i}=\xi^i$ as initial condition. Moreover,  for any $\mu\in\calp_2(\mbb{R}^d)$, 
there exists a measurable mapping $u_i^{\bg{\mu}^j}(t,\cdot,\mu):\mbb{R}^d\ni x\mapsto u_i^{\bg{\mu}^j}(t,x,\mu)$ with
$Y_t^{i,t,\xi^i}=u_i^{\bg{\mu}^j}(t,\xi^i,\call(\xi^i))$ $\mbb{P}$-a.s. such that, for any $\xi^i,\xi^{\prime,i}\in\mbb{L}^2(\Omega,\calf_t,\mbb{P};\mbb{R}^d)$,
\bea
\mbb{E}\bigl[|u_i^{\bg{\mu}^j}(t,\xi^i,\call(\xi^i))-u_i^{\bg{\mu}^j}(t,\xi^{\prime,i},\call(\xi^{\prime,i}))|^2\bigr]^\frac{1}{2}\leq C\mbb{E}\bigl[|\xi^i-\xi^{\prime,i}|^2\bigr]^\frac{1}{2}~, 
\label{master-lip}
\eea
with some constant $C$ depending only on $L$ and $\lambda$.
\begin{proof}
This is a direct result of Lemma 5.6 in \cite{Carmona-Delarue-MFTC}. The Lipschitz constant can be read from the stability estimate
used in the continuation method (Lemma 5.5 in \cite{Carmona-Delarue-MFTC}), which is dependent only on the Lipschitz constant of the FBSDE.
\end{proof}
\end{lemma}
\begin{remark}
\label{remark-Lip-mftc}.
Note that, due to the uniqueness of the solution,  we have for any $t\in[0,T]$, 
$\hat{Y}_t^{i}=Y_t^{i,0,\xi^i}=Y_t^{i,t,\hat{X}_t^{i}}=u_i^{\bg{\mu}^j}(t,\hat{X}_t^i,\call(\hat{X}_t^i))$ $\mbb{P}$-a.s.
Moreover, once again by Lemma 3.3 in \cite{Carmona-Delarue-MFTC}, for any $\mu\in \calp_2(\mbb{R}^d)$, 
there exists a version $\mbb{R}^d\ni x\mapsto u_i^{\bg{\mu}^j}(t,x,\mu)$ in $\mbb{L}^2(\mbb{R}^d,\mu)$ that is 
Lipschitz continuous with the same Lipschitz constant $C$ used in $(\ref{master-lip})$ i.e.,
$|u_i^{\bg{\mu}^j}(t,x,\mu)-u_i^{\bg{\mu}^j}(t,x^\prime, \mu)|\leq C|x-x^\prime|$ for any $x,x^\prime\in \mbb{R}^d$.
In the remainder, we always use this Lipschitz version and often adopt a simpler notation $u_i$ without the superscript $\bg{\mu}^j$. 
\end{remark}
Making use of the Lipschitz continuity in Lemma~\ref{lemma-Lip-mftc},  we can derive the stability relation.

\begin{lemma}
\label{lemma-stability-mkv}
Suppose that the two set of functions $(b_i, \sigma_i, f_i, g_i)$ and $(b_i^\prime, \sigma_i^\prime, f_i^\prime,g_i^\prime)$
satisfy Assumption {\rm{\tbf{(MFTC-a)}}}. For given inputs $\xi^i, \xi^{\prime, i}\in \mbb{L}^2(\Omega, \calf_0, \mbb{P};\mbb{R}^d)$
and $\bg{\mu}^j, \bg{\mu}^{\prime, j}\in \calc([0,T];\calp_2(\mbb{R}^d))$, 
let us denote the corresponding solution to $(\ref{adjoint-mftc})$ by $(X^i_t,Y^i_t,Z^i_t)_{t\in[0,T]}$
and $(X^{\prime,i}_t, Y^{\prime, i}_t, Z^{\prime, i}_t)_{t\in[0,T]}$, respectively.
Then, there exists a constant $C$ depending only on $(L,\lambda)$ such that
\bea
\label{stability-mkv}
&&\mbb{E}\Bigl[\sup_{t\in[0,T]}|X_t^{i}-X_t^{\prime, i}|^2+\sup_{t\in[0,T]}|Y_t^i-Y_t^{\prime, i}|^2
+\int_0^T|Z_t^i-Z_t^{\prime,i}|^2 dt\Bigr]\nn \\
&&\leq C\mbb{E}\Bigl\{|\xi^i-\xi^{\prime,i}|^2+|\part_x g_i(X_T^i, \call(X_T^i), \mu_T^j)-\part_x g_i^\prime(X_T^i, \call(X_T^i),\mu_T^{\prime,j})|^2 \nn \\
&&\hspace{25mm} +\wt{\mbb{E}}\bigl[|\part_\mu g_i(\wt{X}_T^i,\call(X_T^i),\mu_T^j)(X_T^i)-\part_\mu g_i^\prime(\wt{X}_T^i,\call(X_T^i),\mu_T^{j,\prime})(X_T^i)|^2\bigr]\nn \\
&&+\int_0^T|b_i(t,X_t^i, \call(X_t^i), \mu_t^j, \hat{\alpha}_i(t,X_t^i, \call(X_t^i),\mu_t^j,Y_t^i))
-b_i^\prime (t,X_t^i, \call(X_t^i), \mu_t^{\prime, j}, \hat{\alpha}_i^\prime (t,X_t^i,\call(X_t^i),\mu_t^{\prime,j},Y_t^i))|^2dt \nn\\
&&+\int_0^T |\sigma_i(t,X_t^i,\call(X_t^i),\mu_t^j)-\sigma_{i}^\prime(t,X_t^i,\call(X_t^i), \mu_t^{\prime,j})|^2 dt
\nn \\
&&+\int_0^T \Bigl(|\part_x H_i(t,X_t^i, \call(X_t^i), \mu_t^j, Y_t^i, Z_t^i, \hat{\alpha}_i(t,X_t^i, \call(X_t^i), \mu_t^j, Y_t^i)) \nn \\
&&\hspace{25mm}-\part_x H_i^\prime (t,X_t^i, \call(X_t^i),\mu_t^{\prime, j}, Y_t^i, Z_t^i, \hat{\alpha}_i^\prime (t,X_t^i,\call(X_t^i),
\mu_t^{\prime,j}, Y_t^i))|^2\Bigr) dt\nn 
\eea
\bea
&&+\int_0^T \wt{\mbb{E}}\Bigl[|\part_\mu H_i(t,\wt{X}_t^i,\call(X_t^i),\mu_t^j,\wt{Y}_t^i,\wt{Z}_t^i,
\hat{\alpha}_i(t,\wt{X}_t^i,\call(X_t^i),\mu_t^j,\wt{Y}_t^i))(X_t^i)\nn\\
&&\hspace{25mm}-\part_\mu H_i^\prime(t,\wt{X}_t^i,\call(X_t^i),\mu_t^{\prime,j},\wt{Y}_t^i, \wt{Z}_t^i,
\hat{\alpha}_i^\prime(t,\wt{X}_t^i,\call(X_t^i),\mu_t^{\prime,j},\wt{Y}_t^i))(X_t^i)|^2\Bigr]dt\Bigr\}~,
\eea
where the functions $H_i, H_i^\prime:[0,T]\times \mbb{R}^d\times\calp_2(\mbb{R}^d)^2\times \mbb{R}^d\times \mbb{R}^{d\times d}\times A_i
\rightarrow \mbb{R}$
are the Hamiltonians associated with the coefficients $(b_i, \sigma_i, f_i)$ and $(b_i^\prime, \sigma_i^\prime, f_i^\prime)$, respectively, 
and $\hat{\alpha}_i, \hat{\alpha}_i^\prime$ are their minimizers.

In particular, there is another constant $C^\prime$ depending additionally on $K$ such that
\bea
\mbb{E}\Bigl[\sup_{t\in[0,T]}|X_t^i|^2+\sup_{t\in[0,T]}|Y_t^i|^2+\int_0^T |Z_t^i|^2 dt\Bigr]\leq C\Bigl(||\xi^i||_2^2+\sup_{t\in[0,T]}M_2(\mu^j_t)^2\Bigr)
+C^\prime,
\label{mkv-growth}
\eea
and, for any $t\in[0,T]$,
\bea
|Y_t^i|\leq C\Bigl(||\xi^i||_2+|X_t^i|+\sup_{s\in[0,T]}M_2(\mu_s^j)\Bigr)+C^\prime,~\mbb{P}\text{-a.s.}
\label{mkv-y-growth}
\eea
\begin{proof}
It can be proved  in the same way as Lemma~\ref{lemma-stability}.
For small $T\leq c$, where $c$ is dependent only on $(L,\lambda)$, using the inequality $W_2(\call(X),\call(Y))^2\leq \mbb{E}|X-Y|^2$,
one can show the stability relation $(\ref{stability-mkv})$ exactly in the same way as in 
the standard Lipschitz FBSDE of non-MKV type.
For general $T$, we can connect the short-term estimate by the same technique adopted in 
the proof of Lemma~\ref{lemma-stability}. Here, we make use of  the Lipschitz continuity in Lemma~\ref{lemma-Lip-mftc}.

As for the growth conditions, we get, by the same arguments used  to derive $(\ref{mfg-gw-estimate})$,
\bea
&&\mbb{E}\Bigl[\sup_{t\in[0,T]}|X_t^{i}|^2+\sup_{t\in[0,T]}|Y_t^{i}|^2+\int_0^T|Z_t^i|^2dt\Bigr]\nn \\
&&\leq C\Bigl(||\xi^i||_2^2+|\part_x g_i(0,\del_0,\mu_T^j)|^2+|\part_\mu g_i (0,\del_0,\mu_T^j)(0)|^2 \nn \\
&& +\int_0^T \bigl(|b_i(t,0,\del_0,\mu_t^j,\hat{\alpha}_i(t,0,\del_0,\mu_t^j,0))|^2+
|\sigma_i(t,0,\del_0,\mu_t^j)|^2 \bigr)dt\nn \\
&& +\int_0^T \bigl(|\part_x f_i(t,0,\del_0,\mu_t^j,\hat{\alpha}_i(t,0,\del_0,\mu_t^j,0))|^2+
|\part_\mu f_i(t,0,\del_0,\mu_t^j,\hat{\alpha}_i(t,0,\del_0,\mu_t^j,0))(0)|^2 \bigr)dt\Bigr),
\label{mkv-growth-middle}
\eea
where $\del_0$ denotes the distribution with Dirac mass at the origin.  $(\ref{mkv-growth})$ now easily follows.
Finally, since $Y_t^i=Y_t^{i,t,X_t^i}$, we have
$||u_i(t,X_t^i,\call(X_t^i))||_2\leq C\Bigl(||\xi^i||_2+\sup_{t\in[0,T]}M_2(\mu_t^j)\Bigr)+C^\prime$ from $(\ref{mkv-growth})$.
By the Lipschitz continuity in Remark~\ref{remark-Lip-mftc}  and the estimate in $(\ref{mkv-growth})$, we get
\be
|u_i(t,0,\call(X_t^i))| \leq ||u_i(t,X_t^i,\call(X_t^i))||_2+C ||X_t^i||_2\leq C\Bigl(||\xi^i||_2+\sup_{t\in[0,T]}M_2(\mu_t^j)\Bigr)+C^\prime~.\nn
\ee
Using the Lipschitz continuity in Remark~\ref{remark-Lip-mftc} once again, we get the desired estimate $(\ref{mkv-y-growth})$.
\end{proof}
\end{lemma}

\subsection{Nash MFTC equilibrium under additional boundedness}
In preceding subsections, we have seen that, for given 
flows of probability measures $\bg{\mu}^1,\bg{\mu}^2\in \calc([0,T];\calp_2(\mbb{R}^d))$,
the solution to each optimal control problem of $(\ref{mftc-problem})$ is characterized 
by the uniquely solvable FBSDE $(\ref{adjoint-mftc})$.
It follows that finding a solution to a matching problem $(\ref{mftc-matching})$ is equivalent to find a
solution to the coupled systems of  FBSDEs of MKV-type: for $i,j\in\{1,2\}, j\neq i$,
\bea
\label{mkv-mftc-master}
&&dX_t^i=b_i(t,X_t^i,\call(X_t^i),\call(X_t^j),\hat{\alpha}_i(t,X_t^i,\call(X_t^i),\call(X_t^j),Y_t^i))dt+
\sigma_i(t,X_t^i,\call(X_t^i),\call(X_t^j))dW_t^i, \nn \\
&&dY_t^i=-\part_x H_i(t,X_t^i,\call(X_t^i),\call(X_t^j),Y_t^i,Z_t^i,\hat{\alpha}_i(t,X_t^i,\call(X_t^i),\call(X_t^j),Y_t^i))dt \\
&&\qquad\quad-\wt{\mbb{E}}\bigl[\part_\mu H_i(t,\wt{X}_t^i,\call(X_t^i),\call(X_t^j),\wt{Y}_t^i,\wt{Z}_t^i,
\hat{\alpha}_i(t,\wt{X}_t^i,\call(X_t^i),\call(X_t^j),\wt{Y}_t^i))(X_t^i)\bigr]dt+Z_t^i dW_t^i, \nn
\eea
with $X_0^i=\xi^i\in\mbb{L}^2(\Omega,\calf_0,\mbb{P};\mbb{R}^d)$
and $Y_T^i=\part_x g_i(X_T^i,\call(X_T^i),\call(X_T^j))+\wt{\mbb{E}}\bigl[\part_\mu g_i(\wt{X}_T^i,\call(X_T^i),\call(X_T^j))(X_T^i)\bigr]$~.

In this subsection, we prove the existence of a solution to the system of FBSDEs $(\ref{mkv-mftc-master})$
under the additional assumption.
\begin{assumption}{\rm{\tbf{(MFTC-b)}}} 
For each $1\leq i\leq 2$, there exists some constant $\Lambda$ and some point $0_{A_i}\in A_i$ such that,
for any $t\in[0,T]$ and any $\nu\in\calp_2(\mbb{R}^d)$, 
\bea
&&|b_{i,0}(t,\nu)|, |\sigma_{i,0}(t,\nu)|\leq \Lambda,  \nn \\
&&|\part_{(x,\alpha)}f_i(t,0,\del_0,\nu,0_{A_i})|, |\part_x g_i(0,\del_0,\nu)|\leq \Lambda, \nn \\
&&|\part_\mu f_i(t,0,\del_0,\nu,0_{A_i})(0)|, |\part_\mu g_i(0,\del_0, \nu)(0)|\leq \Lambda~. \nn
\eea
\end{assumption}

Here is the main result of this subsection.
\begin{theorem}
\label{th-mftc-bounded}
Under Assumptions {\rm{\tbf{(MFTC-a,b)}}}, the system of FBSDEs $(\ref{mkv-mftc-master})$ 
(and hence the matching problem $(\ref{mftc-matching})$) is solvable for any $\xi^1,\xi^2\in \mbb{L}^2(\Omega,\calf_0,\mbb{P};\mbb{R}^d)$.
\begin{proof}
We let, with $1\leq i\leq 2$,  $(X_t^{i,\bg{\mu}^j}, Y_t^{i,\bg{\mu}^j},Z_t^{i,\bg{\mu}^j})_{t\in[0,T]}$ denote the solution to 
the FBSDE $(\ref{adjoint-mftc})$ for a given flow $\bg{\mu}^j\in \calc([0,T];\calp_2(\mbb{R}^d))$ $j\neq i$
and the initial condition $X_0^{i,\bg{\mu}^j}=\xi^i$.
By Theorem~\ref{th-mftc-existence},  we can define a map:
\bea
\Phi:\calc([0,T];\calp_2(\mbb{R}^d))^2\ni (\bg{\mu}^1,\bg{\mu}^2)\mapsto \bigl(\call(X_t^{1,\bg{\mu}^2})_{t\in[0,T]},
\call(X_t^{2,\bg{\mu}^2})_{t\in[0,T]}\bigr)\in \calc([0,T];\calp_2(\mbb{R}^d))^2~.\nn
\eea
It is easy to see that the solvability of the system of FBSDEs with McKean-Vlasov type $(\ref{mkv-mftc-master})$
is equivalent to the existence of a fixed point of the map $\Phi$.
As in Theorem~\ref{th-mfg-bounded}, we equip the linear space $\calc([0,T];\calm_f^1(\mbb{R}^d))^2$ with 
the supremum of the Kantorovich-Rubinstein norm $(\ref{KR-norm})$ so that we can apply Schauder FPT (Theorem~\ref{th-Schauder}).

We start from studying a priori estimates.
By the estimate in $(\ref{alpha-growth})$, we get
\be
|\hat{\alpha}_i(t,x,\mu,\nu,y)|\leq \lambda^{-1}\bigl(|b_{i,2}(t,\nu)||y|+|\part_\alpha f_i(t,x,\mu,\nu,0_{A_i})|\bigr)+|0_{A_i}|,\nn
\ee
and hence $|\hat{\alpha}_i(t,0,\del_0,\mu_t^j,0)|\leq \lambda^{-1} \Lambda+|0_{A_i}|\leq C(\lambda,\Lambda)$ uniformly in $\bg{\mu}^j$.
The estimate $(\ref{mkv-growth-middle})$ then implies that 
$\mbb{E}\bigl[\sup_{t\in[0,T]}|Y_t^{i,\bg{\mu}^j}|^2\bigr]\leq C\bigl(1+||\xi^i||_2^2\bigr)$ with $C$ independent of $\bg{\mu}^j$.
From the last part of the proof for Lemma~\ref{lemma-stability-mkv}, we get, for any $t\in[0,T]$, 
\be
|Y_t^{i,\bg{\mu}^j}|\leq C \bigl(1+||\xi^i||_2+|X_t^{i,\bg{\mu}^j}|\bigr), \quad \mbb{P}\text{-a.s.}\nn
\ee
and hence
$|\hat{\alpha}_i(t,X_t^{i,\bg{\mu}^j},\call(X_t^{i,\bg{\mu}^j}),\mu_t^j,Y_t^{i,\bg{\mu}^j})|\leq C\bigl(1+||\xi^i||_2+|X_t^{i,\bg{\mu}^j}|+M_2(\call(X_t^{i,\bg{\mu}^j}))\bigr)$ uniformly in $\bg{\mu}^j$. 
Thus it is straightforward to see that there exists some constant $C$ independent of $\bg{\mu}^j$ such that $\mbb{E}\bigl[\sup_{t\in[0,T]}|X_t^{i,\bg{\mu}^j}|^2\bigr]\leq C$
and 
\bea
\mbb{E}\Bigl[\sup_{t\in[0,T]}|X_t^{i,\bg{\mu}^j}|^4|\calf_0\Bigr]^\frac{1}{2}\leq C\bigl(1+|\xi^i|^2\bigr),\quad W_2(\call(X_t^{i,\bg{\mu}^j}),\call(X_s^{i,\bg{\mu}^j}))\leq C|t-s|^\frac{1}{2}, ~\forall t,s,\in[0,T].\nn
\eea
Therefore,  just repeating the arguments used in the proof for Theorem~\ref{th-mfg-bounded},
we can show that $\Phi$ is a self-map on a closed and convex subset $\cale$ of $\calc([0,T];\calm_f^1(\mbb{R}^d))^2$,
\bea
&&\cale:=\Bigl\{(\bg{\mu}^1,\bg{\mu}^2)\in \calc([0,T];\calp_2(\mbb{R}^d))^2; \nn \\
&&\qquad \forall a\geq 1, 1\leq i\leq 2,  ~~ \sup_{t\in[0,T]}\int_{|x|\geq a}|x|^2\mu_t^i(dx)
\leq C\Bigl(a^{-1}+\mbb{E}\bigl[|\xi^i|^2\bold{1}_{\{|\xi^i|\geq \sqrt{a}\}}\bigr]\Bigr)^\frac{1}{2}\Bigr\}, 
\label{domain-E}
\eea
with some constant $C$ and that $\Phi(\cale)$ is a relatively compact subset of $\calc([0,T];\calp_2(\mbb{R}^d))^2$.
The continuity of the map $\Phi$ can be shown by Lemma~\ref{lemma-stability-mkv} just as in Theorem~\ref{th-mfg-bounded}.
Schauder FPT now guarantees the existence of  a fixed point for map $\Phi$, which then establishes
the existence of solution to the system of FBSDEs $(\ref{mkv-mftc-master})$.
\end{proof}
\end{theorem}

\subsection{Nash MFTC equilibrium for small $T$ or small coupling}
Here is the main result of this section.
\begin{theorem}
\label{th-mftc-main}
Under Assumption {\rm{\tbf{(MFTC-a)}}}, there exists some positive constant $c$ depending only on $(L,\lambda)$ such that,
for any $T\leq c$, the system of FBSDEs $(\ref{mkv-mftc-master})$ (and hence the matching problem $(\ref{mftc-matching})$)
is solvable for any $\xi^1,\xi^2\in \mbb{L}^2(\Omega,\calf_0,\mbb{P};\mbb{R}^d)$.
\begin{proof}
As in the proof for Theorem~\ref{th-mfg-main}, we use the push-forward $\phi_n\circ \mu$ of the measure $\mu\in\calp_2(\mbb{R}^d)$
defined by the map $\displaystyle \mbb{R}^d\ni x\mapsto \frac{nx}{\max(M_2(\mu),n)}$.
For eacn $n\in \mbb{N}$, we introduce the approximated coefficient functions by
\bea
&&(b_{i,0}^n, b_{i,2}^n, \sigma_{i,0}^n)(t,\nu):=(b_{i,0},b_{i,2}, \sigma_{i,0})(t,\phi_n\circ \nu), \nn \\
&&f_i^n(t,x,\mu,\nu,\alpha):=f_i(t,x,\mu,\phi_n\circ \nu,\alpha),
\quad g_i^n(x,\mu,\nu):=g_i(x,\mu, \phi_n\circ \nu)~,\nn
\eea
and accordingly define
\bea
&&b_i^n(t,x,\mu,\nu,\alpha):=b_{i,0}^n(t,\nu)+b_{i,1}(t,\nu)x+\bar{b}_{i,1}(t,\nu)\bar{\mu}+b_{i,2}^n(t,\nu)\alpha~,\nn \\
&&\sigma_i^n(t,x,\mu,\nu):=\sigma_{i,0}^n(t,\nu)+\sigma_{i,1}(t,\nu)x+\bar{\sigma}_{i,1}(t,\nu)\bar{\mu}~, \nn \\
&&H_i^n(t,x,\mu,\nu,y,z,\alpha):=\langle b_i^n(t,x,\mu,\nu,\alpha),y\rangle+{\rm{tr}}[\sigma_i^n(t,x,\mu,\nu)^\top z]
+f_i^n(t,x,\mu,\nu,\alpha)~.\nn
\eea
It is obvious to see that the approximated coefficients $(b_i^n,\sigma_i^n,f_i^n,g_i^n)$ satisfy every condition in
Assumptions {\rm{\tbf{(MFTC-a,b)}}}. Moreover, the minimizer $\hat{\alpha}_i^n$ of $H_i^n$ is given by
\be
\hat{\alpha}_i^n(t,x,\mu,\nu,y)=\hat{\alpha}_i(t,x,\mu,\phi_n\circ \nu,y)~,\nn
\ee
where $\hat{\alpha}_i$ is the minimizer of the original Hamiltonian.
The regularization for $b_{i,2}$ is done solely to obtain the simple expression for $\hat{\alpha}_i^n$ as above.
By Theorem~\ref{th-mftc-bounded}, for eacn $n\in \mbb{N}$, there exists a solution 
$(X_t^{i,n},Y_t^{i,n},Z_t^{i,n})_{t\in[0,T]}, 1\leq i\leq 2$ to the system of FBSDEs $(\ref{mkv-mftc-master})$
with the approximated coefficient functions $(b_i^n,\sigma_i^n,f_i^n,g_i^n)_{1\leq i \leq 2}$.
By the estimate $(\ref{mkv-y-growth})$, there exist constants $C$ depending only on $(L,\lambda)$ and $C^\prime$ 
depending additionally on $K$ such that, for any $t\in[0,T]$,
$
|Y_t^{i,n}|\leq C\Bigl(||\xi^i||_2+|X_t^{i,n}|+\sup_{s\in[0,T]}M_2(\call(X_s^{j,n}))\Bigr)+C^\prime,~\mbb{P}{\text{-a.s.}}
$
uniformly in $n$.  Lemma~\ref{lemma-alpha-mftc} then implies that $\hat{\alpha}_i^n(t):=
\hat{\alpha}_i^n(t,X_t^{i,n},\call(X_t^{i,n}),\call(X_t^{j,n}),Y_t^{i,n})$ satisfies
$
|\hat{\alpha}_i^n(t)|\leq 
C\Bigl(||\xi^i||_2+|X_t^{i,n}|+M_2(\call(X_t^{i,n}))+\sup_{s\in[0,T]}M_2(\call(X_s^{j,n}))\Bigr)+C^\prime$.
Thus, for any $t\in[0,T]$, 
\bea
\mbb{E}\bigl[|X_t^{i,n}|^2\bigr] &\leq& C\mbb{E}\Bigl[|\xi^i|^2+\int_0^t\bigl[ |b_i^n(s,X_s^{i,n},\call(X_s^{i,n}),\call(X_s^{j,n}),
\hat{\alpha}_i^n(s))|^2 \nn \\
&&\hspace{20mm}+ |\sigma_i^n(s,X_s^{i,n},\call(X_s^{i,n}),\call(X_s^{j,n}))|^2\bigr]ds\Bigr]\nn \\
&\leq& C\Bigl(||\xi^i||_2^2+T\sup_{s\in[0,T]}M_2(\call(X_s^{j,n}))^2+\int_0^t\sum_{j=1}^2 \mbb{E}[|X_s^{j,n}|^2]ds\Bigr)+C^\prime .
\label{xsq-ineq-mftc}
\eea
Hence Gronwall's inequality gives
$
\sum_{i=1}^2 \sup_{t\in[0,T]}\mbb{E}\bigl[|X_t^{i,n}|^2\bigr]\leq C^\prime+C T \sum_{i=1}^2 \sup_{t\in[0,T]}M_2(\call(X_t^{i,n}))^2$,
where $C^\prime$ now depends also on $||\bg{\xi}||_2$.
Therefore there exists some constant $c$ depending only on $(L,\lambda)$ such that, for any $T\leq c$, 
\bea
\sum_{i=1}^2 \sup_{t\in[0,T]}\mbb{E}\bigl[|X_t^{i,n}|^2\bigr]\leq C(L,\lambda, K, ||\bg{\xi}||_2)
\label{mkv-moment}
\eea
uniformly in $n$. For such $T\leq c$, using the estimate $(\ref{mkv-moment})$, we get by the standard technique that
\bea
&&W_2(\call(X_t^{i,n}),\call(X_s^{i,n}))\leq C(L,\lambda, K, ||\bg{\xi}||_2) |t-s|^\frac{1}{2}~,\nn \\
&&\mbb{E}\Bigl[\sup_{t\in[0,T]}|X_t^{i,n}|^4|\calf_0\Bigr]^\frac{1}{2}\leq C(L,\lambda, K, ||\bg{\xi}||_2)\bigl(1+|\xi^i|^2\bigr)~, \nn
\eea
uniformly in $n$.  We thus see  
that $(\call(X_t^{1,n})_{t\in[0,T]},\call(X_t^{2,n})_{t\in[0,T]})_{n\geq 1}$
is a relatively compact subset of $\calc([0,T];\calp_2(\mbb{R}^d))^2$. 
Upon extracting some subsequence,
there exist $\bg{\mu}^1, \bg{\mu}^2\in \calc([0,T];\calp_2(\mbb{R}^d))$ such that
$\lim_{n\rightarrow \infty}\sup_{t\in[0,T]}W_2(\call(X_t^{i,n}), \mu_t^{i})=0, 1\leq i\leq 2$.
By letting $(X_t^i,Y_t^i,Z_t^i)_{t\in[0,T], 1\leq i\leq 2}$ denote the solution to the FBSDE $(\ref{adjoint-mftc})$
with the flows $(\bg{\mu}^1, \bg{\mu}^2)$,  
we can prove that
$(X_t^i,Y_t^i,Z_t^i)_{t\in[0,T], 1\leq i\leq 2}$ is actually a solution to $(\ref{mkv-mftc-master})$ by 
the stability estimate in Lemma~\ref{lemma-stability-mkv} and
the same 
arguments used in the proof for Theorem~\ref{th-mfg-main}.
\end{proof}
\end{theorem}

As in Section~\ref{sec-mfg}, it is possible to guarantee the existence of an equilibrium for a given $T$ with quadratic cost functions
by making the couplings between FSDE and BSDE small enough. 
\begin{theorem}
\label{th-mftc-small-coupling}
Under Assumption {\rm{\tbf{(MFTC-a)}}} and a given $T$, the system of FBSDEs $(\ref{mkv-mftc-master})$ (and hence he matching problem
$(\ref{mftc-matching})$) is solvable for any $\xi^1,\xi^2\in \mbb{L}^2(\Omega,\calf_0,\mbb{P};\mbb{R}^d)$ if $\lambda^{-1}||b_{i,2}||_{\infty}$,
$1\leq i\leq 2$ are small enough.
\begin{proof}
As in the proof of Theorem~\ref{th-mfg-small-coupling}, the term involving $\sup_{t\in[0,T]}M_2(\call(X_s^{j,n}))^2$ in $(\ref{xsq-ineq-mftc})$
is proportional to $\lambda^{-1}||b_{i,2}||_{\infty}$. Hence,  if we make this factor small enough,
we obtains the estimate $(\ref{mkv-moment})$ for a given $T$. The remaining arguments are the same as in 
the proof for Theorem~\ref{th-mftc-main}.
\end{proof}
\end{theorem}

\begin{remark}
There is no difficulty to generalize all the analyses in Section~\ref{sec-mftc} for handling any finite number of populations $1\leq i\leq m$.
\end{remark}
\section{Games among Cooperative and non-Cooperative Populations}
\label{sec-mftc-mfg}
As a natural extension of Sections~\ref{sec-mfg} and \ref{sec-mftc}, we now study a 
mean-field equilibrium with two populations, where the agents in the first population (P-1) cooperate
by adopting the same feedback strategy while each agent in the second population (P-2)
competes with every other agent.  As before, we assume that the agents in each population share 
the same cost functions as well as the coefficient functions of their state dynamics.
Let us call the large population limit of this problem {\it Nash MFTC-MFG Problem}.
In Section~\ref{sec-mftc-mfg-fa}, under some additional assumptions, we shall see that the mean-field solution 
obtained in this section actually forms an approximate Nash equilibrium for the 
corresponding problem with a large but finite number of agents.
One of the motives to study  this problem
is to treat a situation, for example,  where a large number of oil producers are competing 
to maximize their profits while a part of them are members of a certain association, such as OPEC,  
cooperating within the group to maintain  a favorable level of oil price.  Since the analysis can be generalized to 
any finite number of populations, it may have many interesting applications.

\subsection{Definition of Nash MFTC-MFG problem}
We formulate the problem in the following way.\\
(i) Fix any two deterministic flows of probability measures $(\bg{\mu}^i=(\mu_t^i)_{t\in[0,T]})_{i\in\{1,2\}}$
given on $\mbb{R}^d$.\\
(ii) Solve the two optimal control problems 
\bea
\inf_{\bg{\alpha}^1\in \mbb{A}_1}J_1^{\bg{\mu}^2}(\bg{\alpha}^1),\quad \inf_{\bg{\alpha}^2\in \mbb{A}_2}J_2^{\bg{\mu}^2,\bg{\mu}^1}(\bg{\alpha}_2)
\label{mftc-mfg-problem}
\eea
over some admissible strategies $\mbb{A}_i~(i\in\{1,2\})$, where
\bea
&&J_1^{\bg{\mu}^2}(\bg{\alpha}^1):=\mbb{E}\Bigl[\int_0^T f_1(t,X_t^1,\call(X_t^1), \mu_t^2,\alpha_t^1)dt+g_1(X_T^1,\call(X_T^1),\mu_T^2)\Bigr]~, \nn\\
&&J_2^{\bg{\mu}^2,\bg{\mu}^1}(\bg{\alpha}^2):=\mbb{E}\Bigl[\int_0^T f_2(t,X_t^2,\mu_t^2,\mu_t^1,\alpha_t^2)dt+g_2(X_T^2,\mu_T^2,\mu_T^1)\Bigr]~, \nn
\eea
subject to the dynamic constraints
\bea
&&dX_t^1=b_1(t,X_t^1,\call(X_t^1),\mu_t^2,\alpha_t^1)dt+\sigma_1(t,X_t^1,\call(X_t^1),\mu_t^2)dW_t^1~, \nn \\
&&dX_t^2=b_2(t,X_t^2,\mu_t^2,\mu_t^1,\alpha_t^2)dt+\sigma_2(t,X_t^2,\mu_t^2,\mu_t^1)dW_t^2~, \nn
\eea
for $t\in[0,T]$ 
with $\bigl(X_0^i=\xi^i\in \mbb{L}^2(\Omega,\calf_0,\mbb{P};\mbb{R}^d)\bigr)_{1\leq i\leq 2}$.
Notice that the first control problem is of McKean-Vlasov type which represents the large population limit of 
cooperative agents. For each population $i\in\{1,2\}$, we suppose that $\mbb{A}_i$ is the set of $A_i$-valued $\mbb{F}^i$-progressively 
measurable processes $\bg{\alpha}^i$ satisfying $\mbb{E}\int_0^T |\alpha_t^i|^2 dt<\infty$ and
$A_i\subset \mbb{R}^k$ is closed and convex, as before.  \\
(iii) Find a pair of probability flows $(\bg{\mu}^1, \bg{\mu}^2)$ as a solution to the matching problem:
\bea
\forall t\in[0,T], \quad \mu_t^1=\call(\hat{X}_t^{1,\bg{\mu}^2}), \quad
\mu_t^2=\call(\hat{X}_t^{2,\bg{\mu}^2, \bg{\mu}^1})~, 
\label{mftc-mfg-matching}
\eea
where $(\hat{X}^{1,\bg{\mu}^2})$ and $(\hat{X}^{2,\bg{\mu}^2,\bg{\mu}^1})$ are the solutions to the optimal control problems in (ii).
\\

Throughout Section~\ref{sec-mftc-mfg}, the major assumptions for the coefficients $(b_1,\sigma_1,f_1,g_1)$ of the first population (P-1) 
are given by {\rm{\tbf{(MFTC-a)}}}, and those for the coefficients $(b_2,\sigma_2,f_2,g_2)$ of the second population (P-2) are given 
by {\rm{\tbf{(MFG-a)}}}. 
We have already learned from Theorems~\ref{th-fbsde-existence} and \ref{th-mftc-existence} that the solution 
to  each of the optimal control problems in $(\ref{mftc-mfg-problem})$ for given deterministic flows $\bg{\mu}^1, \bg{\mu}^2\in 
\calc([0,T];\calp_2(\mbb{R}^d))$ is characterized by the uniquely solvable FBSDEs,
\bea
&&dX_t^1=b_1(t,X_t^1,\call(X_t^1), \mu_t^2,\hat{\alpha}_1(t,X_t^1,\call(X_t^1),\mu_t^2,Y_t^1))dt+\sigma_1(t,X_t^1,\call(X_t^1),\mu_t^2)dW_t^1~, \nn\\
&&dY_t^1=-\part_x H_1(t,X_t^1,\call(X_t^1),\mu_t^2,Y_t^1,Z_t^1,\hat{\alpha}_1(t,X_t^1,\call(X_t^1),\mu_t^2,Y_t^1))dt\nn \\
&&\qquad-\wt{\mbb{E}}\bigl[\part_\mu H_1(t,\wt{X}_t^1,\call(X_t^1),\mu_t^2,\wt{Y}_t^1,\wt{Z}_t^1,
\hat{\alpha}_1(t,\wt{X}_t^1,\call(X_t^1),\mu_t^2,\wt{Y}_t^1))(X_t^1)\bigr]dt+Z_t^1 dW_t^1~,
\label{p1-fbsde}
\eea
with $X_0^1=\xi^1$ and $Y_T^1=\part_x g_1(t,X_T^1,\call(X_T^1),\mu_T^2)
+\wt{\mbb{E}}\bigl[\part_\mu g_1(\wt{X}_T^1,\call(X_T^1),\mu_T^2)(X_T^1)\bigr]$, and
\bea
&&dX_t^2=b_2(t,X_t^2,\mu_t^2,\mu_t^1,\hat{\alpha}_2(t,X_t^2,\mu_t^2,\mu_t^1,Y_t^2))dt+\sigma_2(t,X_t^2,\mu_t^2,\mu_t^1)dW_t^2~, \nn \\
&&dY_t^2=-\part_x H_2(t,X_t^2,\mu_t^2,\mu_t^1,Y_t^2,Z_t^2,\hat{\alpha}_2(t,X_t^2,\mu_t^2,\mu_t^1,Y_t^2))dt+Z_t^2 dW_t^2~, 
\label{p2-fbsde}
\eea
with $X_0^2=\xi^2$ and $Y_T^2=\part_x g_2(X_T^2,\mu_T^2,\mu_T^1)$, respectively.
Here, the Hamiltonian $H_i$ and its minimizer $\hat{\alpha}_i$ are defined as before using the corresponding coefficients $(b_i,\sigma_i, f_i)$.

\subsection{MFTC-MFG equilibrium under additional boundedness}
In order to establish the existence of an equilibrium $(\ref{mftc-mfg-matching})$, we have to show the 
existence of a solution to the following  coupled system of FBSDEs:
\bea
&&dX_t^1=b_1(t,X_t^1,\call(X_t^1), \call(X_t^2),\hat{\alpha}_1(t,X_t^1,\call(X_t^1),\call(X_t^2),Y_t^1))dt+
\sigma_1(t,X_t^1,\call(X_t^1),\call(X_t^2))dW_t^1~, \nn\\
&&dY_t^1=-\part_x H_1(t,X_t^1,\call(X_t^1),\call(X_t^2),Y_t^1,Z_t^1,\hat{\alpha}_1(t,X_t^1,\call(X_t^1),\call(X_t^2),Y_t^1))dt\nn \\
&&\qquad-\wt{\mbb{E}}\bigl[\part_\mu H_1(t,\wt{X}_t^1,\call(X_t^1),\call(X_t^2),\wt{Y}_t^1,\wt{Z}_t^1,
\hat{\alpha}_1(t,\wt{X}_t^1,\call(X_t^1),\call(X_t^2),\wt{Y}_t^1))(X_t^1)\bigr]dt+Z_t^1 dW_t^1~, \nn
\eea
with $X_0^1=\xi^1$ and $Y_T^1=\part_x g_1(t,X_T^1,\call(X_T^1),\call(X_T^2))+\wt{\mbb{E}}\bigl[\part_{\mu}g_1(\wt{X}_T^1,
\call(X_T^1),\call(X_T^2))(X_T^1)\bigr]$, and 
\bea
&&dX_t^2=b_2(t,X_t^2,\call(X_t^2),\call(X_t^1),\hat{\alpha}_2(t,X_t^2,\call(X_t^2),\call(X_t^1),Y_t^2))dt+
\sigma_2(t,X_t^2,\call(X_t^2),\call(X_t^1))dW_t^2~, \nn\\
&&dY_t^2=-\part_x H_2(t,X_t^2,\call(X_t^2),\call(X_t^1),Y_t^2,Z_t^2,\hat{\alpha}_2(t,X_t^2,\call(X_t^2),\call(X_t^1),Y_t^2))dt+Z_t^2dW_t^2~, 
\label{mftc-mfg-master}
\eea
with $X_0^2=\xi^2$ and $Y_T^2=\part_x g_2(X_T^2,\call(X_T^2),\call(X_T^1))$.

In this section, our goal is to prove the following result.
\begin{theorem}
\label{th-mixed-bounded}
Under Assumptions {\rm{\tbf{(MFTC-a,b)}}} for the coefficients $(b_1,\sigma_1,f_1,g_1)$
and Assumptions ${\rm{\tbf{(MFG-a,b)}}}$ for  the coefficients $(b_2,\sigma_2,f_2,g_2)$, 
the system of FBSDEs $(\ref{mftc-mfg-master})$ (and hence the matching problem $(\ref{mftc-mfg-matching})$)
is solvable for any $\xi^1,\xi^2\in \mbb{L}^2(\Omega,\calf_0,\mbb{P};\mbb{R}^d)$.
\begin{proof}
We let $(X_t^{1,\bg{\mu}^2}, Y_t^{1,\bg{\mu}^2},Z_t^{1,\bg{\mu}^2})_{t\in[0,T]}$ 
and $(X_t^{2,\bg{\mu}^2,\bg{\mu}^1}, Y_t^{2,\bg{\mu}^2,\bg{\mu}^1}, Z_t^{2,\bg{\mu}^2,\bg{\mu}^1})_{t\in[0,T]}$
denote the solutions to the FBSDE $(\ref{p1-fbsde})$ and $(\ref{p2-fbsde})$ respectively 
for given flows of probability measures $(\bg{\mu}^1,\bg{\mu}^2)$.
By defining the map $\Phi$ as
\bea
\Phi:\calc([0,T];\calp_2(\mbb{R}^d))^2\ni (\bg{\mu}^1,\bg{\mu}^2)\mapsto 
\bigl(\call(X_t^{1,\bg{\mu}^2})_{t\in[0,T]}, \call(X_t^{2,\bg{\mu}^2,\bg{\mu}^1})_{t\in[0,T]}\bigr)\in \calc([0,T];\calp_2(\mbb{R}^d))^2~,
\eea
the claim is proved once we find a fixed point of the map $\Phi$.

It is the direct result of Theorem~\ref{th-mftc-bounded} for (P-1) and Theorem~\ref{th-mfg-bounded} for (P-2)
that there exists a constant $C$ independent of $\bg{\mu}^1$ and $\bg{\mu}^2$ such that, for any $t,s\in[0,T]$, 
\bea
\label{p1-estimate-bounded}
&&\mbb{E}\bigl[\sup_{t\in[0,T]}|X_t^{1,\bg{\mu}^2}|^4|\calf_0\bigr]^\frac{1}{2}
\leq C\bigl(1+|\xi^1|^2\bigr), \quad \mbb{E}\bigl[\sup_{t\in[0,T]}|X_t^{2,\bg{\mu}^2,\bg{\mu}^1}|^4|\calf_0\bigr]^\frac{1}{2}\leq C\bigl(1+|\xi^2|^2\bigr),~
\nn \\
&&W_2(\call(X_t^{1,\bg{\mu}^2}),\call(X_s^{1,\bg{\mu}^2}))\leq C|t-s|^\frac{1}{2}~, \quad
W_2(\call(X_t^{2,\bg{\mu}^2,\bg{\mu}^1}), \call(X_s^{2,\bg{\mu}^2,\bg{\mu}^1}))\leq C|t-s|^\frac{1}{2}.
\eea
Thus we can show that, for the same form of closed and convex subset $\cale$ of $C([0,T];\calm_f^1(\mbb{R}^d))^2$ in $(\ref{domain-E})$
with sufficiently large $C$, that $\Phi$ maps $\cale$ into itself and also that $\Phi(\cale)$ is 
a relatively compact subset of $\calc([0,T];\calp_2(\mbb{R}^d))^2$. 
The continuity of the map $\Phi$ can be shown by Lemmas~\ref{lemma-stability-mkv} and \ref{lemma-stability}
just as in Theorems~\ref{th-mftc-bounded} and \ref{th-mfg-bounded}.
By Schauder FPT, the claim is proved.
\end{proof}
\end{theorem}

\subsection{MFTC-MFG equilibrium for small $T$ or small coupling}
We now give the main result of Section~\ref{sec-mftc-mfg}.
\begin{theorem}
\label{th-mftc-mfg-main}
Under Assumption {\rm{\tbf{(MFTC-a)}}} for the coefficients $(b_1,\sigma_1,f_1,g_1)$
and Assumption {\rm{\tbf{(MFG-a)}}} for the coefficients $(b_2,\sigma_2,f_2,g_2)$, 
there exists some positive constant $c$ depending only on $(L,\lambda)$ such that, for any $T\leq c$,
the system of FBSDEs $(\ref{mftc-mfg-master})$ (and hence matching problem $(\ref{mftc-mfg-matching})$)
is solvable for any $\xi^1,\xi^2\in\mbb{L}^2(\Omega,\calf_0,\mbb{P};\mbb{R}^d)$. 
\begin{proof}
Let us introduce the approximated functions $(b_1^n, \sigma_1^n,f_1^n, g_1^n)_{n\geq 1}$ as in Theorem~\ref{th-mftc-main}
and also $(b_2^n, \sigma_2^n, f_2^n, g_2^n)_{n\geq 1}$ as in Theorem~\ref{th-mfg-main},
which satisfy Assumptions {\rm{\tbf{(MFTC-a,b)}}} and Assumptions {\rm{\tbf{(MFG-a,b)}}} for each $n$, respectively.
Theorem~\ref{th-mixed-bounded} then guarantees that there exists a solution to the system of FBSDEs $(\ref{mftc-mfg-master})$
with the approximated functions $(b_i^n, \sigma_i^n, f_i^n, g_i^n)_{1\leq i\leq 2}$
for each $n$. We let $(X_t^{i,n},Y_t^{i,n},Z_t^{i,n})_{t\in[0,T]}, 1\leq i\leq 2$ denote the corresponding solution.

Since inequalities $(\ref{xsq-ineq-mftc})$ and $(\ref{xsq-ineq-mfg})$ still hold, we can show that there exist
constants $C$ depending only on $(L,\lambda)$ and $C^\prime$ depending additionally on $K$ such that,
\bea
&&\mbb{E}\bigl[|X_t^{1,n}|^2\bigr]\leq C\Bigl(||\xi^1||_2^2+T\sup_{s\in[0,T]}M_2(\call(X_s^{2,n}))^2+\int_0^t\mbb{E}
\sum_{i=1}^2 \bigl[|X_s^{i,n}|^2\bigr]ds\Bigr)
+C^\prime~, \nn \\
&&\mbb{E}\bigl[|X_t^{2,n}|^2\bigr]\leq C\Bigl(||\xi^2||_2^2+T\sup_{s\in[0,T]}\sum_{i=1}^2 M_2(\call(X_s^{i,n}))^2
+\int_0^t \sum_{i=1}^2\mbb{E}\bigl[|X_s^{i,n}|^2\bigr]ds\Bigr)+C^\prime~. \nn
\eea
Hence we get, by Gronwall's inequality, that
\bea
\sup_{t\in[0,T]}\sum_{i=1}^2\mbb{E}\bigl[|X_t^{i,n}|^2\bigr]\leq C\Bigl(||\bg{\xi}||_2^2+T\sup_{s\in[0,T]}\sum_{i=1}^2 M_2(\call(X_s^{i,n}))^2
\Bigr)+C^\prime.\nn
\eea
Therefore there exists a positive constant $c$ depending only on $(L,\lambda)$ such that, for any $T\leq c$, 
\bea
\sup_{t\in[0,T]}\sum_{i=1}^2\mbb{E}\bigl[|X_t^{i,n}|^2\bigr]\leq C^\prime\bigl(1+||\bg{\xi}||_2^2\bigr)~. \nn
\eea
Using the linear growth property of $\hat{\alpha}_i^n$ in $|X_t^{i,n}|$, 
we can show that 
 $(\call(X_t^{1,n})_{t\in[0,T]},\call(X_t^{2,n})_{t\in[0,T]})_{n\geq 1}$ is a relatively compact subset of $\calc([0,T];\calp_2(\mbb{R}^d))^2$. 
The remaining arguments proceed in exactly the same way as in the proofs for Theorems~\ref{th-mftc-main} and \ref{th-mfg-main}.
\end{proof}
\end{theorem}

\begin{theorem}
\label{th-mftc-mfg-small-coupling}
Under Assumption {\rm{\tbf{(MFTC-a)}}} for the coefficients $(b_1,\sigma_1,f_1,g_1)$
and Assumption {\rm{\tbf{(MFG-a)}}} for the coefficients $(b_2,\sigma_2,f_2,g_2)$ and a given $T$,
the system of FBSDEs $(\ref{mftc-mfg-master})$ (and hence matching problem $(\ref{mftc-mfg-matching})$)
is solvable for any $\xi^1,\xi^2\in\mbb{L}^2(\Omega,\calf_0,\mbb{P};\mbb{R}^d)$ if $\lambda^{-1}||b_{i,2}||_{\infty}$, 
$1\leq i\leq 2$ are small enough. 
\begin{proof}
The claim can be proved in a completely parallel way to Theorems~\ref{th-mfg-small-coupling} and \ref{th-mftc-small-coupling}.
\end{proof}
\end{theorem}

\begin{remark}
As in Sections~\ref{sec-mfg} and \ref{sec-mftc}, the analysis can be easily extended for the situation with
 any finite number of cooperative and non-cooperative populations.
\end{remark}

\section{Approximate Equilibrium for MFG with Finite Agents}
\label{sec-mfg-fa}
In the remaining sections, we investigate quantitative relationships
between the solutions to the mean field games 
obtained in the previous three sections and those to their associated games with a finite number of agents.
We make use of  the techniques developed in 
\cite{Sznitman, Bensoussan-LQ, Cardaliaguet-note, Carmona-Delarue-MFG,Carmona-Delarue-MFTC}
and in particular Chapter 6 in \cite{Carmona-Delarue-2} with appropriate generalizations to fit our situation.
First, in this section, we shall study the problem associated with the multi-population mean field game
solved in Section~\ref{sec-mfg}. 
Throughout the section, we assume that the conditions used either in Theorem~\ref{th-mfg-main} or Theorem~\ref{th-mfg-small-coupling}
are satisfied. We let $(\bg{\mu}^1,\bg{\mu}^2)\in \calc([0,T];\calp_2(\mbb{R}^d))^2$ denote a solution to the matching problem $(\ref{mfg-fixed-point})$.
\subsection{Convergence of approximate optimal controls}
For each population $1\leq i \leq 2$, we suppose that  there are  $N_i$ agents who are labeled by $p$.
Let us first introduce  $N_i$ independent and identically distributed (i.i.d.) copies of the state process in the
mean field setup:
\bea
d \ul{X}_t^{i,p}=b_i(t,\ul{X}_t^{i,p},\mu_t^i,\mu_t^j,\ul{\hat{\alpha}}_t^{i,p})dt+
\sigma_i(t,\ul{X}_t^{i,p},\mu_t^i,\mu_t^j)dW_t^{i,p}~,
\label{dynamics-ulX}
\eea
for $1\leq p\leq N_i$, $j \neq i$, $t\in[0,T]$ with $\ul{X}_0^{i,p}=\xi^{i,p}$, and
\bea
\ul{\hat{\alpha}}_t^{i,p}:= \hat{\alpha}_i(t,\ul{X}_t^i,\mu_t^i,\mu_t^j,u_i(t,\ul{X}_t^i)) \nn
\eea
for any $t\in[0,T]$.
Here, $(\xi^{i,p})_{1\leq p\leq N_i}$ is the set of i.i.d random variables with $\call(\xi^{i,p})=\mu_0^i$,
and $(\bg{W}^{i,p}=(W^{i,p}_t)_{t\in[0,T]})_{1\leq p\leq N_i}$ are independent standard Brownian motions, which 
are also independent from $(\xi^{i,p})_{1\leq p\leq N_i}$. Moreover, they are  assumed to be independent from those 
in the other population. In other words, all of the set $(\xi^{i,p},\bg{W}^{i,p})_{1\leq p\leq N_i, 1\leq i\leq 2}$  are assumed to be independent.
$u_i$ is the decoupling field given in Theorem~\ref{th-fbsde-existence} associated with the equilibrium flows of probability 
measures $(\bg{\mu}^i,\bg{\mu}^j)$. $\hat{\alpha}_i$ is the minimizer of the Hamiltonian for the population $i$
defined in $(\ref{def-alpha})$.
By construction, $(\bg{\ul{X}}^{i,p})_{1\leq p\leq N_i}$ are i.i.d. processes satisfying $\call(X_t^{i,p})=\mu_t^i$,
$\forall t\in[0,T]$.
We denote the empirical distribution for $(\bg{\ul{X}}^{i,p})_{1\leq p\leq N_i}$ by
\bea
\ul{\mu}_t^i:=\frac{1}{N_i}\sum_{p=1}^{N_i}\del_{\ul{X}_t^{i,p}}~. \nn
\eea
In the remainder, the complete probability space $(\Omega,\calf,\mbb{P})$ is enlarged accordingly to support $(\xi^{i,p},\bg{W}^{i,p})_{1\leq p\leq N_i, 1\leq i\leq 2}$
and the filtration $\mbb{F}$ is assumed to be generated by $(\xi^{i,p},\bg{W}^{i,p})_{1\leq p\leq N_i, 1\leq i\leq 2}$ 
with complete and right-continuous augmentation.
\begin{lemma}
\label{lemma-cl}
Suppose that the conditions either for Theorem~\ref{th-mfg-main} or Theorem~\ref{th-mfg-small-coupling} are satisfied.
Then, for each population $1\leq i\leq 2$, there exists some sequence $(\ep_{N_i})_{N_i\geq 1}$ that tends to $0$ as $N_i$ tends to $\infty$
and some constant $C$ such that
\bea
\sup_{t\in[0,T]}\mbb{E}\Bigl[W_2(\ul{\mu}_t^i,\mu_t^i)^2\Bigr]\leq C\ep_{N_i}^2~.\nn
\eea
Furthermore, when $\mu_0^i\in \calp_r(\mbb{R}^d)$ with $r>4$, we have
an explicit estimate 
\bea
\ep_{N_i}^2=N_i^{-2/\max(d,4)}(1+\ln(N_i)\bold{1}_{d=4})~.\nn
\eea
\begin{proof}
When $\mu_0^i\in \calp_r(\mbb{R}^d)$, $\forall r\geq 2$, it is standard to check
\bea
\sup_{t\in[0,T]}M_r(\mu_r^i)^r\leq \mbb{E}\Bigl[\sup_{t\in[0,T]}|\ul{X}_t^{i,p}|^r\Bigr]\leq C\Bigl(1+M_r(\mu_0^i)^r\Bigr) \nn
\eea
with some $C$ independent of $N_i$ thanks to the linear growth of the coefficients in $(\ref{dynamics-ulX})$.
Then, the last claim is the direct result of Theorem 5.8 and Remark 5.9 in \cite{Carmona-Delarue-1}.

As for the first claim,  $(5.19)$ in \cite{Carmona-Delarue-1} implies
\bea
\lim_{N_i\rightarrow \infty}\mbb{E}\bigl[W_2(\ul{\mu}_t^i,\mu_t^i)^2\bigr]=0
\label{pointwise}
\eea
for each $t$. In order to prove the uniform convergence in $t$~\footnote{See the arguments 
leading to the estimate $(2.15)$ in the proof of Theorem 2.12 in \cite{Carmona-Delarue-2}.}, it suffices to show that
there exists a compact set $\calk\subset \calc([0,T];\mbb{R})$ such that
\bea
\Bigl(\mbb{E}[W_2(\ul{\mu}_t^i,\mu_t^i)^2]_{t\in[0,T]}\Bigr)_{N_i\geq 1}\subset \calk~.\nn
\eea
In fact, if this is the case, every subsequence has a uniformly convergent subsequence, 
all of which converge to 0 due to the pointwise convergence in $(\ref{pointwise})$.
Hence, the whole sequence must uniformly converges to 0.
The boundedness can be checked by
\bea
\sup_{t\in[0,T]}\mbb{E}\bigl[W_2(\ul{\mu}_t^i,\mu_t^i)^2\bigr]\leq 2\sup_{t\in[0,T]}\Bigl(
\frac{1}{N_i}\sum_{p=1}^{N_i}\mbb{E}[|\ul{X}_t^{i,p}|^2]+M_2(\mu_t^i)^2\Bigr)\leq 4\sup_{t\in[0,T]}M_2(\mu_t^i)^2\leq C~.\nn
\eea
Moreover,  for any $0\leq t, s\leq T$, 
\bea
&&\Bigl|\mbb{E}\bigl[W_2(\ul{\mu}_t^i,\mu_t^i)^2\bigr]-\mbb{E}\bigl[W_2(\ul{\mu}_s^i,\mu_s^i)^2\bigr]\Bigr|\leq C\mbb{E}\bigl[(W_2(\ul{\mu}_t^i,\mu_t^i)-W_2(\ul{\mu}_s^i,\mu_s^i))^2\bigr]^\frac{1}{2}\nn \\
&&\leq C\Bigl(\mbb{E}[W_2(\ul{\mu}_t^i,\ul{\mu}_s^i)^2]+W_2(\mu_t^i,\mu_s^i)^2\Bigr)^\frac{1}{2}
\leq C\mbb{E}\bigl[|\ul{X}_t^{i,p}-\ul{X}_s^{i,p}|^2\bigr]^\frac{1}{2}\leq C|t-s|^\frac{1}{2}~,\nn
\eea
which implies the equicontinuity. Arzela-Ascoli theorem guarantees the desired compactness.
\end{proof}
\end{lemma}

\begin{assumption}{\rm{\tbf{(MFG-FA)}}}
On top of Assumption {\rm{\tbf{(MFG-a)}}},  either $T$ or $(\lambda^{-1}||b_{i,2}||_{\infty})_{1\leq i\leq 2}$ is small enough 
to satisfy the conditions for Theorem~\ref{th-mfg-main} or Theorem~\ref{th-mfg-small-coupling}. Moreover,  for $1\leq i\leq 2$, \\
{\rm{\tbf{(A1)}}} There exists some constant $K$ such that
\bea
|(b_{i,0},\sigma_{i,0})(t,\mu^\prime,\nu^\prime)-(b_{i,0},\sigma_{i,0})(t,\mu,\nu)|\leq K\bigl(W_2(\mu^\prime,\mu)+W_2(\nu^\prime,\nu)\bigr) \nn
\eea 
for any $t\in[0,T], \mu^\prime,\mu,\nu^\prime,\nu\in \calp_2(\mbb{R}^d)$,  and $b_{i,1},b_{i,2}$ as well as $\sigma_{i,1}$ 
are independent of the measure arguments. \\
{\rm{\tbf{(A2)}}} $f_i$ and $g_i$ are local Lipschitz continuous with respect to the measure arguments i.e., there exists some constant $K$
for any $t\in[0,T], x\in\mbb{R}^d, \mu^\prime,\mu, \nu^\prime,\nu\in \calp_2(\mbb{R}^d)$ and $\alpha\in A_i$, such that,
\footnote{Although there is no $(t,\alpha)$ dependence in $g_i$, we slightly abuse the notation to save the space. 
Note that the local Lipschitz property for the arguments $(x,\alpha)$ follows from {\rm{\tbf{(MFG-a)}}}.} 
\bea
&&|(f_i,g_i)(t,x,\mu^\prime,\nu^\prime, \alpha)-(f_i,g_i)(t,x,\mu,\nu,\alpha)|\nn \\
&&\leq K\Bigl(1+|x|+M_2(\mu^\prime)+M_2(\mu)+M_2(\nu^\prime)+M_2(\nu)
+|\alpha|\Bigr)\bigl(W_2(\mu^\prime,\mu)+W_2(\nu^\prime,\nu)\bigr)~.\nn 
\eea
\end{assumption}
\vspace{5mm}
For $1\leq i, j\leq 2$, $j\neq i$, $1\leq p\leq N_i$, let us consider the following state dynamics.
\bea
dX_t^{i,p}=b_i(t,X_t^{i,p},\ol{\mu}_t^i,\ol{\mu}_t^j,\hat{\alpha}_t^{i,p})dt+\sigma_i(t,X_t^{i,p},\ol{\mu}_t^i,\ol{\mu}_t^j)dW_t^{i,p}~, 
\label{mfg-fa-approx}
\eea
for $t\in[0,T]$ with $X_0^{i,p}=\xi^{i,p}$. Here, $\ol{\mu}_t^i:=\frac{1}{N_i}\sum_{p=1}^{N_i}\del_{X_t^{i,p}}$ is the empirical distribution and
\be
\hat{\alpha}_t^{i,p}:=\hat{\alpha}_i(t,X_t^{i,p},\mu_t^i,\mu_t^j,u_i(t,X_t^{i,p}))~.\nn
\ee
Since we have $W_2(\ol{\mu}_t^i, \ol{\mu}_t^{\prime,i})^2\leq \frac{1}{N_i}\sum_{p=1}^{N_i}|X_t^{i,p}-X_t^{\prime,i,p}|^2$, 
Assumption {\rm{\tbf{(MFG-FA) (A1)}}}, the Lipschitz continuity of $\hat{\alpha}_i$ and that of the decoupling field $u_i$
make $(\ref{mfg-fa-approx})$ an $(N_1+N_2)$-dimensional standard Lipschitz SDE. $(\bg{X}^{i,p})_{1\leq i\leq 2, 1\leq p\leq N_i}$
correspond to the state processes of the agents who adopt the feedback control function 
$[0,T]\times \mbb{R}^d\ni (t,x)\mapsto \hat{\alpha}_i(t,x,\mu_t^i,\mu_t^j,u_i(t,x))$.
The cost functional for the agent $p$ in the ith population is given by
\bea
J_i^{N_i,N_j,p}(\bg{\hat{\alpha}}^{i,(N_i)},\bg{\hat{\alpha}}^{j,(N_j)})
:=\mbb{E}\Bigl[\int_0^T f_i(t,X_t^{i,p},\ol{\mu}_t^i, \ol{\mu}_t^j, \hat{\alpha}_t^{i,p})dt+g_i(X_T^{i,p},\ol{\mu}_T^i,\ol{\mu}_T^j)\Bigr]~.\nn
\eea
On the other hand, the optimal cost functional for the mean field game in Section~\ref{sec-mfg} is 
\bea
J_i^p:=\mbb{E}\Bigl[\int_0^T f_i(t,\ul{X}_t^{i,p},\mu_t^i,\mu_t^j,\ul{\hat{\alpha}}_t^{i,p})dt+g_i(\ul{X}_T^{i,p},\mu_T^i,\mu_T^j)\Bigr]~.\nn
\eea

\begin{lemma}
\label{lemma-approx-limit}
Under Assumption {\rm{\tbf{(MFG-FA)}}}, for all $1\leq i\leq 2$, $N_i\in \mbb{N}$, $1\leq p\leq N_i$,
there exists some constant $C$ independent of $(N_i)_{i=1}^2$ such that,
\bea
|J_i^{N_i,N_j,p}(\bg{\hat{\alpha}}^{i,(N_i)},\bg{\hat{\alpha}}^{j,(N_j)})-J_i^p|\leq C\sum_{j=1}^2\ep_{N_j}~, \nn
\eea
where $\ep_{N_j}$ is the one given in Lemma~\ref{lemma-cl}.
\begin{proof}
It suffices to check the case with $(p=1)$. By Lipschitz continuity and the triangle inequality,  we have
$\mbb{E}\bigl[\sup_{s\in[0,t]}|X_s^{i,1}-\ul{X}_s^{i,1}|^2\bigl]\leq C\int_0^t\mbb{E}\Bigl[|X_s^{i,1}-\ul{X}_s^{i,1}|^2+\sum_{i=1}^2 W_2(\ol{\mu}_t^i,
\ul{\mu}_t^i)^2+\sum_{i=1}^2 W_2(\ul{\mu}_t^i,\mu_t^i)^2\Bigr]ds$.
Applying Gronwall's inequality after summing over $1\leq i\leq 2$, we get
\bea
\sum_{i=1}^2\mbb{E}\bigl[\sup_{t\in[0,T]}|X_t^{i,1}-\ul{X}_t^{i,1}|^2\bigr]\leq C\sum_{i=1}^2 \sup_{t\in[0,T]}\mbb{E}
\bigl[W_2(\ul{\mu}_t^i,\mu_t^i)^2\bigr]\leq C\sum_{i=1}^2\ep_{N_i}^2~.\nn
\eea
Hence the triangle inequality 
implies that $\sup_{t\in[0,T]}\mbb{E}\bigl[W_2(\ol{\mu}_t^i, \mu_t^i)^2\bigr]\leq C\sum_{j=1}^2\ep_{N_j}^2$.
We also see $\mbb{E}\bigl[\sup_{t\in[0,T]}|X_t^{i,1}|^2\bigr]\leq C\Bigl(1+\mbb{E}\bigl[\sup_{t\in[0,T]}|\ul{X}_t^{i,1}|^2\bigr]\Bigr)\leq C$.
Using local Lipschitz continuity, it is straightforward to conclude
\bea
&&|J_i^{N_i,N_j,1}(\bg{\hat{\alpha}}^{i,(N_i)},\bg{\hat{\alpha}}^{j,(N_j)})-J_i^1| \leq C\Bigl(1+\sup_{t\in[0,T]}\mbb{E}\Bigl[ |X_t^{i,1}|^2+|\ul{X}_t^{i,1}|^2+\sum_{j=1}^2 M_2(\ol{\mu}_t^j)^2\Bigr]^\frac{1}{2}\Bigr)\nn \\
&&\qquad\times \sup_{t\in[0,T]}\mbb{E}\Bigl[|X_t^{i,1}-\ul{X}_t^{i,1}|^2+\sum_{j=1}^2W_2(\ol{\mu}_t^j,\mu_t^j)^2\Bigr]^\frac{1}{2}\leq C\sum_{j=1}^2\ep_{N_j}~.\nn
\eea
\end{proof}
\end{lemma}
\begin{remark}
From the above analysis, we see that {\rm{\tbf{(MFG-FA) (A2)}}} is unnecessary if we only need the convergence
$J_i^{N_i,N_j,p}(\bg{\hat{\alpha}}^{i,(N_i)},\bg{\hat{\alpha}}^{j,(N_j)})\rightarrow J_i^p$ when $N_1,N_2\rightarrow \infty$.
{\rm{\tbf{(A2)}}} is just used to derive the explicit order of convergence in terms of  $(\ep_{N_i})_{i=1}^2$.
\end{remark}
\subsection{Approximate Nash Equilibrium}
In order to investigate an approximate Nash equilibrium, we suppose that one agent 
deviates from the feedback control function $[0,T]\times \mbb{R}^d\ni (t,x)\mapsto \hat{\alpha}_i(t,x,\mu_t^i,\mu_t^j,u_i(t,x))$.
By symmetry, we may assume that this is the first agent in the ith population. The state dynamics of the 
agents is now given by, $1\leq i, j \leq 2, j\neq i$, $1\leq p\leq N_i$, 
\bea
dU_t^{i,p}=b_i(t,U_t^{i,p},\ol{\nu}_t^i,\ol{\nu}_t^j,\beta_t^{i,p})dt+\sigma_i(t,U_t^{i,p},\ol{\nu}_t^i,\ol{\nu}_t^j)dW_t^{i,p}~,\nn
\eea
for $t\in[0,T]$ with $U_0^{i,p}=\xi^{i,p}$. $\bg{\beta}^{i,1}\in \mbb{H}^2$ is any $A_i$-valued $\mbb{F}$-progressively
measurable process, and 
\bea
&&\beta_t^{i,p}:=\hat{\alpha}_i(t,U_t^{i,p},\mu_t^i,\mu_t^j,u_i(t,U_t^{i,p})), ~2\leq p\leq N_i~, \nn \\
&&\beta_t^{j,q}:=\hat{\alpha}_j(t,U_t^{j,q},\mu_t^j,\mu_t^i,u_j(t,U_t^{j,q})), ~1\leq q\leq N_j~.
\label{feedback-mfg}
\eea
$\ol{\nu}_t^i:=\frac{1}{N_i}\sum_{p=1}^{N_i} \del_{U_t^{i,p}}$, $1\leq i\leq 2$  is the empirical distribution.
This is an $(N_1+N_2)$-dimensional Lipschitz SDE and hence well-defined.
The cost functional associated with the deviating agent is given by
\bea
J_i^{N_i,N_j,1}(\bg{\beta}^{i,1}, \bg{\hat{\alpha}}^{i,(N_i)^{-1}},\bg{\hat{\alpha}}^{j,(N_j)})
:=\mbb{E}\Bigl[\int_0^T f_i(t,U_t^{i,1},\ol{\nu}_t^i,\ol{\nu}_t^j,\beta_t^{i,1})dt+g_i(U_T^{i,1},\ol{\nu}_T^i,\ol{\nu}_T^j)\Bigr]~.\nn
\eea

\begin{remark}
\label{remark-open}
As we can see from the above definition of control strategies, we shall focus on the approximate Nash equilibrium 
in the sense of the closed loop framework.  The analysis for the open loop framework can be done in almost the same (actually slightly simpler) manner,
which just requires to replace the feedback forms in $(\ref{feedback-mfg})$ by 
\bea
\beta_t^{i,p}=\ul{\hat{\alpha}}_t^{i,p}, \quad 2\leq p\leq N_i, \quad
\beta_t^{j,q}=\ul{\hat{\alpha}}_t^{j,q}, \quad 1\leq q\leq N_j~.\nn
\eea
\end{remark}

Here is the main result of this section.
\begin{theorem}
\label{th-nash-mfg}
Under Assumption {\rm{\tbf{(MFG-FA)}}} with sufficiently large $(N_i)_{i=1}^2$, 
the feedback control functions $\bigl([0,T]\times \mbb{R}^d\ni (t,x)\mapsto \hat{\alpha}_i(t,x,\mu_t^i,\mu_t^j,u_i(t,x))\bigr)_{1\leq i, j\leq 2, j\neq i}$
form an $(\sum_{j=1}^2\vep_{N_j})$-approximate Nash equilibrium i.e.,
there exists some constant $C$ independent of $(N_i)_{i=1}^2$ such that
\bea
J_i^{N_i,N_j,1}(\bg{\beta}^{i,1}, \bg{\hat{\alpha}}^{i,(N_i)^{-1}},\bg{\hat{\alpha}}^{j,(N_j)})
\geq J_i^{N_i,N_j,1}(\bg{\hat{\alpha}}^{i,(N_i)},\bg{\hat{\alpha}}^{j,(N_j)})-C \sum_{j=1}^2\vep_{N_j}~\nn
\eea
for any $A_i$-valued $\mbb{F}$-progressively measurable process
$\bg{\beta}^{i,1}\in \mbb{H}^2$,  where $\bigl(\vep_{N_j}:=\max(\ep_{N_j},N_j^{-\frac{1}{2}})\bigr)_{1\leq j\leq 2}$.
\begin{proof}
{\rm{\tbf{(first step)}}} Let us introduce another dynamics
\bea
d\ol{U}_t^{i,1}=b_i(t,\ol{U}_t^{i,1},\mu_t^i,\mu_t^j,\beta_t^{i,1})dt+\sigma_i(t,\ol{U}_t^{i,1},\mu_t^i,\mu_t^j)dW_t^{i,1}\nn
\eea
for $t\in[0,T]$ with $U_0^{i,1}=\xi^{i,1}$. It is immediate to see that the following estimate holds:
\be
\mbb{E}\bigl[\sup_{t\in[0,T]}|\ol{U}_t^{i,1}|^2\bigr]\leq C\bigl(1+||\beta^{i,1}||^2_{\mbb{H}^2}\bigr).
\label{olU-norm}
\ee
The associated cost functional 
\bea
\ol{J}_i^{1}(\bg{\beta}^{i,1}):=\mbb{E}\Bigl[\int_0^T f_i(t,\ol{U}_t^{i,1},\mu_t^i,\mu_t^j,\beta_t^{i,1})dt+g_i(\ol{U}_T^{i,1},\mu_T^i,\mu_T^j)\Bigr]~\nn
\eea
satisfies, by Theorem~\ref{th-fbsde-existence},   an inequality
\bea
\ol{J}_i^{1}(\bg{\beta}^{i,1})\geq J_i^1+\lambda \mbb{E}\int_0^T |\beta_t^{i,1}-\ul{\hat{\alpha}}_t^{i,1}|^2 dt~.
\label{olJ-J-estimate}
\eea
\\
{\rm{\tbf{(second step)}}} By the linear  growth property of the coefficients, we get
\bea
&&\mbb{E}\Bigl[ |U_t^{i,p}|^2\Bigr]\leq C\Bigl(1+\int_0^t\mbb{E}\Bigl[|U_s^{i,p}|^2+\frac{1}{N_i}\sum_{p=1}^{N_i}|U_s^{i,p}|^2
+\frac{1}{N_j}\sum_{q=1}^{N_j}|U_s^{j,q}|^2+\bold{1}_{\{p=1\}}|\beta_s^{i,1}|^2\Bigr]ds\Bigr)~, \nn \\
&&\mbb{E}\Bigl[ |U_t^{j,q}|^2\Bigr]\leq C\Bigl(1+\int_0^t\mbb{E}\Bigl[|U_s^{j,q}|^2+\frac{1}{N_i}\sum_{p=1}^{N_i}|U_s^{i,p}|^2
+\frac{1}{N_j}\sum_{q=1}^{N_j}|U_s^{j,q}|^2\Bigr]ds\Bigr)~, \nn
\eea
Taking the average and the  applying Gronwall's inequality, we get
\be
\sup_{t\in[0,T]}\Bigl(\frac{1}{N_i}\sum_{p=1}^{N_i}\mbb{E}[|U_t^{i,p}|^2]+\frac{1}{N_j}
\sum_{q=1}^{N_j}\mbb{E}[|U_t^{j,q}|^2\bigr]\Bigr)\leq C\Bigl(1+\frac{1}{N_i}||\beta^{i,1}||^2_{\mbb{H}^2}\Bigr).~
\label{M2olnu}
\ee
Using the above estimate and Burkholder-Davis-Gundy (BDG) inequality, we get
\bea
\label{U-norm}
&&\mbb{E}\Bigl[\sup_{t\in[0,T]}|U_t^{i,p}|^2\Bigr]\leq C\Bigl(1+\bigl(\frac{1}{N_i}+\bold{1}_{\{p=1\}}\bigr)||\beta^{i,1}||^2_{\mbb{H}^2}\Bigr)~, \\
&&\mbb{E}\Bigl[\sup_{t\in[0,T]}|U_t^{j,q}|^2\Bigr]\leq C\Bigl(1+\frac{1}{N_i}||\beta^{i,1}||^2_{\mbb{H}^2}\Bigr)~.\nn
\eea
By similar calculation, we see
\bea
\mbb{E}\bigl[|U_t^{i,p}-\ul{X}_t^{i,p}|^2\bigr]\leq C\int_0^t \mbb{E}\Bigl[|U_s^{i,p}-\ul{X}_s^{i,p}|^2
+\sum_{j=1}^2 W_2(\ol{\nu}_s^j,\mu_s^j)^2+\bold{1}_{\{p=1\}}|\beta_s^{i,1}-\ul{\hat{\alpha}}_s^{i,1}|^2\Bigr]ds~. \nn
\eea
Combining the same estimate for the jth population, we get from Gronwall's inequality that
\bea
&&\sup_{t\in[0,T]}\Bigl(\frac{1}{N_i}\sum_{p=1}^{N_i}\mbb{E}\bigl[|U_t^{i,p}-\ul{X}_t^{i,p}|^2\bigr]+
\frac{1}{N_j}\sum_{p=1}^{N_j}\mbb{E}\bigl[|U_t^{j,q}-\ul{X}_t^{j,q}|^2\bigr]\Bigr) \nn \\
&&\quad \leq C\Bigl(\sum_{i=1}^2 \sup_{t\in[0,T]}\mbb{E}[W_2(\ul{\mu}_t^i,\mu_t^i)^2]+\frac{1}{N_i}\mbb{E}\int_0^T |\beta_t^{i,1}-
\ul{\hat{\alpha}}_t^{i,1}|^2 dt\Bigr)\nn \\
&&\quad \leq C\Bigl((1+||\beta^{i,1}||^2_{\mbb{H}^2})\vep_{N_i}^2+\vep_{N_j}^2\Bigr)~.\nn
\eea
Here, we have used  $||\ul{\hat{\alpha}}^{i,1}||^2_{\mbb{H}^2}\leq C$ and the result of Lemma~\ref{lemma-cl}.
By the triangle inequality, 
\bea
\sup_{t\in[0,T]} \sum_{i=1}^2\mbb{E}\bigl[W_2(\ol{\nu}_t^i,\mu_t^i)^2\bigr]
&\leq& 2\sup_{t\in[0,T]} \sum_{i=1}^2\Bigl(\mbb{E}\bigl[W_2(\ol{\nu}_t^i,\ul{\mu}_t^i)^2\bigr]+\mbb{E}\bigl[W_2(\ul{\mu}_t^i,\mu_t^i)^2\bigr]\Bigr)\nn \\
&\leq& C\Bigl((1+||\beta^{i,1}||^2_{\mbb{H}^2})\vep_{N_i}^2+\vep_{N_j}^2\Bigr)~
\label{olnu-mu-diff}
\eea
holds. Similarly, we  have
\bea
\mbb{E}\bigl[\sup_{s\in[0,t]}|U_s^{i,1}-\ol{U}_s^{i,1}|^2\bigr]\leq C\mbb{E}\int_0^t \Bigl[|U_s^{i,1}-\ol{U}_s^{i,1}|^2+W_2(\ol{\nu}_s^i,\mu_s^i)^2+W_2(\ol{\nu}_s^j,\mu_s^j)^2
\Bigr]ds \nn
\eea
and hence from $(\ref{olnu-mu-diff})$ we get
\bea
\mbb{E}\bigl[\sup_{t\in[0,T]}|U_t^{i,1}-\ol{U}_t^{i,1}|^2\bigr]^\frac{1}{2}\leq C\Bigl((1+||\beta^{i,1}||_{\mbb{H}^2})\vep_{N_i}+\vep_{N_j}\Bigr)~.
\label{U-olU}
\eea
\\
{\rm{\tbf{(third step)}}}
Finally, we get from the local Lipschitz continuity of the cost functions,
\bea
&&|J_i^{N_i,N_j,1}(\bg{\beta}^{i,1},\bg{\hat{\alpha}}^{i,(N_i)^{-1}},\bg{\hat{\alpha}}^{j,(N_j)})-\ol{J}_i^1(\bg{\beta}^{i,1})|\nn \\
&&=\Bigl|\mbb{E}\Bigl[\int_0^T\bigl(f_i(t,U_t^{i,1},\ol{\nu}_t^i,\ol{\nu}_t^j,\beta_t^{i,1})-
f_i(t,\ol{U}_t^{i,1},\mu_t^i,\mu_t^j,\beta_t^{i,1})\bigr)dt+g_i(U_T^{i,1},\ol{\nu}_T^i,\ol{\nu}_T^j)
-g_i(\ol{U}_T^{i,1},\mu_T^i,\mu_T^j)\Bigr]\Bigr|\nn\\
&&\leq C\Bigl(1+\sup_{t\in[0,T]}\mbb{E}\Bigl[|U_t^{i,1}|^2+|\ol{U}_t^{i,1}|^2+M_2(\ol{\nu}_t^i)^2+M_2(\ol{\nu}_t^j)^2\Bigr]^\frac{1}{2}
+||\beta^{i,1}||_{\mbb{H}^2}\Bigr) \nn \\
&&\qquad\qquad\times \Bigl(\sup_{t\in[0,T]}\mbb{E}\bigl[|U_t^{i,1}-\ol{U}_t^{i,1}|^2\bigr]^\frac{1}{2}+
\sum_{i=1}^2 \sup_{t\in[0,T]}\mbb{E}\bigl[W_2(\ol{\nu}_t^i,\mu_t^i)^2\bigr]^\frac{1}{2}\Bigr)~.\nn
\eea
From $(\ref{olU-norm})$,  $(\ref{M2olnu})$, $(\ref{U-norm})$, $(\ref{olnu-mu-diff})$ and $(\ref{U-olU})$, we get
\bea
&&|J_i^{N_i,N_j,1}(\bg{\beta}^{i,1},\bg{\hat{\alpha}}^{i,(N_i)^{-1}},\bg{\hat{\alpha}}^{j,(N_j)})-\ol{J}_i^1(\bg{\beta}^{i,1})| \nn \\
&&\qquad\qquad \leq C\bigl(1+||\beta^{i,1}||_{\mbb{H}^2}\bigr)\Bigl((1+||\beta^{i,1}||_{\mbb{H}^2})\vep_{N_i}+\vep_{N_j}\Bigr)~.\nn
\eea
By the estimate in Lemma~\ref{lemma-approx-limit}, $(\ref{olJ-J-estimate})$
and the fact that $||\ul{\hat{\alpha}}^{i,1}||^2_{\mbb{H}}\leq C$, we see
\bea
&&J_i^{N_i,N_j,1}(\bg{\beta}^{i,1},\bg{\hat{\alpha}}^{i,(N_i)^{-1}},\bg{\hat{\alpha}}^{j,(N_j)}) -J_i^{N_i,N_j,1}(\bg{\hat{\alpha}}^{i,(N_i)},\bg{\hat{\alpha}}^{j,(N_j)})\nn \\
&&\qquad \geq \lambda ||\beta^{i,1}-\ul{\hat{\alpha}}^{i,1}||^2_{\mbb{H}^2}
-C\bigl(1+||\beta^{i,1}||_{\mbb{H}^2}\bigr)\Bigl((1+||\beta^{i,1}||_{\mbb{H}^2})\vep_{N_i}+\vep_{N_j}\Bigr)~ \nn \\
&&\qquad \geq (\lambda-C\sum_{j=1}^2 \vep_{N_j})||\beta^{i,1}-\ul{\hat{\alpha}}^{i,1}||^2_{\mbb{H}^2}-C\sum_{j=1}^2\vep_{N_j}~.\nn
\eea
For large $N_1$ and $N_2$ with $C\sum_{j=1}^2 \vep_{N_j}\leq \lambda$, we get the desired result.
\end{proof}
\end{theorem}

\begin{remark}
The above analysis can be generalized straightforwardly to the setup with any finite number of populations, $1\leq i\leq m$.
\end{remark}

\section{Approximate Equilibrium for MFTC with Finite Agents}
\label{sec-mftc-fa}
In this section, we shall show how the solution to the Nash MFTC problem studied in Section~\ref{sec-mftc}
can provide an approximate Nash equilibrium among the two competing populations of finite agents who are cooperative
within each population.
In the last section dealing with the non-cooperative agents, the effect to the interactions 
from the agent deviating from the optimal strategy was shown to vanish in the large population limit.
This does not happen in the current case, because all the agents in one population adopt the common strategy different 
from the  optimal one. We shall see that this feature requires us more stringent assumptions to obtain an approximate Nash equilibrium.
Throughout the section, we assume that the conditions used either in Theorem~\ref{th-mftc-main} or Theorem~\ref{th-mftc-small-coupling} 
are satisfied.
We let $(\bg{\mu}^1,\bg{\mu}^2)\in\calc([0,T];\calp_2(\mbb{R}^d))^2$ denote a solution to the matching problem $(\ref{mftc-matching})$.
Moreover, unless otherwise stated, we use the same notation in the last section.
\subsection{Convergence of approximate optimal controls}
For each population $1\leq i\leq 2$, we first consider $(N_i)$ i.i.d. copies of the sate process in the mean field setup
\bea
d\ul{X}_t^{i,p}=b_i(t,\ul{X}_t^{i,p},\mu_t^i,\mu_t^j,\ul{\hat{\alpha}}_t^{i,p})dt+\sigma_i(t,\ul{X}_t^{i,p},\mu_t^i,\mu_t^j)dW_t^{i,p}, \nn
\eea
for $1\leq p\leq N_i$, $j\neq i$, $t\in[0,T]$ with $\ul{X}_0^{i,p}=\xi^{i,p}$, and 
\bea
\ul{\hat{\alpha}}_t^{i,p}:=\hat{\alpha}_i(t,\ul{X}_t^{i,p},\mu_t^i,\mu_t^j, u_i(t,\ul{X}_t^{i,p},\mu_t^i))~, \nn
\eea
where $u_i$ is the function defined in Lemma~\ref{lemma-Lip-mftc} associated with the equilibrium flows of 
probability measures $(\bg{\mu}^i,\bg{\mu}^j)$.
As in the last section, $(\xi^{i,p},\bg{W}^{i,p})_{1\leq p\leq N_i, 1\leq i\leq 2}$ are assumed to be independent
with $\call(\xi^{i,p})=\mu_0^i$. 
By construction, $(\bg{\ul{X}}^{i,p})_{1\leq p\leq N_i}$ are i.i.d. processes satisfying $\call(\ul{X}_t^{i,p})=\mu_t^i$, $\forall t\in[0,T]$. 
$\ul{\mu}_t^i$ denotes the empirical distribution of $(\ul{X}_t^{i,p})_{1\leq p\leq N_i}$. 
\begin{lemma}
\label{lemma-cl-mftc}
Suppose that the conditions either for Theorem~\ref{th-mftc-main} or Theorem~\ref{th-mftc-small-coupling} are satisfied.
Then, for each population $1\leq i\leq 2$, there exists some sequence $(\ep_{N_i})_{N_i\geq 1}$ that tends to $0$ as $N_i$ tends to $\infty$
and some constant $C$ such that
\bea
\sup_{t\in[0,T]}\mbb{E}\Bigl[W_2(\ul{\mu}_t^i,\mu_t^i)^2\Bigr]\leq C\ep_{N_i}^2~.\nn
\eea
Furthermore, when $\mu_0^i\in \calp_r(\mbb{R}^d)$ with $r>4$, we have
an explicit estimate 
\bea
\ep_{N_i}^2=N_i^{-2/\max(d,4)}(1+\ln(N_i)\bold{1}_{d=4})~.\nn
\eea
\begin{proof}
It can be proved in the same way as Lemma~\ref{lemma-cl}.
\end{proof}
\end{lemma}
Let us introduce the following assumptions.
\begin{assumption}{\rm{\tbf{(MFTC-FA-a)}}}
On top of Assumption {\rm{\tbf{(MFTC-a)}}}, either $T$ or $(\lambda^{-1}||b_{i,2}||_{\infty})_{1\leq i\leq 2}$
is small enough to satisfy the conditions for Theorem~\ref{th-mftc-main} or Theorem~\ref{th-mftc-small-coupling}.
Moreover, for $1\leq i\leq 2$, \\
{\rm{\tbf{(A1)}}} There exists some constant $K$ such that
\be
|(b_{i,0},\sigma_{i,0})(t,\nu^\prime)-(b_{i,0},\sigma_{i,0})(t,\nu)|\leq KW_2(\nu^\prime, \nu) \nn
\ee
for any $t\in[0,T]$, $\nu^\prime, \nu\in \calp_2(\mbb{R}^d)$, and $b_{i,1},\bar{b}_{i,1},b_{i,2}, \sigma_{i,1}$ as well as $\bar{\sigma}_{i,1}$
are independent of the measure argument. \\
{\rm{\tbf{(A2)}}}  $f_i$ and $g_i$ are local Lipschitz continuous with respect to the second measure argument i.e.,
there exists some constant $K$ for any $t\in[0,T],x\in \mbb{R}^d, \mu,\nu^\prime, \nu\in \calp_2(\mbb{R}^d)$ and $\alpha\in A_i$, 
such that
\bea
&&|(f_i,g_i)(t,x,\mu,\nu^\prime,\alpha)-(f_i,g_i)(t,x,\mu,\nu,\alpha)|\nn \\
&&\quad \leq K\Bigl(1+|x|+M_2(\mu)+M_2(\nu^\prime)+M_2(\nu)+|\alpha|\Bigr)W_2(\nu^\prime,\nu)~.\nn
\eea
\end{assumption}
\begin{assumption}{\rm{\tbf{(MFTC-FA-b)}}}
On top of Assumption {\rm{\tbf{(MFTC-FA-a)}}}, for $1\leq i\leq 2$, $b_{i,0}$ and $\sigma_{i,0}$
are independent of the measure argument.
\end{assumption}
\begin{remark}
Assumption {\rm{\tbf{(MFTC-FA-b)}}} will be used in the last part where we prove the property of the approximate Nash equilibrium.
Under this stringent assumption, the mutual interactions among the agents belonging to the 
different populations are induced only through the cost functions and can appear only in their control strategies.
\end{remark}

As in the last section, we introduce for $1\leq i,j\leq 2$, $j\neq i$, $1\leq p\leq N_i$ the state dynamics
\bea
dX_t^{i,p}=b_i(t,X_t^{i,p},\ol{\mu}_t^i,\ol{\mu}_t^j,\hat{\alpha}_t^{i,p})dt+\sigma_i(t,X_t^{i,p},\ol{\mu}_t^i,\ol{\mu}_t^j)dW_t^{i,p}~, \nn
\eea
for $t\in[0,T]$ with $X_0^{i,p}=\xi^{i,p}$. Here, $\ol{\mu}_t^i:=\frac{1}{N_i}\sum_{p=1}^{N_i}\del_{X_t^{i,p}}$ and 
\bea
\hat{\alpha}_t^{i,p}:=\hat{\alpha}_i(t,X_t^{i,p},\mu_t^i,\mu_t^j,u_i(t,X_t^{i,p},\mu_t^i))~.\nn
\eea
Under Assumption {\rm{\tbf{(MFTC-FA-a)}}}, it is an $(N_1+N_2)$-dimensional Lipschitz SDE and hence is 
well-defined. This corresponds to the situation where all the agents in each population adopt
the common feedback control given by the solution to the problem in Section~\ref{sec-mftc}.
Let us write the cost functional for the agent in the ith population as
\bea
J_i^{N_i,N_j}(\bg{\hat{\alpha}}^{i,(N_i)},\bg{\hat{\alpha}}^{j,(N_j)}):=
\mbb{E}\Bigl[\int_0^T f_i(t,X_t^{i,p},\ol{\mu}_t^i,\ol{\mu}_t^j,\hat{\alpha}_t^{i,p})dt+g_i(X_T^{i,p},\ol{\mu}_T^i,\ol{\mu}_T^j)\Bigr]~.\nn
\eea
The corresponding cost functional in the mean field problem is given by
\bea
J_i:=\mbb{E}\Bigl[\int_0^T f_i(t,\ul{X}_t^{i,p},\mu_t^i,\mu_t^j,\ul{\hat{\alpha}}_t^{i,p})dt+g_i(\ul{X}_T^{i,p},\mu_T^i,\mu_T^j)\Bigr]~.\nn
\eea
where $1\leq p\leq N_i$ is arbitrary in both cases. With the word cooperative, we mean that the agents use the
common feedback control function. Hence, even when we consider general strategy later, all the cost functionals among the agents within 
each population are the same.

\begin{lemma}
\label{lemma-approx-limit-mftc}
Under Assumption {\rm{\tbf{(MFTC-FA-a)}}}, for all $1\leq i\leq 2$, $N_i\in \mbb{N}$, 
there exists some constant $C$ independent of $(N_i)_{i=1}^2$ such that,
\bea
|J_i^{N_i,N_j}(\bg{\hat{\alpha}}^{i,(N_i)},\bg{\hat{\alpha}}^{j,(N_j)})-J_i|\leq C\sum_{j=1}^2\ep_{N_j}~, \nn
\eea
where $\ep_{N_j}$ is the one given in Lemma~\ref{lemma-cl-mftc}.
\begin{proof}
It is straightforward to get the estimate
\bea
\sum_{i=1}^2\mbb{E}\bigl[\sup_{t\in[0,T]}|X_t^{i,p}-\ul{X}_t^{i,p}|^2\bigr]\leq C
\sum_{i=1}^2\sup_{t\in[0,T]}\mbb{E}\bigl[W_2(\ol{\mu}_t^i,\mu_t^i)^2\bigr]\leq C\sum_{i=1}^2\ep_{N_i}^2~,\nn
\eea
and  $\mbb{E}\bigl[\sup_{t\in[0,T]}|X_t^{i,p}|^2\bigr]+\mbb{E}\bigl[\sup_{t\in[0,T]}|\ul{X}_t^{i,p}|^2\bigr]\leq C$.
Using the local Lipschitz continuity for $f_i$ and $g_i$, we can prove the convergence of the cost functional exactly in the same way as in Lemma~\ref{lemma-approx-limit}.
\end{proof}
\end{lemma}

\subsection{Approximate Nash Equilibrium}
We now consider the general state dynamics for the agents $1\leq i, j\leq 2$, $j\neq i$, $1\leq p\leq N_i$,
\bea
dU_t^{i,p}=b_i(t,U_t^{i,p}, \ol{\nu}_t^i,\ol{\nu}_t^j,\beta_t^{i,p})dt+\sigma_i(t,U_t^{i,p},\ol{\nu}_t^i,\ol{\nu}_t^j)dW_t^{i,p} \nn
\eea
with $t\in[0,T]$, $U_0^{i,p}=\xi^{i,p}$ and $\bg{\beta}^{i,p} \in \mbb{H}^2$ is an $A_i$-valued $\mbb{F}$-progressively
measurable process.
Since we suppose that the agents are cooperative within each population,
we force the set of strategies $(\bg{\beta}^{i,p})_{1\leq p\leq N_i}$ to satisfy the 
condition so that $(\xi^{i,p},\bg{\beta}^{i,p},\bg{W}^{i,p})_{1\leq  p\leq N_i}$ is exchangeable i.e. the distribution is invariant
under the permutation with $p$. As before $\ol{\nu}_t^i:=\frac{1}{N_i}\sum_{p=1}^{N_i}\del_{U^{i,p}_t}$ denotes the 
empirical distribution.
The associated cost functional for the ith population is now given by, with any $1\leq p\leq N_i$,
\bea
J_i^{N_i,N_j}(\bg{\beta}^{i,(N_i)},\bg{\beta}^{j,(N_j)}):=\mbb{E}\Bigl[\int_0^T f_i(t,U_t^{i,p},\ol{\nu}_t^i,\ol{\nu}_t^j,\beta_t^{i,p})dt+
g_i(U_T^{i,p},\ol{\nu}_T^i,\ol{\nu}_T^j)\Bigr]~.\nn
\eea
Let us also introduce $(\ul{Y}^{i,p}_t,\ul{Z}^{i,p}_t)_{t\in[0,T]}$ the solution to $(\ref{mkv-mftc-master})$
associated with the forward component $(\ul{X}_t^{i,p})_{t\in[0,T]}$:
\bea
d\ul{Y}_t^{i,p}&=&-\part_x H_i(t,\ul{X}_t^{i,p},\mu_t^i,\mu_t^j,\ul{Y}_t^{i,p},\ul{Z}_t^{i,p},\hat{\alpha}_i(t,\ul{X}_t^{i,p},\mu_t^i,\mu_t^j,
\ul{Y}_t^{i,p}))dt\nn \\
&&-\wt{\mbb{E}}\bigl[\part_\mu H_i(t,\wt{\ul{X}}_t^{i,p},\mu_t^i,\mu_t^j,\wt{\ul{Y}}_t^{i,p},\wt{\ul{Z}}_t^{i,p},
\hat{\alpha}_i(t,\wt{\ul{X}}_t^{i,p},\mu_t^i,\mu_t^j,\wt{\ul{Y}}_t^{i,p}))(\ul{X}_t^{i,p})\bigr]dt+\ul{Z}_t^{i,p}dW_t^{i,p},
\label{bar-bsde}
\eea
with $\ul{Y}_T^{i,p}=\part_x g_i(\ul{X}_T^{i,p},\mu_T^i,\mu_T^j)+\wt{\mbb{E}}\bigl[\part_\mu g_i(\wt{\ul{X}}_T^{i,p},\mu_T^i,\mu_T^j)(\ul{X}_T^{i,p})]$.
Note that, $\ul{Y}_t^{i,p}=u_i(t,\ul{X}_t^{i,p},\mu_t^i)$ a.s. for any $t\in[0,T]$.

\begin{proposition}
\label{prop-mftc-fa}
Under Assumption {\rm{\tbf{(MFTC-FA-a)}}}, for $1\leq i\leq 2$, $N_i\in \mbb{N}$, 
there exists some constant $C$ independent of $(N_i)_{i=1}^2$ such that
\bea
&&J_i^{N_i,N_j}(\bg{\beta}^{i,(N_i)},\bg{\beta}^{j,(N_j)})-J_i \geq \lambda \mbb{E}\int_0^T |\beta_t^{i,p}-\ul{\hat{\alpha}}_t^{i,p}|^2dt\nn \\
&&\qquad -C\Bigl(1+\sup_{t\in[0,T]}\mbb{E}\bigl[|U_t^{i,p}-\ul{X}_t^{i,p}|^2\bigr]^\frac{1}{2}
+\mbb{E}\bigl[\int_0^T |\beta_t^{i,p}-\ul{\hat{\alpha}}_t^{i,p}|^2 dt\bigr]^\frac{1}{2}\Bigr)\vep_{N_i}\nn \\
&&\qquad+\mbb{E}\Bigl[\int_0^T\bigl( H_i(t,U_t^{i,p},\ol{\nu}_t^i,\ol{\nu}_t^j,\ul{Y}_t^{i,p},\ul{Z}_t^{i,p},\beta_t^{i,p})
-H_i(t,U_t^{i,p},\ol{\nu}_t^i,\mu_t^j, \ul{Y}_t^{i,p},\ul{Z}_t^{i,p},\beta_t^{i,p})\bigr)dt\nn \\
&&\qquad\qquad+g_i(U_T^{i,p},\ol{\nu}_T^i,\ol{\nu}_T^j)-g_i(U_T^{i,p},\ol{\nu}_T^i,\mu_T^j)\Bigr]~,\nn
\eea
where $\vep_{N_i}:=\max(N_i^{-\frac{1}{2}},\ep_{N_i})$, and $1\leq p\leq N_i$ is arbitrary. 
\begin{proof}
We can show the claim by following the same arguments used in the proof for Theorem 6.16 \cite{Carmona-Delarue-2}.
Since it is rather technical and lengthy, we give the details in Appendix~\ref{app-prop-mftc-fa}.
\end{proof}
\end{proposition}

For investigating the approximate Nash equilibrium property, we now suppose that the agents in the ith population 
use general strategy $(\bg{\beta}^{i,p})_{1\leq p\leq N_i}$ under the restriction that $(\xi^{i,p},\bg{\beta}^{i,p},\bg{W}^{i,p})_{1\leq p\leq N_i}$
is exchangeable, and that the agents in the jth population $1\leq q\leq N_j$ adopt the strategy
\bea
\beta_t^{j,q}:=\hat{\alpha}_j(t,U_t^{j,q},\mu_t^j,\mu_t^i,u_j(t,U_t^{j,q},\mu_t^j))\nn
\eea
for any $t\in[0,T]$.  The cost functional for the ith population is now given by
\bea
J_i^{N_i,N_j}(\bg{\beta}^{i,(N_i)},\bg{\hat{\alpha}}^{j,(N_j)})
=\mbb{E}\Bigl[\int_0^T f_i(t,U_t^{i,p},\ol{\nu}_t^i,\ol{\nu}_t^j,\beta_t^{i,p})dt+g_i(U_T^{i,p},\ol{\nu}_T^i,\ol{\nu}_T^j)\Bigr]\nn
\eea
with the above specified control strategies.
We now proceed as in the proof for Theorem~\ref{th-nash-mfg}. The crucial problem
is the term $\mbb{E}\bigl[W_2(\ol{\nu}_t^j,\mu_t^j)^2\bigr]^\frac{1}{2}$ arising from the last line in the estimate of Proposition~\ref{prop-mftc-fa}.
Although this term is suppressed in the non-cooperative game  as in $(\ref{olnu-mu-diff})$, 
it does not happen in the current situation.
Under Assumption {\rm{\tbf{(MFTC-FA-a)}}}, the deviation from the strategy $\hat{\alpha}_i$ for the agents in the ith population induces the term 
$||\beta^{i,1}-\ul{\hat{\alpha}}_i||_{\mbb{H}^2}$ with no suppression of $\vep_{N_i}$.
This is why we need Assumption {\rm{\tbf{(MFTC-FA-b)}}}.

\begin{theorem}
\label{th-mftc-fa}
Under Assumption {\rm{\tbf{(MFTC-FA-b)}}} with sufficiently large $(N_i)_{i=1}^2$, 
the feedback control functions $\bigl([0,T]\times \mbb{R}^d\ni (t,x)\mapsto \hat{\alpha}_i(t,x,\mu_t^i,\mu_t^j,u_i(t,x,\mu_t^i))\bigr)_{1\leq i, j\leq 2, j\neq i}$ form an $(\sum_{j=1}^2\vep_{N_j})$-approximate Nash equilibrium i.e., there exists some constant $C$ independent of 
$(N_i)_{i=1}^2$ such that
\bea
J_i^{N_i,N_j}(\bg{\beta}^{i,(N_i)},\bg{\hat{\alpha}}^{j,(N_j)})\geq J_i^{N_i,N_j}(\bg{\hat{\alpha}}^{i,(N_i)},\bg{\hat{\alpha}}^{j,(N_j)})
-C\sum_{j=1}^2\vep_{N_j}\nn
\eea
for any $A_i$-valued $\mbb{F}$-progressively measurable processes $(\bg{\beta}^{i,p}\in \mbb{H}^2)_{1\leq p\leq N_i}$
so that $(\xi^{i,p},\bg{\beta}^{i,p},\bg{W}^{i,p})_{1\leq p\leq N_i}$ is exchangeable. Here, $\bigl(\vep_{N_j}:=\max(N_j^{-\frac{1}{2}},\ep_{N_j})\bigr)_{1\leq j\leq 2}$.
\begin{proof}
Under {\rm{\tbf{(MFTC-FA-b)}}}, we can write the coefficients for both populations $1\leq j\leq 2$ as,
\bea
b_j(t,x,\mu,\alpha)&=&b_{j,0}(t)+b_{j,1}(t)x+\bar{b}_{j,1}(t)\bar{\mu}+b_{j,2}(t)\alpha\nn \\
\sigma_j(t,x,\mu)&=&\sigma_{j,0}(t)+\sigma_{j,1}(t)x+\bar{\sigma}_{j,1}(t)\bar{\mu}~.
\label{bsigma-FA-b}
\eea
For the ith population, we get
\bea
\mbb{E}\bigl[\sup_{s\in[0,t]}|U_s^{i,p}|^2\bigr]\leq C\mbb{E}\Bigl[|\xi^{i,p}|^2+\int_0^t\bigl(1+|U_s^{i,p}|^2+M_2(\ol{\nu}_s^i)^2+|\beta_s^{i,p}|^2\bigr)ds
\Bigr]~,\nn 
\eea
which yields $\mbb{E}\bigl[\sup_{t\in[0,T]}|U_t^{i,p}|^2\bigr]\leq C\bigl(1+||\beta^{i,p}||^2_{\mbb{H}^2}\bigr)$.
We also get $\mbb{E}\bigl[\sup_{t\in[0,T]}|U_t^{j,q}|^2\bigr]\leq C$ for the jth population.
Similarly,  we see for the ith population,
\bea
\mbb{E}\bigl[\sup_{s\in[0,t]}|U_s^{i,p}-\ul{X}_s^{i,p}|^2\bigr]
\leq C\mbb{E}\int_0^t\Bigl[|U_s^{i,p}-\ul{X}_s^{i,p}|^2+W_2(\ol{\nu}_s^i,\ul{\mu}_s^i)^2+W_2(\ul{\mu}_s^i,\mu_s^i)^2+
|\beta_s^{i,p}-\ul{\hat{\alpha}}_s^{i,p}|^2\Bigr]ds, \nn
\eea
and then, by Gronwall's inequality, we get
\bea
\mbb{E}\bigl[\sup_{t \in[0,T]}|U_t^{i,p}-\ul{X}_t^{i,p}|^2\bigr]&\leq& C\Bigl(\sup_{t\in[0,T]}\mbb{E}\bigl[W_2(\ul{\mu}_t^i,\mu_t^i)^2\bigr]
+||\beta^{i,p}-\ul{\hat{\alpha}}^{i,p}||^2_{\mbb{H}^2}\Bigr)~\nn \\
&\leq &C\bigl(\vep_{N_i}^2+||\beta^{i,p}-\ul{\hat{\alpha}}^{i,p}||^2_{\mbb{H}^2}\bigr)~\nn.
\eea
Similarly,  we have
\bea
\mbb{E}\bigl[\sup_{t\in[0,T]}|U_t^{j,q}-\ul{X}_t^{j,q}|^2\bigr]\leq C\Bigl(\sup_{t\in[0,T]}\mbb{E}\bigl[W_2(\ul{\mu}_t^j,\mu_t^j)^2\bigr]\Bigr)
\leq C\vep_{N_j}^2~.\nn
\eea
In particular, by the triangle inequality, $\sup_{t\in[0,T]}\mbb{E}\bigl[W_2(\ol{\nu}_t^j,\mu_t^j)^2\bigr]\leq C\vep_{N_j}^2$ holds.

From $(\ref{bsigma-FA-b})$, it is easy to see 
\bea
&&|H_i(t,U_t^{i,p},\ol{\nu}_t^i,\ol{\nu}_t^j,\ul{Y}_t^{i,p},\ul{Z}_t^{i,p},\beta_t^{i,p})
-H_i(t,U_t^{i,p},\ol{\nu}_t^i,\mu_t^j, \ul{Y}_t^{i,p},\ul{Z}_t^{i,p},\beta_t^{i,p})| \nn \\
&&=|f_i(t,U_t^{i,p},\ol{\nu}_t^i,\ol{\nu}_t^j,\beta_t^{i,p})-f_i(t,U_t^{i,p},\ol{\nu}_t^i,\mu_t^j,\beta_t^{i,p})| \nn \\
&&\leq C(1+|U_t^{i,p}|+|M_2(\ol{\nu}_t^i)|+|M_2(\ol{\nu}_t^j)|+|\beta_t^{i,p}|)W_2(\ol{\nu}_t^j,\mu_t^j)
\eea
Using similar estimate for $g_i$ and exchangeability, we get from Proposition~\ref{prop-mftc-fa},
\bea
&&J_i^{N_i,N_j}(\bg{\beta}^{i,(N_i)},\bg{\hat{\alpha}}^{j,(N_j)})-J_i \nn \\
&&\geq  \lambda \mbb{E}\int_0^T |\beta_t^{i,p}-\ul{\hat{\alpha}}_t^{i,p}|^2dt 
-C\bigl(1+||\beta^{i,p}-\ul{\hat{\alpha}}^{i,p}||_{\mbb{H}^2} \bigr)\vep_{N_i}\nn \\
&&\quad -C\Bigl(1+\sup_{t\in[0,T]}\mbb{E}\bigl[|U_t^{i,p}|^2+|U_t^{j,q}|^2\bigr]^\frac{1}{2}+||\beta^{i,p}||_{\mbb{H}^2}\Bigr)
\sup_{t\in[0,T]}\mbb{E}\bigl[W_2(\ol{\nu}_t^j,\mu_t^j)^2\bigr]^\frac{1}{2}\nn \\
&&\geq  \lambda \mbb{E}\int_0^T |\beta_t^{i,p}-\ul{\hat{\alpha}}_t^{i,p}|^2dt 
-C\bigl(1+||\beta^{i,p}-\ul{\hat{\alpha}}^{i,p}||_{\mbb{H}^2} \bigr)\vep_{N_i} -C\Bigl(1+||\beta^{i,p}||_{\mbb{H}^2}\Bigr)\vep_{N_j}~.\nn
\eea
Since $||\ul{\hat{\alpha}}^{i,p}||_{\mbb{H}^2}\leq C$, and using the fact that
\bea
||\beta^{i,p}-\ul{\hat{\alpha}}^{i,p}||_{\mbb{H}^2}\vep_{N_i}\leq C\bigl(||\beta^{i,p}-\ul{\hat{\alpha}}^{i,p}||_{\mbb{H}^2}^2\vep_{N_i}+\vep_{N_i}\bigr)\nn
\eea
we get
\bea
J_i^{N_i,N_j}(\bg{\beta}^{i,(N_i)},\bg{\hat{\alpha}}^{j,(N_j)})-J_i \geq \bigl(\lambda-C\sum_{j=1}^2\vep_{N_j}\bigr)||\beta^{i,p}-
\ul{\hat{\alpha}}^{i,p}||^2_{\mbb{H}^2}-C\sum_{j=1}^2\vep_{N_j}~.\nn
\eea
We now get the desired estimate from Lemma~\ref{lemma-approx-limit-mftc} for sufficiently large $N_1$ and $N_2$.
\end{proof}
\end{theorem}

\begin{remark}
We have investigated the approximate Nash equilibrium in the closed loop framework.  The analysis for the open loop 
framework can be done in a quite similar manner as explained in Remark~\ref{remark-open}. 
Generalization to an arbitrary number of populations $1\leq i\leq m$
can be done straightforwardly.
\end{remark}

\section{Approximate Equilibrium for MFTC-MFG with Finite Agents}
\label{sec-mftc-mfg-fa}
In this section, we shall see how the solution to the Nash MFTC-MFG problem studied in Section~\ref{sec-mftc-mfg} can provide
an approximate Nash equilibrium among the two competing populations of finite agents, where the agents in the first population
are cooperative but those in the second population are not.
As we have seen in the last section, the effect of deviation from the optimal strategy in the first population 
will not be suppressed by $\vep_{N_1}$. In order to obtain an approximate Nash equilibrium, this feature
implies that we have to cut the direct interaction with the first population in the state dynamics of the second one.
On the other hand, the agents in the second population are non-cooperative,
and hence, the effect of the deviation of a single agent will be suppressed by $\vep_{N_2}$.
This suggests that we may include the direct interaction with the second population in  the state dynamics of the first one.
Throughout the section, we assume that the conditions used either in Theorem~\ref{th-mftc-mfg-main} or Theorem~\ref{th-mftc-mfg-small-coupling} 
are satisfied.
We let $(\bg{\mu}^1,\bg{\mu}^2)\in\calc([0,T];\calp_2(\mbb{R}^d))^2$ denote a solution to the matching problem $(\ref{mftc-mfg-matching})$.
Moreover, unless otherwise stated, we use the same notation in the last two sections.

\subsection{Convergence of approximate optimal controls}
For each population $1\leq i \leq 2$,  we give i.i.d. copies of the state process in the mean field setup:
\bea
d\ul{X}_t^{i,p}=b_i(t,\ul{X}_t^{i,p},\mu_t^i,\mu_t^j, \ul{\hat{\alpha}}_t^{i,p})dt+\sigma_i(t,\ul{X}_t^{i,p},\mu_t^i,\mu_t^j)dW_t^{i,p}~, \nn
\eea
for $1\leq p\leq N_i$, $j\neq i$, $t\in[0,T]$ with $X_0^{i,p}=\xi^{i,p}$, and 
\bea
&&\ul{\hat{\alpha}}^{1,p}:=\hat{\alpha}_1(t,\ul{X}_t^{1,p},\mu_t^1,\mu_t^2,u_1(t,\ul{X}_t^{1,p},\mu_t^1))~, \nn \\
&&\ul{\hat{\alpha}}^{2,p}:=\hat{\alpha}_2(t,\ul{X}_t^{2,p},\mu_t^2,\mu_t^1,u_2(t,\ul{X}_t^{2,p}))~,\nn
\eea
where $u_1$ is the master field defined in Lemma~\ref{lemma-Lip-mftc} applied to $(\ref{p1-fbsde})$
and $u_2$ the decoupling field defined in Theorem~\ref{th-fbsde-existence} applied to $(\ref{p2-fbsde})$
where the equilibrium flow of probability measures $(\bg{\mu}^1,\bg{\mu}^2)$ are used as inputs 
in both cases. As before, $(\xi^{i,p},\bg{W}^{i,p})_{1\leq p\leq N_i, 1\leq i\leq 2}$ are assumed to be independent
with $\call(\xi^{i,p})=\mu_0^i$. 
By construction, $(\bg{\ul{X}}^{i,p})_{1\leq p\leq N_i}$ are i.i.d. processes satisfying $\call(\ul{X}_t^{i,p})=\mu_t^i$, $\forall t\in[0,T]$. 
$\ul{\mu}_t^i:=\frac{1}{N_i}\sum_{p=1}^{N_i}\del{\ul{X}_t^{i,p}}$ denotes the empirical distribution of $(\ul{X}_t^{i,p})_{1\leq p\leq N_i}$.

\begin{lemma}
\label{lemma-cl-mftc-mfg}
Suppose that the conditions either for Theorem~\ref{th-mftc-mfg-main} or Theorem~\ref{th-mftc-mfg-small-coupling} are satisfied.
Then, for each population $1\leq i\leq 2$, there exists some sequence $(\ep_{N_i})_{N_i\geq 1}$ that tends to $0$ as $N_i$ tends to $\infty$
and some constant $C$ such that
\bea
\sup_{t\in[0,T]}\mbb{E}\Bigl[W_2(\ul{\mu}_t^i,\mu_t^i)^2\Bigr]\leq C\ep_{N_i}^2~.\nn
\eea
Furthermore, when $\mu_0^i\in \calp_r(\mbb{R}^d)$ with $r>4$, we have
an explicit estimate 
\bea
\ep_{N_i}^2=N_i^{-2/\max(d,4)}(1+\ln(N_i)\bold{1}_{d=4})~.\nn
\eea
\begin{proof}
It is the direct result of Lemma~\ref{lemma-cl} and Lemma~\ref{lemma-cl-mftc}.
\end{proof}
\end{lemma}
We introduce the following assumptions.
\begin{assumption}{\rm{\tbf{(MFTC-MFG-FA-a)}}} On top of Assumption {\rm{\tbf{(MFTC-a)}}}
for the coefficients $(b_1, \sigma_1,f_1,g_1)$ and Assumption {\rm{\tbf{(MFG-a)}}} 
for the coefficients $(b_2,\sigma_2,f_2,g_2)$, either $T$ or $(\lambda^{-1}||b_{i,2}||_{\infty})_{1\leq i\leq 2}$
is small enough to satisfy the conditions for Theorem~\ref{th-mftc-mfg-main} or Theorem~\ref{th-mftc-mfg-small-coupling}.
Moreover,  the coefficients $(b_1,\sigma_1,f_1,g_1)$ satisfy  {\rm{\tbf{(MFTC-FA-a) (A1-A2)}}},
and $(b_2,\sigma_2,f_2,g_2)$  satisfy
{\rm{\tbf{(MFG-FA) (A1-A2)}}}.
\end{assumption}
\begin{assumption}{\rm{\tbf{(MFTC-MFG-FA-b)}}}
On top of Assumption {\rm{\tbf{(MFTC-MFG-FA-a)}}}, the coefficients $b_{2,0}$ and $\sigma_{2,0}$
are independent of the second measure argument, i.e. $(b_{2,0},\sigma_{2,0})(t,\mu,\nu)=(b_{2,0},\sigma_{2,0})(t,\mu)$.
\end{assumption}

For $1\leq i,j\leq 2$, $j\neq i$, $1\leq p\leq N_i$, we introduce the state processes
\bea
dX_t^{i,p}=b_i(t,X_t^{i,p},\ol{\mu}_t^i,\ol{\mu}_t^j,\hat{\alpha}_t^{i,p})dt+\sigma_i(t,X_t^{i,p},\ol{\mu}_t^i,\ol{\mu}_t^j)dW_t^{i,p}~, \nn
\eea
for $t\in[0,T]$ with $X_0^{i,p}=\xi^{i,p}$. Here, $\ol{\mu}_t^i:=\frac{1}{N_i}\sum_{p=1}^{N_i}\del_{X_t^{i,p}}$ denotes the empirical
distribution and
\bea
&&\hat{\alpha}_t^{1,p}:=\hat{\alpha}_1(t,X_t^{1,p},\mu_t^1,\mu_t^2,u_1(t,X_t^{1,p},\mu_t^1))~, \nn\\
&&\hat{\alpha}_t^{2,p}:=\hat{\alpha}_2(t,X_t^{2,p},\mu_t^2,\mu_t^1,u_2(t,X_t^{2,p}))~.\nn
\eea
Under Assumption {\rm{\tbf{(MFTC-MFG-FA-a)}}}, it is an $(N_1+N_2)$-dimensional Lipschitz SDE and hence is 
well-defined. The corresponding cost functional for any agent $p$ in the first population is given by
\bea
J_1^{N_1,N_2}(\bg{\hat{\alpha}}^{1,(N_1)},\bg{\hat{\alpha}}^{2,(N_2)}):=\mbb{E}\Bigl[
\int_0^T f_1(t,X_t^{1,p},\ol{\mu}_t^1,\ol{\mu}_t^2,\hat{\alpha}_t^{1,p})dt+g_1(X_T^{1,p},\ol{\mu}_T^1,\ol{\mu}_T^2)\Bigr]~,\nn
\eea
and the for the agent $q$ in the second population,
\bea
J_2^{N_2,N_1,q}(\bg{\hat{\alpha}}^{2,(N_2)},\bg{\hat{\alpha}}^{1,(N_1)}):=
\mbb{E}\Bigl[\int_0^T f_2(t,X_t^{2,q},\ol{\mu}_t^2,\ol{\mu}_t^1,\hat{\alpha}_t^{2,q})dt+g_2(X_T^{2,q},\ol{\mu}_T^2,\ol{\mu}_T^1)\Bigr]~.\nn
\eea
We also introduce the optimal cost functionals in the mean field setup:
\bea
&&J_1:=\mbb{E}\Bigl[\int_0^T f_1(t,\ul{X}_t^{1,p},\mu_t^1,\mu_t^2,\ul{\hat{\alpha}}_t^{1,p})dt+g_1(\ul{X}_T^{1,p},\mu_T^1,\mu_T^2)\Bigr]~, \nn \\
&&J_2^q:=\mbb{E}\Bigl[\int_0^T f_2(t,\ul{X}_t^{2,q},\mu_t^2,\mu_t^1,\ul{\hat{\alpha}}_t^{2,q})dt+g_2(\ul{X}_T^{2,q},\mu_T^2,\mu_T^1)\Bigr]~.\nn
\eea
\begin{remark}
Under the given control strategy, the value of $J_2^q$ is independent of $q$. However, since each agent in the second population can choose
his/her own strategy in general, we need to specify the agent when we discuss the approximate Nash equilibrium later. 
This is why we keep the index $q$ in the cost functional.
\end{remark}

\begin{lemma}
\label{lemma-approx-mftc-mfg}
Under Assumption~{\rm{\tbf{(MFTC-MFG-FA-a)}}}, for any $N_1, N_2\in \mbb{N}$ and $1\leq q\leq N_2$,
there exists some constant C independent of $(N_i)_{i=1}^2$ such that,
\bea
&&|J_1^{N_1,N_2}(\bg{\hat{\alpha}}^{1,(N_1)},\bg{\hat{\alpha}}^{2,(N_2)})-J_1|\leq C\sum_{j=1}^2 \ep_{N_j},~~|J_2^{N_2,N_1,q}(\bg{\hat{\alpha}}^{2,(N_2)},\bg{\hat{\alpha}}^{1,(N_1)})-J_2^q|\leq C\sum_{j=1}^2 \ep_{N_j}~ \nn
\eea
where $\ep_{N_j}$ is the one given in Lemma~\ref{lemma-cl-mftc-mfg}. 
\begin{proof}
As in the last two sections, we can show, by the same arguments, 
\bea
\sum_{i=1}^2\mbb{E}\bigl[\sup_{t\in[0,T]}|X_t^{i,p}-\ul{X}_t^{i,p}|^2\bigr]\leq C
\sum_{i=1}^2\sup_{t\in[0,T]}\mbb{E}\bigl[W_2(\ol{\mu}_t^i,\mu_t^i)^2\bigr]\leq C\sum_{i=1}^2\ep_{N_i}^2~,\nn
\eea
and  $\mbb{E}\bigl[\sup_{t\in[0,T]}|X_t^{i,p}|^2\bigr]+\mbb{E}\bigl[\sup_{t\in[0,T]}|\ul{X}_t^{i,p}|^2\bigr]\leq C$.
Thus the convergence of the cost functionals is the direct result of Lemma~\ref{lemma-approx-limit}
and Lemma~\ref{lemma-approx-limit-mftc}.
\end{proof}
\end{lemma}

\subsection{Approximate Nash Equilibrium}
We now consider the general state dynamics for the agents $1\leq i,j \leq 2$, $j\neq i$, $1\leq p\leq N_i$,
\bea
dU_t^{i,p}=b_i(t,U_t^{i,p},\ol{\nu}_t^i,\ol{\nu}_t^j,\beta_t^{i,p})dt+\sigma_i(t,U_t^{i,p},\ol{\nu}_t^i,\ol{\nu}_t^j)dW_t^{i,p} \nn
\eea
with $t\in[0,T]$, $U_0^{i,p}=\xi^{i,p}$ and $\bg{\beta}^{i,p}\in \mbb{H}^2$ is an 
$A_i$-valued $\mbb{F}$-progressively measurable process.
For the first population, we impose the condition so that $(\xi^{1,p},\bg{\beta}^{1,p},\bg{W}^{1,p})_{1\leq p\leq N_1}$
is exchangeable.  As usual, $\ol{\nu}_t^i:=\frac{1}{N_i}\sum_{p=1}^{N_i}\del_{U_t^{i,p}}$ denotes
the empirical distribution. We shall investigate the following two situations: \\\\
{\rm{\tbf{(setup-1)}}}: The agents in the first population adopt the general exchangeable strategy $(\bg{\beta}^{1,p})_{1\leq p\leq N_1}$
and the agents in the second population  adopt, for any $t\in[0,T]$,
\bea
\beta_t^{2,q}:=\hat{\alpha}_2(t,U_t^{2,q},\mu_t^2,\mu_t^1,u_2(t,U_t^{2,q}))~, 1\leq q\leq N_2.\nn
\eea
{\rm{\tbf{(setup-2)}}}: The first agent in the second population adopts the general strategy $\bg{\beta}^{2,1}$, and 
the remaining agents in the second population as well as the agents in the first population adopt, for any $t\in[0,T]$, 
\bea
&&\beta_t^{1,p}:=\hat{\alpha}_1(t,U_t^{1,p},\mu_t^1,\mu_t^2,u_1(t,U_t^{1,p},\mu_t^1))~, 1\leq p\leq N_1, \nn \\
&&\beta_t^{2,q}:=\hat{\alpha}_2(t,U_t^{2,q},\mu_t^2,\mu_t^1,u_2(t,U_t^{2,q}))~, 2\leq q\leq N_2~. \nn
\eea
\\

The cost functional for any agent in the first population in {\rm{\tbf{(setup-1)}}} is given by
\bea
J_1^{N_1,N_2}(\bg{\beta}^{1,(N_1)},\bg{\hat{\alpha}}^{2,(N_2)}):=
\mbb{E}\Bigl[\int_0^T f_1(t,U_t^{1,p},\ol{\nu}_t^1,\ol{\nu}_t^2,\beta_t^{1,p})dt+g_1(U_T^{1,p},\ol{\nu}_T^1,\ol{\nu}_T^2)\Bigr]~\nn
\eea
with an arbitrary $1\leq p\leq N_1$. On the other hand, 
the cost functional for the first agent in the second population in {\rm{\tbf{(setup-2)}}} is given by
\bea
J_2^{N_2,N_1,1}(\bg{\beta}^{2,1},\bg{\hat{\alpha}}^{2,(N_2)^{-1}},\bg{\hat{\alpha}}^{1,(N_1)})
:=\mbb{E}\Bigl[\int_0^T f_2(t,U_t^{2,1},\ol{\nu}_t^2,\ol{\nu}_t^1,\beta_t^{2,1})dt+g_2(U_t^{2,1},\ol{\nu}_T^2,\ol{\nu}_T^1)\Bigr]~.\nn
\eea
The main result of this section is as follows.

\begin{theorem} Under Assumption {\rm{\tbf{(MFTC-MFG-FA-b)}}} with sufficiently large $(N_i)_{i=1}^2$, 
the feedback control functions $\bigl([0,T]\times \mbb{R}^d\ni (t,x)\mapsto \hat{\alpha}_1(t,x,\mu_t^1,\mu_t^2,u_1(t,x,\mu_t^1))\bigr)$
for the first population and  $\bigl([0,T]\times \mbb{R}^d\ni (t,x)\mapsto \hat{\alpha}_2(t,x,\mu_t^2,\mu_t^1,u_2(t,x))\bigr)$ for
the second one form an $(\sum_{j=1}^2\vep_{N_j})$-approximate Nash equilibrium i.e., 
there exists some constant $C$ independent of $(N_i)_{i=1}^2$ such that
\bea
J_1^{N_1,N_2}(\bg{\beta}^{1,(N_1)},\bg{\hat{\alpha}}^{2,(N_2)})\geq J_1^{N_1,N_2}(\bg{\hat{\alpha}}^{1,(N_1)},\bg{\hat{\alpha}}^{2,(N_2)})
-C\sum_{j=1}^2\vep_{N_j}~\nn
\eea
under  {\rm{\tbf{(setup-1)}}},  and also that
\bea
J_2^{N_2,N_1,1}(\bg{\beta}^{2,1},\bg{\hat{\alpha}}^{2,(N_2)^{-1}},\bg{\hat{\alpha}}^{1,(N_1)})
\geq J_2^{N_2,N_1,1}(\bg{\hat{\alpha}}^{2,(N_2)},\bg{\hat{\alpha}}^{1,(N_1)})-C\sum_{j=1}^2\vep_{N_j} \nn
\eea
under  ${\rm{\tbf{(setup-2)}}}$. In both cases, $\bigl(\vep_{N_j}:=\max(N_j^{-\frac{1}{2}},\ep_{N_j})\bigr)_{1\leq j\leq 2}$.
\begin{proof}
{\rm{\tbf{(first step)}}}:
Let us first prove the claim under  {\rm{\tbf{(setup-1)}}}. 
In contrast to the assumptions used in Theorem~\ref{th-mftc-fa}, the agents in the first population 
now have direct interactions  with those in the second population in their state processes. 
Applying the result of Proposition~\ref{prop-mftc-fa}
to the first population, we get
\bea
&&J_1^{N_1,N_2}(\bg{\beta}^{1,(N_1)},\bg{\hat{\alpha}}^{2,(N_2)})-J_1\geq 
\lambda \mbb{E}\int_0^T |\beta_t^{1,p}-\ul{\hat{\alpha}}_t^{1,p}|^2dt\nn \\
&&\qquad -C\Bigl(1+\sup_{t\in[0,T]}\mbb{E}\bigl[|U_t^{1,p}-\ul{X}_t^{1,p}|^2\bigr]^\frac{1}{2}
+\mbb{E}\bigl[\int_0^T |\beta_t^{1,p}-\ul{\hat{\alpha}}_t^{1,p}|^2 dt\bigr]^\frac{1}{2}\Bigr)\vep_{N_1}\nn \\
&&\qquad+\mbb{E}\Bigl[\int_0^T\bigl( H_1(t,U_t^{1,p},\ol{\nu}_t^1,\ol{\nu}_t^2,\ul{Y}_t^{1,p},\ul{Z}_t^{1,p},\beta_t^{1,p})
-H_1(t,U_t^{1,p},\ol{\nu}_t^1,\mu_t^2, \ul{Y}_t^{1,p},\ul{Z}_t^{1,p},\beta_t^{1,p})\bigr)dt\nn \\
&&\qquad\qquad+g_1(U_T^{1,p},\ol{\nu}_T^1,\ol{\nu}_T^2)-g_1(U_T^{1,p},\ol{\nu}_T^1,\mu_T^2)\Bigr]~,\nn
\eea
where $(\ul{Y}^{1,p},\ul{Z}^{1,p})$ is defined in the same way as $(\ref{bar-bsde})$.
Since $(b_2,\sigma_2)$ are independent from the second measure argument, it is straightforward to confirm 
that $\mbb{E}\bigl[\sup_{t\in[0,T]}|U_t^{2,q}|^2\bigr]\leq C$ for any $q$.
Then we get
$\mbb{E}\bigl[\sup_{t\in[0,T]}|U_t^{1,p}|^2\bigr]\leq C\bigl(1+||\beta^{1,p}||_{\mbb{H}^2}^2\bigr)$ for any $p$.
Moreover, it  is immediate to
obtain, for any $q$,
\bea
\mbb{E}\bigl[\sup_{t\in[0,T]}|U_t^{2,q}-\ul{X}_t^{2,q}|^2\bigr]\leq C\sup_{t\in[0,T]}\mbb{E}\bigl[W_2(\ul{\mu}_t^2,\mu_t^2)^2\bigr]
\leq C\vep_{N_2}^2~.\nn
\eea
In particular, this also implies $\sup_{t\in[0,T]}\mbb{E}\bigl[W_2(\ol{\nu}_t^2,\mu_t^2)^2\bigr]\leq C\vep_{N_2}^2$.
Since, for any $p$, 
\bea
\mbb{E}\bigl[\sup_{s\in[0,t]}|U_s^{1,p}-\ul{X}_s^{1,p}|^2\bigr]
\leq C\mbb{E}\int_0^t \Bigl[|U_s^{1,p}-\ul{X}_s^{1,p}|^2+\sum_{j=1}^2W_2(\ol{\nu}^j_s,\mu_s^j)^2+|\beta_s^{1,p}-\ul{\hat{\alpha}}_s^{1,p}|^2
\Bigr]ds~ \nn
\eea
we get, from the last estimate and the triangle inequality,
\bea
\mbb{E}\bigl[\sup_{t\in[0,T]}|U_t^{1,p}-\ul{X}_t^{1,p}|^2\bigr]\leq C\Bigl(\sum_{j=1}^2\vep_{N_j}^2+||\beta^{1,p}-\ul{\hat{\alpha}}^{1,p}||^2_{\mbb{H}^2}
\Bigr)~. \nn
\eea
By the standard calculation, we see
\bea
&&\Bigl|\mbb{E}\Bigl[\int_0^T\bigl( H_1(t,U_t^{1,p},\ol{\nu}_t^1,\ol{\nu}_t^2,\ul{Y}_t^{1,p},\ul{Z}_t^{1,p},\beta_t^{1,p})
-H_1(t,U_t^{1,p},\ol{\nu}_t^1,\mu_t^2, \ul{Y}_t^{1,p},\ul{Z}_t^{1,p},\beta_t^{1,p})\bigr)dt\nn \\
&&\qquad\qquad+g_1(U_T^{1,p},\ol{\nu}_T^1,\ol{\nu}_T^2)-g_1(U_T^{1,p},\ol{\nu}_T^1,\mu_T^2)\Bigr]\Bigr|\nn \\
&&\leq C\Bigl(1+\sup_{t\in[0,T]}\mbb{E}\Bigl[\sum_{i=1}^2|U_t^{i,p}|^2+|\ul{Y}_t^{1,p}|^2\Bigr]^\frac{1}{2}
+\mbb{E}\Bigl[\int_0^T \bigl[|\ul{Z}_t^{1,p}|^2+|\beta_t^{1,p}|^2\bigr]dt\Bigr]^\frac{1}{2}\Bigr)
\sup_{t\in[0,T]}\mbb{E}\bigl[W_2(\ol{\nu}_t^2,\mu_t^2)^2\bigr]^\frac{1}{2}\nn \\
&&\leq C\bigl(1+||\beta^{1,p}||_{\mbb{H}^2}\bigr)\vep_{N_2}~.\nn
\eea
Since $||\ul{\hat{\alpha}}^{1,p}||_{\mbb{H}^2}\leq C$, we get
\bea
&&J_1^{N_1,N_2}(\bg{\beta}^{1,(N_1)},\bg{\hat{\alpha}}^{2,(N_2)})-J_1\nn \\
&&\qquad \geq 
\lambda||\beta^{1,p}-\ul{\hat{\alpha}}^{1,p}||^2_{\mbb{H}^2}-C\Bigl(1+||\beta_t^{1,p}-\ul{\hat{\alpha}}^{1,p}||_{\mbb{H}^2}\Bigr)
\sum_{j=1}^2 \vep_{N_j}\nn \\
&&\qquad  \geq \Bigl(\lambda-C\sum_{j=1}^2\vep_{N_j}\Bigr)||\beta^{1,p}-\ul{\hat{\alpha}}^{1,p}||^2_{\mbb{H}^2}-C\sum_{j=1}^2\vep_{N_j}~.\nn
\eea
Now Lemma~\ref{lemma-approx-mftc-mfg} gives the desired estimate.
\\\\
{\rm{\tbf{(second step)}}}: Let us now prove the claim under  {\rm{\tbf{(setup-2)}}}.
By putting $i=2$ and $j=1$, all of the arguments in the proof for Theorem~\ref{th-nash-mfg}
work as they are. In fact, due to the independence of $(b_2,\sigma_2)$ from the second measure argument,
some of the estimates become slightly simpler.
In particular, $(\ref{U-norm})$ and $(\ref{olnu-mu-diff})$ hold with $(i=2,j=1)$. The estimate
$(\ref{U-olU})$ is now given by
\bea
\mbb{E}\bigl[\sup_{t\in[0,T]}|U_T^{2,1}-\ol{U}_t^{2,1}|^2\bigr]^\frac{1}{2}\leq C(1+||\beta^{2,1}||_{\mbb{H}^2})\vep_{N_2}~\nn
\eea
without the term $\vep_{N_1}$. The estimate for $|J_2^{N_2,N_1,1}(\bg{\beta}^{2,1},\bg{\hat{\alpha}}^{2,(N_2)^{-1}},\bg{\hat{\alpha}}^{1,(N_1)})
-\ol{J}_2^1(\bg{\beta}^{2,1})|$ is given by exactly the same formula 
as in {\rm{\tbf{(third step)}}} of the proof for Theorem~\ref{th-nash-mfg} with $(i=2,j=1)$.
Now, combining the result in Lemma~\ref{lemma-approx-mftc-mfg}, we get the desired estimate.
\end{proof}
\end{theorem}

\begin{remark}
We have investigated the approximate Nash equilibrium in the closed loop framework.  The analysis for the open loop 
framework can be done in a quite similar manner as explained in Remark~\ref{remark-open}.
Generalization to an arbitrary number of populations $1\leq i\leq m$
can be done straightforwardly, but the direct interactions in the state processes must be carefully arranged.
The empirical distribution of the state of the agents who are in a cooperative population must not appear 
in the coefficients of the state process of the agents in any other populations. The empirical distribution
can appear only in the control strategies, which is indirectly induced by the interactions in the cost functions.
On the other hand, the distribution of the state of the agents who are in a non-cooperative population
can directly appear in the coefficients of the state processes of the agents in any populations. 
\end{remark}

\section{Conclusion and Discussion}
\label{sec-conclusion}
In this work, we have systematically investigated mean field games and mean field type control problems with 
multiple populations for three different situations: (i) every agent is non-cooperative,
(ii) the agents within each population are cooperative, and (iii) the agents in 
some populations are cooperative but not in the other populations.  The relevant adjoint equations were shown 
to be given 
in terms of a coupled system of forward-backward stochastic differential equations of McKean-Vlasov type.
In each case, we have provided several sets of sufficient conditions for the existence of an equilibrium, in particular the one which allows the cost functions of  quadratic growth both in the state variable as well as in its distribution
so that it is applicable to some of the popular setups of linear quadratic problems. 
In the second half of the paper, under additional assumptions, we have proved that each
solution to the mean field problems solved in the first half of the paper actually
provides an approximate Nash equilibrium for the corresponding game
with a large but finite number of agents.

As future works, we may study similar problems by adopting HJB type approach using so-called quadratic growth BSDEs
as in \cite{Carmona-Delarue-1},  where the backward component directly represents the value function.  
Although we need the boundedness of the coefficients and the non-degeneracy for the 
diffusion function, the resultant boundedness of the solution to the BSDEs
will make the analysis simpler.
When each agent is subject to independent 
random Poisson measure, we may use the recent developments of the quadratic growth BSDEs with jumps such as in \cite{Morlais, Kazi, fujii-qg}.
Finally, developing an efficient numerical method for mean field games and 
mean field type control problems remains as a very important issue.
For a general problem, due to 
its infinite dimensionality,  machine learning techniques (such as in \cite{Jentzen}) 
are promising candidates.
If the problem can be approximated by a linear quadratic setup, its solution may help to accelerate the speed of 
convergence for the learning process in the spirit of the work \cite{fujii-deep}.

\begin{appendix}
\section{Proof for Proposition~\ref{prop-mftc-fa}}
\label{app-prop-mftc-fa}
\small
In the following, we use the Landau notation $\calo(\cdot)$ in the sense that $|\calo(x)|\leq C|x|$ with some constant $C$ independent of 
the population sizes $(N_1,N_2)$. Let us define
\bea
&&T_1^{i,p}:=\mbb{E}\bigl[\langle U_T^{i,p}-\ul{X}_T^{i,p}, \ul{Y}_T^{i,p}\rangle\bigr]+
\mbb{E}\Bigl[\int_0^T \bigl[f_i(t,U_t^{i,p},\ol{\nu}_t^i,\ol{\nu}_t^j,\beta_t^{i,p})-f_i(t,\ul{X}_t^{i,p},\mu_t^i,\mu_t^j,\ul{\hat{\alpha}}_t^{i,p})\bigr]dt\Bigr]~,\nn \\
&&T_{2,1}^{i,p}:=\mbb{E}\bigl[g_i(U_T^{i,p},\ol{\nu}_T^i,\ol{\nu}_T^j)-g_i(\ul{X}_T^{i,p},\mu_T^i,\mu_T^j)\bigr]~,\nn \\
&&T_{2,2}^{i,p}:=\mbb{E}\bigl[\langle U_T^{i,p}-\ul{X}_T^{i,p},\part_x g_i(\ul{X}_T^{i,p},\mu_T^i,\mu_T^j)\rangle\bigr]~, \nn \\
&&T_{2,3}^{i,p}:=\mbb{E}\wt{\mbb{E}}\bigl[\langle \wt{U}_T^{i,p}-\wt{\ul{X}}_T^{i,p},\part_\mu g_i(\ul{X}_T^{i,p},\mu_T^i,\mu_T^j)(\wt{\ul{X}}_T^{i,p})\bigr]~, \nn
\eea
which satisfy $J_i^{N_i,N_j}(\bg{\beta}^{i,(N_i)},\bg{\beta}^{j,(N_j)})-J_i=T_1^{i,p}+T_{2}^{i,p}$ 
with $T_2^{i,p}:=T_{2,1}^{i,p}-T_{2,2}^{i,p}-T_{2,3}^{i,p}$.
\subsection{Estimate for $T_2^{i,p}$}
Consider the difference
\bea
&&\Bigl|\mbb{E}\wt{\mbb{E}}\bigl[\langle \wt{U}_T^{i,p}-\wt{\ul{X}}_T^{i,p},\part_\mu g_i(\ul{X}_T^{i,p},\mu_T^i,\mu_T^j)(\wt{\ul{X}}_T^{i,p})\bigr]
-\frac{1}{N_i}\sum_{q=1}^{N_i}\wt{\mbb{E}}\bigl[\langle \wt{U}_T^{i,p}-\wt{\ul{X}}_T^{i,p},
\part_\mu g_i(\wt{\ul{X}}_T^{i,q},\mu_T^i,\mu_T^j)(\wt{\ul{X}}_T^{i,p})\bigr]\Bigr|\nn \\
&&\leq \mbb{E}\bigl[|U_T^{i,p}-\ul{X}_T^{i,p}|^2\bigr]^\frac{1}{2}
\wt{\mbb{E}}\Bigl[\Bigl| \mbb{E}\bigl[\part_\mu g_i(\ul{X}_T^{i,p},\mu_T^i,\mu_T^j)(\wt{\ul{X}}_T^{i,p})\bigr]
-\frac{1}{N_i}\sum_{q=1}^{N_i}\part_\mu g_i(\wt{\ul{X}}_T^{i,q},\mu_T^i,\mu_T^j)(\wt{\ul{X}}_T^{i,p})\Bigr|^2\Bigr]^\frac{1}{2}\nn \\
&&\leq C\mbb{E}\bigl[|U_T^{i,p}-\ul{X}_T^{i,p}|^2\bigr]^\frac{1}{2}N_i^{-\frac{1}{2}}. 
\label{app-weak-LofL}
\eea
The last estimate is from the law of large numbers with the finite second moment of $\part_\mu g_i$-term 
and the independence of $(\wt{\ul{X}}_T^{i,q})_{1\leq q\leq N_i}$.
Taking the average in $p$, we get,
\bea
\frac{1}{N_i}\sum_{p=1}^{N_i}T_{2,3}^{i,p}&=&\mbb{E}\bigl[|U_T^{i,1}-\ul{X}_T^{i,1}|^2\bigr]^\frac{1}{2}\calo(N_i^{-\frac{1}{2}})
+\frac{1}{N_i^2}\sum_{p,q=1}^{N_i}\mbb{E}\bigl[\langle U_T^{i,p}-\ul{X}_T^{i,p},\part_\mu g_i(\ul{X}_T^{i,q},\mu_T^i,\mu_T^j)(\ul{X}_T^{i,p})\bigr]~ \nn \\
&=&\mbb{E}\bigl[|U_T^{i,1}-\ul{X}_T^{i,1}|^2\bigr]^\frac{1}{2}\calo(N_i^{-\frac{1}{2}})+\frac{1}{N_i}\sum_{p=1}^{N_i}
\mbb{E}\wt{\mbb{E}}\bigl[\langle U_T^{i,\theta}-\ul{X}_T^{i,\theta},\part_\mu g_i(\ul{X}_T^{i,p},\mu_T^i,\mu_T^j)(\ul{X}_T^{i,\theta})\bigr]~, \nn
\eea
where $\theta$ is a random variable on $(\wt{\Omega},\wt{\calf},\wt{\mbb{P}})$ with uniform distribution on 
the set $\{1,\cdots,N_i\}$. Using Lemma~\ref{lemma-cl-mftc} and the Lipschitz property of $\part_\mu g_i$ with respect to the 
first measure argument, we get
\bea
\frac{1}{N_i}\sum_{p=1}^{N_i}T_{2,3}^{i,p}=\mbb{E}\bigl[|U_T^{i,1}-\ul{X}_T^{i,1}|^2\bigr]^\frac{1}{2}\calo(\vep_{N_i})+\frac{1}{N_i}\sum_{p=1}^{N_i}
\mbb{E}\wt{\mbb{E}}\bigl[\langle U_T^{i,\theta}-\ul{X}_T^{i,\theta},\part_\mu g_i(\ul{X}_T^{i,p},\ul{\mu}_T^i,\mu_T^j)(\ul{X}_T^{i,\theta})\bigr]~. 
\label{app-T23}
\eea
Using Lemma~\ref{lemma-cl-mftc} also for $T_{2,2}^{i,p}, T_{2,3}^{i,p}$, we obtain
\bea
\frac{1}{N_i}\sum_{p=1}^{N_i}T_2^{i,p}&=&\frac{1}{N_i}\sum_{p=1}^{N_i}\Bigl\{\mbb{E}\bigl[g_i(U_T^{i,p},\ol{\nu}_T^i,\mu_T^j)-g_i(\ul{X}_T^{i,p},\ul{\mu}_T^i,\mu_T^j)\bigr]\nn \\
&&-\mbb{E}\bigl[\langle U_T^{i,p}-\ul{X}_T^{i,p},\part_x g_i(\ul{X}_T^{i,p},\ul{\mu}_T^i,\mu_T^j)\rangle\bigr]
-\mbb{E}\wt{\mbb{E}}\bigl[\langle U_T^{i,\theta}-\ul{X}_T^{i,\theta},\part_\mu g_i(\ul{X}_T^{i,p},\ul{\mu}_T^i,\mu_T^j)(
\ul{X}_T^{i,\theta})\rangle\bigr]\Bigr\}\nn \\
&&+\Bigl(1+\mbb{E}\bigl[|U_T^{i,1}-\ul{X}_T^{i,1}|^2\bigr]^\frac{1}{2}\Bigr)\calo(\vep_{N_i})
+\frac{1}{N_i}\sum_{p=1}^{N_i}\mbb{E}\bigl[g_i(U_T^{i,p},\ol{\nu}_T^i,\ol{\nu}_T^j)-g_i(U_T^{i,p},\ol{\nu}_T^i,\mu_T^j)\bigr]~.\nn
\eea
Using the fact that the conditional law of $U_T^{i,\theta}$ (respectively $\ul{X}_T^{i,\theta}$) under $\wt{\mbb{P}}$ 
is given by the empirical distribution $\ol{\nu}_T^i$ (respectively $\ul{\mu}_T^i$), the convexity in {\rm{\tbf{(MFTC-a) (A6)}}}
implies
\bea
\frac{1}{N_i}\sum_{p=1}^{N_i}T_2^{i,p}&\geq &\Bigl(1+\mbb{E}\bigl[|U_T^{i,1}-\ul{X}_T^{i,1}|^2\bigr]^\frac{1}{2}\Bigr)\calo(\vep_{N_i})
+\frac{1}{N_i}\sum_{p=1}^{N_i}\mbb{E}\bigl[g_i(U_T^{i,p},\ol{\nu}_T^i,\ol{\nu}_T^j)-g_i(U_T^{i,p},\ol{\nu}_T^i,\mu_T^j)\bigr]~.\nn
\eea 

\subsection{Estimate for $T_1^{i,p}$}
Using Ito formula to evaluate $\mbb{E}\bigl[\langle U_T^{i,p}-\ul{X}_T^{i,p}, \ul{Y}_T^{i,p}\rangle \bigr]$, we can rewrite $T_1^{i,p}$ as 
\bea
T_1^{i,p}&=&\mbb{E}\int_0^T \Bigl\{ H_i(t,U_t^{i,p},\ol{\nu}_t^i,\ol{\nu}_t^j,\ul{Y}_t^{i,p},\ul{Z}_t^{i,p},\beta_t^{i,p})-H_i(t,\ul{X}_t^{i,p},\mu_t^i,\mu_t^j,
\ul{Y}_t^{i,p},\ul{Z}_t^{i,p},\ul{\hat{\alpha}}_t^{i,p})\nn \\
&&\qquad-\langle U_t^{i,p}-\ul{X}_t^{i,p}, \part_x H_i(t,\ul{X}_t^{i,p},\mu_t^i,\mu_t^j,\ul{Y}_t^{i,p},\ul{Z}_t^{i,p},\ul{\hat{\alpha}}_t^{i,p})\rangle \nn \\
&&\qquad-\wt{\mbb{E}}\bigl[\langle \wt{U}_t^{i,p}-\wt{\ul{X}}_t^{i,p}, \part_\mu H_i (t,\ul{X}_t^{i,p},\mu_t^i,\mu_t^j,
\ul{Y}_t^{i,p},\ul{Z}_t^{i,p},\ul{\hat{\alpha}}_t^{i,p})(\wt{\ul{X}}_t^{i,p})\rangle\bigr]\Bigr\}dt\nn \\
&=:&T_{1,1}^{i,p}-T_{1,2}^{i,p}-T_{1,3}^{i,p}~.\nn
\eea 
Using local Lipschitz continuity of $H_i$ with respect to the first measure argument, the estimate 
\be
\mbb{E}\Bigl[\sup_{t\in[0,T]}|\ul{X}_t^{i,p}|^2+\sup_{t\in[0,T]}|\ul{Y}_t^{i,p}|^2+\int_0^T(|\ul{Z}_t^{i,p}|^2+|\ul{\hat{\alpha}}_t^{i,p}|^2) dt\Bigr]\leq C~, \nn
\ee
and the result in Lemma~\ref{lemma-cl-mftc}, we get
\bea
T_{1,1}^{i,p}=\mbb{E}\int_0^T \Bigl\{H_i(t,U_t^{i,p},\ol{\nu}_t^i,\ol{\nu}_t^j,\ul{Y}_t^{i,p},\ul{Z}_t^{i,p},\beta_t^{i,p})
-H_i(t,\ul{X}_t^{i,p},\ul{\mu}_t^i,\mu_t^j,\ul{Y}_t^{i,p},\ul{Z}_t^{i,p},\ul{\hat{\alpha}}_t^{i,p})\Bigr\}dt+\calo(\vep_{N_i})~.\nn
\eea
We also get by similar calculation that
\bea
T_{1,2}^{i,p}&=&\mbb{E}\int_0^T \bigl[\langle U_t^{i,p}-\ul{X}_t^{i,p}, \part_x H_i(t,\ul{X}_t^{i,p},\ul{\mu}_t^i,\mu_t^j,\ul{Y}_t^{i,p},\ul{Z}_t^{i,p},\ul{\hat{\alpha}}_t^{i,p}) \rangle \bigr]dt \nn \\
&&+\sup_{t\in[0,T]}\mbb{E}\bigl[|U_t^{i,p}-\ul{X}_t^{i,p}|^2\bigr]^\frac{1}{2}\calo(\vep_{N_i})~.\nn
\eea
By the same arguments used in $(\ref{app-weak-LofL})$, we get
\bea
T_{1,3}^{i,p}&=&\frac{1}{N_i}\sum_{q=1}^{N_i}\int_0^T \mbb{E}\bigl[
\langle U_t^{i,p}-\ul{X}_t^{i,p}, \part_\mu H_i (t,\ul{X}_t^{i,q},\mu_t^i,\mu_t^j, \ul{Y}_t^{i,q},\ul{Z}_t^{i,q},\ul{\hat{\alpha}}_t^{i,q})(\ul{X}_t^{i,p})
\rangle \bigr]dt\nn \\
&&+\sup_{t\in[0,T]}\mbb{E}\bigl[|U_t^{i,1}-\ul{X}_t^{i,1}|^2\bigr]^\frac{1}{2}\calo(N_i^{-\frac{1}{2}})~, \nn
\eea
and then same analysis used for $(\ref{app-T23})$ gives
\bea
\frac{1}{N_i}\sum_{p=1}^{N_i} T_{1,3}^{i,p}&=&\frac{1}{N_i}\sum_{p=1}^{N_i}\int_0^T
\mbb{E}\wt{\mbb{E}}\bigl[\langle U_t^{i,\theta}-\ul{X}_t^{i,\theta}, \part_\mu H_i(t,\ul{X}_t^{i,p},\ul{\mu}_t^i,\mu_t^j,
\ul{Y}_t^{i,p},\ul{Z}_t^{i,p},\ul{\hat{\alpha}}_t^{i,p})(\ul{X}_t^{i,\theta})\rangle \bigr]dt \nn \\
&&+\sup_{t\in[0,T]}\mbb{E}\bigl[|U_t^{i,1}-\ul{X}_t^{i,1}|^2\bigr]^\frac{1}{2}\calo(\vep_{N_i})~. \nn
\eea
Finally, we can see that
\bea
&&\Bigl|\mbb{E}\int_0^T\bigl[\langle  \beta_t^{i,p}-\ul{\hat{\alpha}}_t^{i,p},
\part_\alpha H_i(t,\ul{X}_t^{i,p},\ul{\mu}_t^i,\mu_t^j,\ul{Y}_t^{i,p},\ul{Z}_t^{i,p},\ul{\hat{\alpha}}_t^{i,p})
-\part_\alpha H_i(t,\ul{X}_t^{i,p},\mu_t^i,\mu_t^j,\ul{Y}_t^{i,p},\ul{Z}_t^{i,p},\ul{\hat{\alpha}}_t^{i,p})\rangle \bigr]dt \Bigr|\nn \\
&&\leq C\mbb{E}\Bigl[\int_0^T |\beta_t^{i,p}-\ul{\hat{\alpha}}_t^{i,p}|^2dt\Bigr]^\frac{1}{2}\vep_{N_i}~.\nn
\eea
Using the optimality condition, exchangeability,  and the results obtained above, we get
\bea
&&\frac{1}{N_i}\sum_{p=1}^{N_i}T_1^{i,p}\geq \frac{1}{N_i}\sum_{p=1}^{N_i}(T_{1,1}^{i,p}-T_{1,2}^{i,p}-T_{1,3}^{i,p})-\frac{1}{N_i}\sum_{p=1}^{N_i}\mbb{E}\int_0^T \bigl[\langle \beta_t^{i,p}-\ul{\hat{\alpha}}_t^{i,p},
\part_\alpha H_i (t,\ul{X}_t^{i,p},\mu_t^i,\mu_t^j,\ul{Y}_t^{i,p},\ul{Z}_t^{i,p},\ul{\hat{\alpha}}_t^{i,p})\rangle \bigr]dt\nn \\
&&=\frac{1}{N_i}\sum_{p=1}^{N_i}\mbb{E}\int_0^T \Bigl\{
H_i(t,U_t^{i,p},\ol{\nu}_t^i,\mu_t^j,\ul{Y}_t^{i,p},\ul{Z}_t^{i,p},\beta_t^{i,p})-H_i(t,\ul{X}_t^{i,p},\ul{\mu}_t^i,\mu_t^j,
\ul{Y}_t^{i,p},\ul{Z}_t^{i,p},\ul{\hat{\alpha}}_t^{i,p})\nn \\
&&\quad-\langle U_t^{i,p}-\ul{X}_t^{i,p}, \part_x H_i(t,\ul{X}_t^{i,p},\ul{\mu}_t^i,\mu_t^j,\ul{Y}_t^{i,p},\ul{Z}_t^{i,p},\ul{\hat{\alpha}}_t^{i,p})\rangle 
-\wt{\mbb{E}}\bigl[\langle U_t^{i,\theta}-\ul{X}_t^{i,\theta}, \part_\mu H_i(t,\ul{X}_t^{i,p},\ul{\mu}_t^i,\mu_t^j,\ul{Y}_t^{i,p},
\ul{Z}_t^{i,p},\ul{\hat{\alpha}}_t^{i,p})(\ul{X}_t^{i,\theta})\rangle \bigr]\nn \\
&&\quad -\langle \beta_t^{i,p}-\ul{\hat{\alpha}}_t^{i,p},\part_\alpha H_i(t,\ul{X}_t^{i,p},\ul{\mu}_t^i,\mu_t^j,\ul{Y}_t^{i,p},\ul{Z}_t^{i,p},
\ul{\hat{\alpha}}_t^{i,p})\rangle \Bigr\}dt\nn \\
&&\quad+\Bigl(1+\sup_{t\in[0,T]}\mbb{E}\bigl[|U_t^{i,1}-\ul{X}_t^{i,1}|^2\bigr]^\frac{1}{2}+\mbb{E}\Bigl[\int_0^T
|\beta_t^{i,1}-\ul{\hat{\alpha}}_t^{i,1}|^2dt\Bigr]^\frac{1}{2}\Bigr)\calo(\vep_{N_i})\nn \\
&&\quad+\mbb{E}\int_0^T\bigl[ H_i(t,U_t^{i,1},\ol{\nu}_t^i,\ol{\nu}_t^j,\ul{Y}_t^{i,1},\ul{Z}_t^{i,1},\beta_t^{i,1})
-H_i(t,U_t^{i,1},\ol{\nu}_t^i,\mu_t^j,\ul{Y}_t^{i,1},\ul{Z}_t^{i,1},\beta_t^{i,1})\bigr]dt\nn \\
&&\geq \lambda \mbb{E}\int_0^T |\beta_t^{i,1}-\ul{\hat{\alpha}}_t^{i,1}|^2 dt\nn \\
&&\quad+\Bigl(1+\sup_{t\in[0,T]}\mbb{E}\bigl[|U_t^{i,1}-\ul{X}_t^{i,1}|^2\bigr]^\frac{1}{2}+\mbb{E}\Bigl[\int_0^T
|\beta_t^{i,1}-\ul{\hat{\alpha}}_t^{i,1}|^2dt\Bigr]^\frac{1}{2}\Bigr)\calo(\vep_{N_i})\nn \\
&&\quad+\mbb{E}\int_0^T\bigl[ H_i(t,U_t^{i,1},\ol{\nu}_t^i,\ol{\nu}_t^j,\ul{Y}_t^{i,1},\ul{Z}_t^{i,1},\beta_t^{i,1})
-H_i(t,U_t^{i,1},\ol{\nu}_t^i,\mu_t^j,\ul{Y}_t^{i,1},\ul{Z}_t^{i,1},\beta_t^{i,1})\bigr]dt~.\nn 
\eea
\subsection{Final Step}
By exchangeability,  we have
\bea
&&J_i^{N_i,N_j}(\bg{\beta}^{i,(N_i)},\bg{\beta}^{j,(N_j)})-J_i=T_1^{i,p}+T_{2}^{i,p} =\frac{1}{N_i}\sum_{p=1}^{N_i}(T_1^{i,p}+T_2^{i,p})\nn \\
&&\quad \geq  \lambda \mbb{E}\int_0^T |\beta_t^{i,1}-\ul{\hat{\alpha}}_t^{i,1}|^2 dt
+\Bigl(1+\sup_{t\in[0,T]}\mbb{E}\bigl[|U_t^{i,1}-\ul{X}_t^{i,1}|^2\bigr]^\frac{1}{2}+\mbb{E}\Bigl[\int_0^T
|\beta_t^{i,1}-\ul{\hat{\alpha}}_t^{i,1}|^2dt\Bigr]^\frac{1}{2}\Bigr)\calo(\vep_{N_i})\nn \\
&&\quad+\mbb{E}\int_0^T\bigl[ H_i(t,U_t^{i,1},\ol{\nu}_t^i,\ol{\nu}_t^j,\ul{Y}_t^{i,1},\ul{Z}_t^{i,1},\beta_t^{i,1})
-H_i(t,U_t^{i,1},\ol{\nu}_t^i,\mu_t^j,\ul{Y}_t^{i,1},\ul{Z}_t^{i,1},\beta_t^{i,1})\bigr]dt~ \nn \\
&&\quad+\mbb{E}\bigl[g_i(U_T^{i,1},\ol{\nu}_T^i,\ol{\nu}_T^j)-g_i(U_T^{i,1},\ol{\nu}_T^i,\mu_T^j)\bigr]~.\nn 
\eea
This gives the desired result.
\end{appendix}



\end{document}